\pageno=1
\input pictex
\input xy
\xyoption{all}

\footline={}
\headline={\ifnum\pageno=1{\hfil}\else
 {\ifodd\pageno{{\tensl An introduction to algebraic surgery}\hfil{\tenrm\folio}}
  \else{{\tenrm\folio}\hfil{\tensl Andrew Ranicki}}
\fi}\fi}
\hoffset=1in
\voffset=1in
\hsize=27pc
\vsize=520pt
\baselineskip=13pt
\hfuzz=15pt
\parskip=1pc plus 1pt minus 1pt
\parindent=10pt
\baselineskip=12truept 
\displaywidth=\hsize
\displayindent=0pt
\abovedisplayskip=4.5pt plus 1pt minus 1pt
\belowdisplayskip=4.5pt plus 1pt minus 1pt
\abovedisplayshortskip=4.5pt plus 1pt minus 1pt
\belowdisplayshortskip=4.5pt plus 1pt minus 1pt

\font\titlefont=cmr17
\font\authorfont=cmr12
\font\sectionfont=cmbx12

\font\tenrm=cmr10 at 10pt

\font\sc=cmcsc10

\font\it=cmti10

\font\tenmsa=msam10
\font\sevenmsa=msam7
\font\fivemsa=msam5
\newfam\msafam
\textfont\msafam=\tenmsa \scriptfont\msafam=\sevenmsa
 \scriptscriptfont\msafam=\fivemsa

\font\tenmsb=msbm10
\font\sevenmsb=msbm7
\font\fivemsb=msbm5
\newfam\msbfam
\textfont\msbfam=\tenmsb \scriptfont\msbfam=\sevenmsb
 \scriptscriptfont\msbfam=\fivemsb

\def\hexnumber@#1{\ifnum#1<10 \number#1\else
 \ifnum#1=10 A\else\ifnum#1=11 B\else\ifnum#1=12 C\else
 \ifnum#1=13 D\else\ifnum#1=14 E\else\ifnum#1=15 F\fi\fi\fi\fi\fi\fi\fi}

\def\msa@{\hexnumber@\msafam}
\def\msb@{\hexnumber@\msbfam}
\mathchardef\boxtimes="2\msa@02
\mathchardef\subsetneqq="3\msb@24
\mathchardef\looparrowright="3\msa@23

\def\ccr{\cr\noalign{\vskip2pt}}

\def\dirlim#1{\mathop{\vtop{\ialign{##\crcr
  $\hfil\displaystyle{\rm lim}\hfil$\crcr\noalign{\kern0pt\nointerlineskip}
  \rightarrowfill\crcr\noalign{\kern0pt\nointerlineskip}
  $\hfil\scriptstyle{#1}\hfil$\crcr
  }}}}
\def\hocolim#1{\mathop{\vtop{\ialign{##\crcr
  $\hfil\displaystyle{\rm hocolim}\hfil$\crcr\noalign{\kern0pt\nointerlineskip}
  \rightarrowfill\crcr\noalign{\kern0pt\nointerlineskip}
  $\hfil\scriptstyle{#1}\hfil$\crcr
  }}}}
\def\invlim#1{\mathop{\vtop{\ialign{##\crcr
  $\hfil\displaystyle{\rm lim}\hfil$\crcr\noalign{\kern0pt\nointerlineskip}
  \leftarrowfill\crcr\noalign{\kern0pt\nointerlineskip}
  $\hfil\scriptstyle{#1}\hfil$\crcr
  }}}}
\def\holim#1{\mathop{\vtop{\ialign{##\crcr
  $\hfil\displaystyle{\rm holim}\hfil$\crcr\noalign{\kern0pt\nointerlineskip}
  \leftarrowfill\crcr\noalign{\kern0pt\nointerlineskip}
  $\hfil\scriptstyle{#1}\hfil$\crcr
  }}}}

\font\tenmsb=msbm10
\font\sevenmsb=msbm7
\font\fivemsb=msbm5
\newfam\msbfam
\def\msb{\fam\msbfam\tenmsb}
\textfont\msbfam=\tenmsb
\scriptfont\msbfam=\sevenmsb
\scriptscriptfont\msbfam=\fivemsb

\font\tenscr=eusm10
\font\sevenscr=eusm7
\font\fivescr=eusm5
\newfam\scrfam
\def\scr{\fam\scrfam\tenscr}
\textfont\scrfam=\tenscr
\scriptfont\scrfam=\sevenscr
\scriptscriptfont\scrfam=\fivescr

\def\bbc{{\msb C}}

\def\bbl{{\msb L}}

\def\bbp{{\msb P}}
\def\bbq{{\msb Q}}
\def\bbr{{\msb R}}

\def\bbz{{\msb Z}}

\def\la#1{\hbox to #1pc{\leftarrowfill}}
\def\ra#1{\hbox to #1pc{\rightarrowfill}}

\def\fract#1#2{\raise4pt\hbox{$ #1 \atop #2 $}}

\def\hsmash{\triangleright\kern-4.5pt\raise.66pt\hbox{$\scriptstyle{<}$}}

\def\sqr#1#2{{\vcenter{\vbox{\hrule height.#2pt
  \hbox{\vrule width.#2pt height#1pt \kern#1pt 
   \vrule width.#2pt}
  \hrule height.#2pt}}}}

\def\square{\mathchoice\sqr44\sqr44\sqr{2.9}3\sqr{1.95}3}
\def\endstate{\vskip1pt \hbox to \hsize{\hfil $\square$}\smallskip}

\def\mapto{\mapstochar\!\ra{1}~}
\def\isa{\raise5pt\hbox{$\lower2pt\hbox{$\cong$}\atop\to$}}

\def\underleftarrow#1{\vbox{\ialign{##\crcr
  $\hfil\displaystyle{#1}\hfil$\crcr
  \noalign{\kern-1pt\nointerlineskip}
  \leftarrowfill\crcr}}}
\def\underrightarrow#1{\vbox{\ialign{##\crcr
  $\hfil\displaystyle{#1}\hfil$\crcr
  \noalign{\kern-1pt\nointerlineskip}
  \rightarrowfill\crcr}}}
\def\t{\thinspace}
\def\bysame{$\hbox to 25pt{\hrulefill}$~}

\tolerance=500
\overfullrule=0pt
\null\vskip 0.5 in
\centerline{\titlefont An introduction to algebraic surgery}
\bigskip
\centerline{\authorfont Andrew Ranicki}
\vskip 0.5 in
\noindent{\sectionfont Introduction}

Surgery theory investigates the homotopy types of manifolds, using a
combination of algebra and topology. It is the aim of these notes to
provide an introduction to the more algebraic aspects of the theory,
without losing sight of the geometric motivation.

\noindent {\bf 0.1 Historical background}

A closed $m$-dimensional topological manifold $M$ has Poincar\'e duality
isomorphisms
$$H^{m-*}(M)~\cong~H_*(M)~.$$
In order for a space $X$ to be homotopy equivalent to an $m$-dimensional
manifold it is thus necessary (but not in general sufficient) for $X$ to 
be an $m$-dimensional Poincar\'e duality space, with
$H^{m-*}(X)\cong H_*(X)$. The {\it topological structure set}\ ${\scr S}^{TOP}(X)$ 
is defined to be the set of equivalence classes of pairs
$$(\hbox{\rm $m$-dimensional manifold}~M~,~\hbox{\rm homotopy equivalence}~h:M \to X)$$
subject to the equivalence relation 
$$\eqalign{(M,h) \sim (M',h')~~&\hbox{\rm if there exists a homeomorphism}\ccr
&\hbox{\rm $f:M \to M'$ such that $h'f \simeq h:M \to X$}.}$$
The basic problem of surgery theory is to decide if a Poincar\'e
complex $X$ is homotopy equivalent to a manifold (i.e. if ${\scr S}^{TOP}(X)$ 
is non-empty), and if so to compute ${\scr S}^{TOP}(X)$ in terms of the 
algebraic topology of $X$.

Surgery theory was first developed for differentiable manifolds, and
then extended to $PL$ and topological manifolds. 

The classic Browder--Novikov--Sullivan--Wall obstruction theory
for deciding if a Poincar\'e complex $X$ is homotopy equivalent to
a manifold has two stages :
\smallskip
{\parindent=23pt
\parskip=2pt
\item{\rm (i)} the primary topological $K$-theory obstruction $\nu_X
\in [X,B(G/TOP)]$ to a $TOP$ reduction $\widetilde{\nu}_X:X \to BTOP$
of the Spivak normal fibration $\nu_X:X \to BG$, which vanishes if and
only if there exists a manifold $M$ with a normal map $(f,b):M \to X$,
that is a degree 1 map $f:M \to X$ with a bundle map $b:\nu_M \to
\widetilde{\nu}_X$,
\vskip2pt
\item{\rm (ii)} a secondary algebraic $L$-theory obstruction
$$\sigma_*(f,b) \in L_m(\bbz[\pi_1(X)])$$
in the surgery obstruction group of Wall [28], which is defined if the
obstruction in (i) vanishes, and which depends on the choice of 
$TOP$ reduction $\widetilde{\nu}_X$, or equivalently on the bordism class
of the normal map $(f,b):M \to X$.
The surgery obstruction is such that $\sigma_*(f,b)=0$
if (and for $m\geq 5$ only if) $(f,b)$ is normal bordant to a homotopy equivalence.\par}
\noindent There exists a $TOP$ reduction $\widetilde{\nu}_X$ of $\nu_X$
for which the corresponding normal map $(f,b):M \to X$ has zero surgery
obstruction if (and for $m \geq 5$ only if) the structure set 
${\scr S}^{TOP}(X)$ is non-empty. A relative version of the theory gives a
two-stage obstruction for deciding if a homotopy equivalence $M \to X$
from a manifold $M$ is homotopic to a homeomorphism, which is traditionally
formulated as the {\it surgery exact sequence}
$$\dots \to L_{m+1}(\bbz[\pi_1(X)]) \to {\scr S}^{TOP}(X) \to [X,G/TOP]
\to L_m(\bbz[\pi_1(X)])~.$$
See the paper by Browder [2] for an account of the 
original Sullivan-Wall surgery exact sequence in the differentiable category
in the case when $X$ has the homotopy type of a differentiable manifold
$$\dots \to L_{m+1}(\bbz[\pi_1(X)]) \to {\scr S}^O(X) \to [X,G/O] \to
L_m(\bbz[\pi_1(X)])~.$$
\indent 
The algebraic $L$-groups $L_*(\Lambda)$ of a ring with involution 
$\Lambda$ are defined
using quadratic forms over $\Lambda$ and their automorphisms, and are 4-periodic
$$L_m(\Lambda)~=~L_{m+4}(\Lambda)~.$$ 
The surgery classification of exotic spheres of Kervaire and
Milnor [7] included the first computation of the $L$-groups, namely
$$L_m(\bbz)~=~\cases{\bbz&if $m \equiv 0$ $(\bmod\,4)$\ccr
0&if $m \equiv 1$ $(\bmod\,4)$\ccr
\bbz_2&if $m \equiv 2$ $(\bmod\,4)$\ccr
0&if $m \equiv 3$ $(\bmod\,4)$~.}$$
\indent
The relationship between topological and $PL$ manifolds was
investigated using surgery methods in the 1960's by Novikov, Casson,
Sullivan, Kirby and Siebenmann [8] (cf. Ranicki [23]), culminating in
a disproof of the manifold Hauptvermutung\t : there exist
homeomorphisms of $PL$ manifolds which are not homotopic to $PL$
homeomorphisms, and in fact there exist topological manifolds without
$PL$ structure. The surgery exact sequence for the $PL$ manifold
structure set ${\scr S}^{PL}(M)$ for a $PL$ manifold $M$ was related to
the exact sequence for ${\scr S}^{TOP}(M)$ by a commutative braid of
exact sequences
$$\hskip5pt\xymatrix@!C@C-48pt@R-10pt{
H^3(M;\bbz_2)
\ar[dr]\ar@/^2pc/[rr]^{}&&
[M,G/PL]\ar[dr] \ar@/^2pc/[rr]&&L_m(\bbz[\pi_1(M)]) \\&
{\scr S}^{PL}(M)\ar[ur] \ar[dr] && [M,G/TOP]
\ar[ur] \ar[dr]&&\\
L_{m+1}(\bbz[\pi_1(M)])
\ar[ur]\ar@/_2pc/[rr]&&{\scr S}^{TOP}(M)
\ar[ur]\ar@/_2pc/[rr]_{}&&H^4(M;\bbz_2)
}$$
with
$$\pi_*(G/TOP)~=~L_*(\bbz)~.$$
\indent Quinn [17] gave a geometric construction of a 
spectrum of simplicial sets for any group $\pi$ 
$$\bbl_{\bullet}(\bbz[\pi])~=~\{\bbl_n(\bbz[\pi])\,\vert\,\Omega\bbl_n(\bbz[\pi]) \simeq \bbl_{n+1}(\bbz[\pi])\}$$
with homotopy groups
$$\pi_n(\bbl_{\bullet}(\bbz[\pi]))~=~\pi_{n+k}(\bbl_{-k}(\bbz[\pi]))~=~
L_n(\bbz[\pi])~,$$
and 
$$\bbl_0(\bbz)~\simeq~L_0(\bbz) \times G/TOP~.$$
The construction included an assembly map
$$A~:~H_*(X;\bbl_{\bullet}(\bbz)) \to L_*(\bbz[\pi_1(X)])$$
and for a manifold $X$ the surgery obstruction function is given by
$$[X,G/TOP]\subset [X,L_0(\bbz)\times G/TOP]~\cong~H_n(X;\bbl_{\bullet}(\bbz))~
\raise4pt\hbox{$\displaystyle{A} \atop\to$}~L_m(\bbz[\pi_1(X)])~.$$
The surgery classifying spectra $\bbl_{\bullet}(\Lambda)$ 
and the assembly map $A$ were constructed algebraically
in Ranicki [22] for any ring with involution $\Lambda$,
using quadratic Poincar\'e complex $n$-ads over $\Lambda$.
The spectrum $\bbl_{\bullet}(\bbz)$ is appropriate for
the surgery classification of homology manifold structures (Bryant,
Ferry, Mio and Weinberger [3]); for topological manifolds it is necessary
to work with the 1-connective spectrum
$\bbl_{\bullet}=\bbl_{\bullet}(\bbz)\langle 1 \rangle$,
such that $\bbl_n$ is $n$-connected with $\bbl_0 \simeq G/TOP$. 
The relative homotopy groups of the spectrum-level assembly map
$${\scr S}_m(X)~=~\pi_m(A:X_+\wedge \bbl_{\bullet} \to \bbl_{\bullet}(\bbz[\pi_1(X)]))$$
fit into the {\it algebraic surgery exact sequence}
$$\eqalign{\dots \to L_{m+1}(\bbz[\pi_1(X)])&\to {\scr S}_{m+1}(X)\cr
&\to H_m(X;\bbl_{\bullet})~
\raise4pt\hbox{$\displaystyle{A} \atop\ra{1.2}$}~L_m(\bbz[\pi_1(X)]) \to \dots~.}$$

The algebraic surgery theory of Ranicki [20],\t [22] provided one-stage
obstructions\t :
{\parindent=23pt
\parskip=2pt
\item{\rm (i)}
An $m$-dimensional Poincar\'e duality space $X$ has a {\it total surgery 
obstruction} $s(X) \in {\scr S}_m(X)$ such that $s(X)=0$ if (and for $m \geq 5$
only if) $X$ is homotopy equivalent to a manifold.
\item{\rm (ii)}
A homotopy equivalence of $m$-dimensional manifolds $h:M' \to M$ has a 
{\it total surgery obstruction}
$s(h) \in {\scr S}_{m+1}(M)$ such that $s(h)=0$ if (and for $m \geq 5$
only if) $h$ is homotopic to a homeomorphism.\par}
\noindent Moreover, if $X$ is an $m$-dimensional manifold and $m \geq 5$
the geometric surgery exact sequence is isomorphic to the algebraic
surgery exact sequence
$$\xymatrix@C-11pt@R+5pt
{\dots \ar[r] & L_{m+1}(\bbz[\pi_1(X)]) \ar@{=}[d] \ar[r]&
{\scr S}^{TOP}(X) \ar[d]^{\displaystyle{\cong}} \ar[r] 
& [X,G/TOP] \ar[d]^{\displaystyle{\cong}} \ar[r] &
L_m(\bbz[\pi_1(X)]) \ar@{=}[d] \\
\dots \ar[r] & L_{m+1}(\bbz[\pi_1(X)]) \ar[r]&
{\scr S}_{m+1}(X) \ar[r] & H_m(X;\bbl_{\bullet}) \ar[r]^-{\displaystyle{A}} &
L_m(\bbz[\pi_1(X)])}$$
with
$${\scr S}^{TOP}(X) \to {\scr S}_{m+1}(X) ~;~(M,h:M \to X) \mapto s(h)~.$$

Given a normal map $(f,b):M^m \to X$ it is possible to {\it kill} an element
$x \in \pi_r(f)$ by surgery if and only if $x$ can be represented by an
embedding $S^{r-1} \times D^{n-r+1} \hookrightarrow M$ with
a null-homotopy in $X$, 
in which case the {\it trace} of the surgery is a normal bordism
$$((g,c);(f,b),(f',b'))~:~(N;M,M') \to X \times ([ 0,1];\{0\},\{1\})$$
with
$$\eqalign{&N^{m+1}~=~M \times I \cup D^r \times D^{m-r+1}~,\ccr
&M^{\prime m}~=~(M \backslash S^{r-1} \times D^{m-r+1})\cup D^r \times S^{m-r}~.}$$
The normal map $(f',b'):M'\to X$ is the {\it geometric effect} of the surgery on $(f,b)$. 
Surgery theory investigates the extent to which a normal map can be made
bordant to a homotopy equivalence by killing as much of $\pi_*(f)$
as possible. The original definition of the non-simply-connected surgery obstruction
$\sigma_*(f,b) \in L_m(\bbz[\pi_1(X)])$ (Wall [28]) was obtained
after preliminary surgeries below the middle dimension, to kill the relative
homotopy groups $\pi_r(f)$ for $2r \leq m$. It could thus be assumed
that $(f,b):M \to X$ is $[m/2]$-connected, with $\pi_r(f)=0$ for $2r
\leq m$, and $\sigma_*(f,b)$ was defined using the Poincar\'e duality
structure on the middle-dimensional homotopy kernel(s). The surgery
obstruction theory is much easier in the even-dimensional case $m=2n$
when $\pi_r(f)$ can be non-zero at most for $r=m+1$
than in the odd-dimensional case $m=2n+1$ when 
$\pi_r(f)$ can be non-zero for $r=m+1$ and $r=m+2$.

Wall [28,\S18G] asked for a chain complex formulation of surgery, in
which the obstruction groups $L_m(\Lambda)$ would appear as the
cobordism groups of chain complexes with $m$-dimensional quadratic
Poincar\'e duality, by analogy with the cobordism groups of manifolds
$\Omega_*$. Mishchenko [15] initiated such a theory of
``$m$-dimensional symmetric Poinc\-ar\'e complexes'' $(C,\phi)$ with
$C$ an $m$-dimensional f.\ g. free $\Lambda$-module chain complex
$$C~:~C_m \raise4pt\hbox{$d \atop \ra{1.5} $} C_{m-1} 
\raise4pt\hbox{$d \atop \ra{1.5} $} C_{m-2} \to \dots 
 \to C_1 \raise4pt\hbox{$d \atop \ra{1.5} $} C_0$$
and $\phi$ a quadratic structure inducing $m$-dimensional Poincar\'e
duality isomorphisms $\phi_0:H^*(C)\to H_{m-*}(C)$. 
The cobordism groups $L^m(\Lambda)$ (which are covariant in $\Lambda$) 
are such that for any $m$-dimensional geometric Poincar\'e
complex $X$ there is defined a symmetric signature invariant
$$\sigma^*(X)~=~(C(\widetilde{X}),\phi) \in L^m(\bbz[\pi_1(X)])~.$$
The corresponding quadratic theory was developed
in Ranicki [19]; the $m$-dimensional quadratic $L$-groups
$L_m (\Lambda)$ for any $m \ge 0$ were obtained as the groups of equivalence classes
of ``$m$-dimensional quadratic Poincar\'e complexes'' $(C,\psi)$.
The surgery obstruction of an $m$-dimensional normal map
$(f,b):M^m \to X$ was expressed as a cobordism class
$$\sigma_*(f,b)~=~(C,\psi) \in L_m(\bbz[\pi_1(X)])$$
with
$$H_*(C)~=~K_*(M)~=~H_{*+1}( \widetilde f : \widetilde M \to 
\widetilde X)~.$$
The symmetrization maps
$1+T:L_*(\Lambda) \to L^*(\Lambda)$ are isomorphisms modulo 8-torsion, and the
symmetrization of the surgery obstruction is the difference of the
symmetric signatures
$$(1+T)\sigma_*(f,b)~=~ \sigma^*(M)-\sigma^*(X)\in L^m(\bbz[\pi_1(X)])~.$$
However, the theory of [19] is fairly elaborate. 
The algebra of [18] and [19] is used in these notes
to simplify the original theory of Wall [28]
in the odd-dimensional case, without invoking the full theory of [19].

\noindent {\bf 0.2 What is in these notes}

These notes give an elementary account of the construction of the
$L$-groups $L_*$ and the surgery obstruction $\sigma_*$ for
differentiable manifolds.  For the more computational aspects of the
$L$-groups see the papers by Hambleton and Taylor [4] and Stark [25].

The even-dimensional $L$-groups $L_{2n}(\Lambda)$ are the Witt groups
of nonsingular $(-1)^n$-quadratic forms over $\Lambda$. 
It is relatively easy to pass from an $n$-connected 
$2n$-dimensional normal map $(f,b):M^{2n} \to X$ to a $(-1)^n$-quadratic
form representing $\sigma_*(f,b)$, and to see how the form changes
under a surgery on $(f,b)$. This will be done in \S\S 1--5 of these notes.

The odd-dimensional $L$-groups $L_{2n+1}(\Lambda)$ are the stable 
automorphism groups of nonsingular $(-1)^n$-quadratic forms over $\Lambda$. 
It is relatively hard to pass from an $n$-connected 
$(2n+1)$-dimensional normal map $(f,b):M^{2n+1} \to X$ to an
automorphism of a $(-1)^n$-quadratic form representing $\sigma_*(f,b)$,
and even harder to follow through in algebra the effect of a surgery on 
$(f,b)$. Novikov [16] suggested the reformulation of the odd-dimensional 
theory in terms of the language of hamiltonian physics, and to
replace the automorphisms by ordered pairs of lagrangians (= maximal
isotropic subspaces). This reformulation was carried out in Ranicki [18],
where such pairs were called `formations', but it was still hard to follow 
the algebraic effects of individual surgeries. This became easier after
the further reformulation of Ranicki [19] in terms of chain complexes with 
Poincar\'e duality -- see \S\S8,9 for a description of how the kernel 
formation changes under a surgery on $(f,b)$.

The original definition of $L_*(\Lambda)$ in Wall [28] was for the category of based f.\ g. free
$\Lambda$-modules and simple isomorphisms, for surgery up to simple homotopy
equivalence. Here, f.\ g. stands for finitely generated and simple
means that the Whitehead torsion is trivial, as in the hypothesis of
the $s$-cobordism theorem. These notes will only deal with
free $L$-groups $L_*(\Lambda)=L^h_*(\Lambda)$, the obstruction groups for surgery
up to homotopy equivalence. 

The algebraic theory of $\epsilon$-quadratic forms $(K,\lambda,\mu)$ 
over a ring $\Lambda$ with an involution $\Lambda \to \Lambda ; a \mapsto \bar a$ 
is developed in \S\S1,2, with $\epsilon = \pm 1$ and
$$\lambda ~:~ K \times K \to \Lambda~;~(x,y)~\mapto \lambda(x,y)$$
an $\epsilon$-symmetric pairing on a $\Lambda$-module $K$ such that
$$\lambda(x,y) ~=~ \epsilon \overline {\lambda(y,x)} \in \Lambda ~~ (x,y \in K)$$
and
$$\mu ~:~ K \to Q_{\epsilon}(\Lambda)~=~\Lambda/\{a - \epsilon \bar a \,\vert\, a \in \Lambda\}$$
an $\epsilon$-quadratic refinement of $\lambda$ such that
$$\lambda(x,x) ~=~ \mu(x) + \epsilon \overline{\mu(x)} \in \Lambda~~ (x \in K)~.$$ 
For an $n$-connected $2n$-dimensional normal map $(f,b):M^{2n}\to X$ geometric
(intersection,\t self-intersection) numbers define a $(-1)^n$-quadratic
form $(K_n(M),\lambda,\mu)$ on the kernel stably f.\ g. free
$\bbz[\pi_1(X)]$-module
$$K_n(M)~=~{\rm ker}(\widetilde f_* : H_n(\widetilde M) \to H_n(\widetilde X))$$
with $\widetilde X$ the universal cover of $X$ and $\widetilde M= f^*\widetilde X$
the pullback of $\widetilde X$ to $M$. 

The hyperbolic $\epsilon$-quadratic form on a f.\ g. free $\Lambda$-module $\Lambda^k$
$$H_{\epsilon}(\Lambda^k)~=~(\Lambda^{2k}, \lambda, \mu)$$
is defined by
$$\eqalign{&\lambda ~:~ \Lambda^{2k} \times \Lambda^{2k} \to \Lambda ~;\ccr
&\indent ((a_1,a_2,\dots,a_{2k}),(b_1,b_2,\dots,b_{2k}))~ \mapto 
\sum^k _{i=1}(b_{2i-1} \bar a_{2i} + \epsilon b_{2i} \bar a_{2i-1})~,\ccr
&\mu ~:~ \Lambda^{2k} \to Q_{\epsilon}(\Lambda) ~;~ 
(a_1,a_2,\dots,a_{2k}) \mapto \sum^k _{i=1}a_{2i-1} \bar a_{2i}~.}$$
\indent
The even-dimensional $L$-group $L_{2n}(\Lambda)$ is defined in \S3 to
be the abelian group of stable isomorphism classes of nonsingular
$(-1)^n$-quadratic forms on (stably) f.\ g. free $\Lambda$-modules, where
stabilization is with respect to the hyperbolic forms
$H_{(-1)^n}(\Lambda^k)$. A nonsingular $(-1)^n$-quadratic form
$(K,\lambda,\mu)$ represents 0 in $L_{2n}(\Lambda)$ if and only if
there exists an isomorphism
$$(K,\lambda,\mu) \oplus H_{(-1)^n}(\Lambda^k)~\cong~H_{(-1)^n}(\Lambda^{k'})$$
for some integers $k,k' \ge 0$. The surgery obstruction of an
$n$-connected $2n$-dimensional normal map $(f,b):M^{2n}\to X$ is
defined by
$$\sigma_*(f,b)~=~(K_n(M),\lambda,\mu) \in L_{2n}(\bbz[\pi_1(X)])~.$$
\indent
The algebraic effect of a geometric surgery on an $n$-connected
$2n$-dimension\-al normal map $(f,b)$ is given in \S5. Assuming that
the result of the surgery is still $n$-connected, the effect on the
kernel form of a surgery on $S^{n-1} \times D^{n+1} \hookrightarrow M$ (resp. 
$S^n \times D^n \hookrightarrow M$) is to add on (resp. split off) a
hyperbolic $(-1)^n$-quadratic form $H_{(-1)^n}(\bbz[\pi_1(X)])$.

\S6 introduces the notion of a ``$(2n+1)$-complex'' $(C,\psi)$,
which is a f.\ g. free $\Lambda$-module chain complex of the type
$$C~:~\dots \to 0\to C_{n+1}~ \raise4pt\hbox{$d \atop \ra{1.5} $}~ C_n~ 
\to 0\to \dots $$
together with a quadratic structure $\psi$ inducing Poincar\'e
duality isomorphisms $(1+T)\psi: H^{2n+1-*}(C) \to H_*(C)$. 
(This is just a $(2n+1)$-dimensional quadratic Poincar\'e complex 
$(C,\psi)$ in the sense of [19], with $C_r=0$ for $r \neq n,n+1$.)
An $n$-connected $(2n+1)$-dimensional normal map 
$(f,b):M^{2n+1}\to X$ determines a kernel $(2n+1)$-complex $(C,\psi)$ 
(or rather a homotopy equivalence class of such complexes) with
$$H_*(C)~=~K_*(M)~=~ {\rm ker} (\widetilde f_* : H_*(\widetilde M) \to H_*(\widetilde X))~.$$
The cobordism of $(2n+1)$-complexes is defined in \S7.
The odd-dimensional $L$-group $L_{2n+1}(\Lambda)$ is defined in \S8
as the cobordism group of $(2n+1)$-complexes. The surgery obstruction 
of an $n$-connected normal map
$(f,b):M^{2n+1}\to X$ is the cobordism class of the kernel complex
$$\sigma_*(f,b) ~=~(C,\psi) \in L_{2n+1}(\bbz[\pi_1(X)])~.$$
\indent
The odd-dimensional $L$-group $L_{2n+1}(\Lambda)$ was originally
defined in [28] as a potentially non-abelian quotient of the stable
unitary group of the matrices of automorphisms of hyperbolic
$(-1)^n$-quadratic forms over $\Lambda$
$$L_{2n+1}(\Lambda)~=~U_{(-1)^n}(\Lambda)/EU_{(-1)^n}(\Lambda)$$
with
$$U_{(-1)^n}(\Lambda)~=~
\bigcup\limits^{\infty}_{k=1}{\rm Aut}_{\Lambda}H_{(-1)^n}(\Lambda^k)$$
and $EU_{(-1)^n}(\Lambda)\triangleleft\, U_{(-1)^n}(\Lambda)$ 
the normal subgroup generated by the elementary matrices of the type
$$\pmatrix{\alpha & 0\ccr0 & \alpha^{*-1}}~~,~~
\pmatrix{1& 0 \ccr \beta +(-1)^{n+1} \beta^* & 1}~~,~~\pmatrix{0 & 1 \ccr (-1)^n & 0 }$$ 
for any invertible matrix $\alpha$, and any square matrix $\beta$.
The group $L_{2n+1}(\Lambda)$ is abelian, since
$$[U_{(-1)^n}(\Lambda),U_{(-1)^n}(\Lambda)] \subseteq EU_{(-1)^n}(\Lambda)~.$$
The surgery obstruction $\sigma_*(f,b) \in L_{2n+1}(\bbz[\pi_1(X)])$ of an
$n$-connected $(2n+1)$-dimensional normal map $(f,b):M^{2n+1}\to X$ is
represented by an automorphism of a hyperbolic $(-1)^n$-quadratic form
obtained from a high-dimensional generalization of the Heegaard
decompositions of 3-dimensional manifolds as twisted doubles.

\S8,\t \S9 and \S10 describe three equivalent ways of defining
$L_{2n+1}(\Lambda)$, using unitary matrices, formations and chain
complexes. In each case it is necessary to make some choices in
passing from the geometry to the algebra, and to verify that the
equivalence class in the $L$-group is independent of the choices.

The definition of $L_{2n+1}(\Lambda)$ using complexes given in \S8 is a
special case of the general theory of chain complexes with Poincar\'e
duality of Ranicki [19].  The 4-periodicity in the quadratic $L$-groups
$$L_m(\Lambda)~=~L_{m+4}(\Lambda)$$
(given geometrically by taking product with $\bbc\,\bbp^{\,2}$, as in
Chapter 9 of Wall [28]) was proved in [19] using an algebraic analogue
of surgery below the middle dimension: it is possible to represent
every element of $L_m (\Lambda)$ by a quadratic Poincar\'e complex
$(C,\psi)$ which is ``highly-connected'', meaning that
$$C_r~=~0~{\rm for}~\cases{r \ne n & {\rm if}~$m=2n$ \ccr
r\ne n,n+1 & {\rm if}~$m=2n+1$~.}$$
In these notes only the highly-connected 
$(2n+1)$-dimensional quadratic Poincar\'e complexes are considered, 
namely the ``$(2n+1)$-complexes'' of \S6. 

I am grateful to the referee for suggesting several improvements.

\indent	The titles of the sections are:
\vskip 1pc
\centerline{\vbox{
\hbox{\S1. Duality}
\hbox{\S2. Quadratic forms}
\hbox{\S3. The even-dimensional $L$-groups}
\hbox{\S4. Split forms}
\hbox{\S5. Surgery on forms}
\hbox{\S6. Short odd complexes}
\hbox{\S7. Complex cobordism}
\hbox{\S8. The odd-dimensional $L$-groups}
\hbox{\S9. Formations}
\hbox{\S10. Automorphisms}}}

\noindent {\sectionfont \S1. Duality}

\S1 considers rings $\Lambda$ equipped with an ``involution'' reversing
the order of multiplication.  An involution allows right
$\Lambda$-modules to be regarded as left $\Lambda$-modules, especially
the right $\Lambda$-modules which arise as the duals of left
$\Lambda$-modules.  In particular, the group ring $\bbz[\pi_1(M)]$ of
the fundamental group $\pi_1(M)$ of a manifold $M$ has an involution,
which allows the Poincar\'e duality of the universal cover
$\widetilde{M}$ to be regarded as $\bbz[\pi_1(M)]$-module isomorphisms.

Let $X$ be a connected space, and let $\widetilde{X}$ be a regular
cover of $X$ with group of covering translations $\pi$.  The action of
$\pi$ on $\widetilde{X}$ by covering translations
$$\pi \times \widetilde X \to \widetilde X ~;~ (g,x) \mapto gx$$
induces a left action of the group ring $\bbz[\pi]$ on the homology of
$\widetilde{X}$
$$\bbz[\pi] \times H_*(\widetilde X) \to H_*(\widetilde X)~;~
(\sum_{g \in \pi}n_g g, x) \mapto \sum_{g \in \pi}n_g gx$$
so that the homology groups $H_*(\widetilde{X})$ are left $\bbz[\pi]$-modules.
In dealing with cohomology let 
$$H^*(\widetilde{X})~=~H^*_{cpt}(\widetilde{X})$$
be the compactly supported cohomology groups, regarded as 
left $\bbz[\pi]$-modules by
$$\bbz[\pi] \times H^*(\widetilde X) \to H^*(\widetilde X)~;~
(\sum_{g \in \pi}n_g g,x) ~\mapto \sum_{g \in \pi}n_g xg^{-1}~ .$$
(For finite $\pi$ $H^*(\widetilde{X})$ is just the ordinary cohomology
of $\widetilde{X}$.)
Cap product with any homology class $[X] \in H_m(X)$ 
defines $\bbz[\pi]$-module morphisms
$$[X] \cap -~:~ H^*(\widetilde X) \to H_{m-*}(\widetilde X)~.$$
\noindent{\sc Definition 1.1} An 
{\it oriented $m$-dimensional geometric Poincar\'e complex} (Wall [27])
is a finite $CW$ complex $X$ with a fundamental class $[X] \in H_m(X)$ 
such that cap product defines $\bbz[\pi_1(X)]$-module isomorphisms
$$[X] \cap -~:~ H^*(\widetilde X) ~\isa~ H_{m-*}(\widetilde X)$$
with $\widetilde{X}$ the universal cover of $X$.\hfill$\square$

See 1.14 below for the general definition of a geometric Poincar\'e
complex, including the nonorientable case.

\noindent{\sc Example 1.2} A compact oriented $m$-dimensional manifold is an
oriented $m$-dimensional geometric Poincar\'e complex.\hfill$\square$

In order to also deal with nonorientable manifolds and Poincar\'e complexes
it is convenient to have an involution: 

\noindent{\sc Definition 1.3}
Let $\Lambda$ be an associative ring with 1. An {\it involution}
on $\Lambda$ is a function 
$$\Lambda \to \Lambda ~;~a \mapto \overline{a}$$
satisfying
$${\overline {(a+b)}} ~ =~ \overline{a} + {\overline b}~,~
 {\overline {(ab)}}~=~{\overline b}.\overline{a} ~,~
{\overline {\overline{a}}}~=~a~,~{\overline 1}~=~1 \in \Lambda~.\eqno{\square}$$

\noindent{\sc Example 1.4}
A commutative ring $\Lambda$ admits the identity involution
 $$\Lambda \to \Lambda~;~a \mapto \overline{a}~=~a~.\eqno{\square}$$

\noindent{\sc Definition 1.5} Given a group $\pi $ and a group morphism 
$$w~:~\pi \to\bbz_2 ~=~\{\pm 1\}$$ 
define the $w$-{\it twisted involution} on the integral group ring 
$\Lambda=\bbz[\pi]$
$$\hfil \Lambda \to \Lambda~;~a~=~\sum_{g \in \pi} n_g g \mapto 
\bar a~=~\sum_{g \in \pi} w(g)n_gg^{-1}~ (n_g \in\bbz)~.\eqno{\square}$$

In the topological application $\pi$ is the fundamental group of a
space with $w:\pi \to \bbz_2$ an orientation character. 
In the oriented case $w(g)=+1$ for all $g \in \pi$.

\noindent {\sc Example 1.6}
Complex conjugation defines an involution on the ring of 
complex numbers $\Lambda=\bbc$
$$\bbc \to \bbc~;~ z~=~a+ib \mapto {\overline z}=a-ib~.\eqno{\square}$$

A ``hermitian'' form is a symmetric form on a (finite-dimensional)
vector space over $\bbc$ with respect to this involution. The study of
forms over rings with involution is sometimes called ``hermitian
$K$-theory'', although ``algebraic $L$-theory'' seems preferable.

The dual of a left $\Lambda$-module $K$ is the right $\Lambda$-module
$$K^*~=~{\rm Hom}_\Lambda(K,\Lambda)$$
with
$$K^* \times \Lambda \to K^*~;~(f,a) \mapto (x \mapto f(x).a)~.$$
An involution $\Lambda\to \Lambda;a \mapto \bar a$ determines an
isomorphism of categories
$$\{\hbox{\rm right $\Lambda$-modules}\}~\isa~
\{\hbox{\rm left $\Lambda$-modules}\}~;~L \mapto L^{op}~,$$
with $L^{op}$ the left $\Lambda$-module with the same additive group as
the right $\Lambda$-module $L$ and $\Lambda$ acting by
$$\Lambda \times L^{op} \to L^{op}~;~ (a,x) \mapto x\overline{a}~.$$

From now on we shall work with a ring $\Lambda$ which is equipped with
a particular choice of involution $\Lambda\to \Lambda$.  Also,
$\Lambda$-modules will always be understood to be left $\Lambda$-modules.  

For any $\Lambda$-module $K$ the $\Lambda$-module $(K^*)^{op}$ is written 
as $K^*$. Here is the definition of $K^*$ all at once:

\noindent {\sc Definition 1.7} 
The {\it dual} of a $\Lambda$-module $K$ is the $\Lambda$-module 
$$K^* ~=~{\rm Hom}_\Lambda (K,\Lambda)~,$$
with $\Lambda$ acting by
$$\Lambda \times K^* \to K^*~;~ (a,f) \mapto (x \mapto f(x).\overline {a})$$
for all $a\in \Lambda$, $f \in K^*$, $x\in K$.\hfill$\square$

There is a corresponding notion for morphisms:

\noindent {\sc Definition 1.8} The {\it dual}
of a $\Lambda$-module morphism $f:K\to L$ is the $\Lambda$-module morphism
$$f^*~:~ L^* \to K^*~;~ g \mapto \big(x \mapto 
g(f(x))\big)~.\eqno{\square}$$

Thus duality is a contravariant functor
$$*~:~\{\Lambda\hbox{\rm -modules}\} \to 
\{\Lambda\hbox{\rm-modules}\}~;~K \mapto K^*~.$$

\noindent {\sc Definition 1.9}
For any $\Lambda$-module $K$ define the $\Lambda$-module morphism
$$e_K~:~K \to K^{**}~;~ x \mapto (f \mapto \overline{f(x)})~.\eqno{\square}$$

The morphism $e_K$ is natural in the sense that for any $\Lambda$-module
 morphism $f:K\to L$ there is defined a commutative diagram
$$\xymatrix@C+4pt@R+4pt{
K \ar[r]^{\displaystyle{f}} \ar[d]_{\displaystyle{e_K}} & L \ar[d]^{\displaystyle{e_L}} \\
K^{**} \ar[r]^{\displaystyle{f^{**}}} & L^{**}}$$

\noindent {\sc Definition 1.10}
(i) A $\Lambda$-module $K$ is {\it f.\ g. projective}
if there exists a $\Lambda$-module $L$ such that $K\oplus L$ is isomorphic 
to the f.\ g. free $\Lambda$-module $\Lambda^n$, for some $n\ge 0$. \hfil\break
(ii) A $\Lambda$-module $K$ is {\it stably f.\ g. free}
if $K\oplus \Lambda^m $ is isomorphic to $\Lambda^n$, for some $m,n\ge 0$.\hfill$\square$

In particular, f.\ g. free $\Lambda$-modules are stably f.\ g. free, and stably
f.\ g. free $\Lambda$-modules are f.\ g. projective.

\noindent {\sc Proposition 1.11}
{\it The dual of a f.\ g. projective $\Lambda$-module $K$ is a f.\ g.
 projective $\Lambda$-module $K^*$, and $e_K:K\to K^{**}$ is an isomorphism.
 Moreover, if $K$ is stably f.\ g. free then so is $K^*$.}\hfil\break
{\sc Proof}\t :
For any $\Lambda$-modules $K,L$ there are evident identifications
$$\eqalign{&(K\oplus L)^*~=~K^* \oplus L^*~,\ccr
&e_{K\oplus L}~=~e_K\oplus e_{L}~:~K\oplus L \to (K\oplus L)^{**}~
=~K^{**} \oplus L^{**}~,}$$
so it suffices to consider the special case $K=\Lambda$.
The $\Lambda$-module isomorphism
$$f~:~\Lambda ~\isa~ \Lambda^*~;~a \mapto (b \mapto b\overline{a})~,$$
can be used to construct an explicit inverse for $e_\Lambda $
$$(e_\Lambda )^{-1}~:~\Lambda^{**} \to \Lambda~;~ g \mapto g(f(1))~.\eqno{\square}$$

In dealing with f.\ g. projective $\Lambda$-modules $K$ use 
the natural isomorphism $e_K:K\cong K^{**}$ to identify
$K^{**}=K$~. For any morphism $f:K\to L$ of f.\ g.
 projective $\Lambda$-modules there is a corresponding identification
$$f^{**}~=~f~:~K^{**}~=~K \to L^{**}~=~L~.$$

\noindent {\sc Remark 1.12}
The additive group ${\rm Hom}_\Lambda (\Lambda^m,\Lambda^n)$ of the
morphisms\break $\Lambda^m \to \Lambda^n$ between f.\ g.  free
$\Lambda$-modules $\Lambda^m,\Lambda^n$ may be identified with the
additive group $M_{m,n}(\Lambda)$ of $m\times n$ matrices
$(a_{ij})_{1\le i\le m,1\le j\le n}$ with entries $a_{ij} \in \Lambda$,
using the isomorphism
$$\eqalign{&M_{m,n}(\Lambda) ~\isa~ {\rm Hom}_\Lambda(\Lambda^m,\Lambda^n)~;\ccr
&(a_{ij}) \mapto ((x_1,x_2,\dots ,x_m) \mapto (\sum_{i=1}^m  x_ia_{i1},
\sum_{i=1}^m x_ia_{i2},\dots ,\sum_{i=1}^m x_ia_{in}))~.}$$
The composition of morphisms
$$\eqalign{&{\rm Hom}_\Lambda(\Lambda^m,\Lambda^n) \times {\rm Hom}_\Lambda(\Lambda^n,\Lambda^p) \to 
{\rm Hom}_\Lambda(\Lambda^m,\Lambda^p)~;\ccr
& \hbox to 1in{\hss}(f,g) \mapto \big(gf:x\mapto (gf)(x)=g(f(x))\big)}$$
corresponds to the multiplication of matrices
$$\displaylines{M_{m,n}(\Lambda) \times M_{n,p}(\Lambda) \to M_{m,p}(\Lambda)~;~
((a_{ij}),(b_{jk})) \mapto (c_{ik})\ccr
(\,c_{ik}~=~\sum^n_{j=1}a_{ij}b_{jk}~~(1\le i\le m,~1\le k\le p)\,)~.}$$
Use the isomorphism of f.\ g. free $\Lambda$-modules
$$\Lambda^m ~\isa~ (\Lambda^m)^*~;~(x_1,x_2,\dots ,x_m) \mapto ((y_1,y_2,\dots ,y_m) ~
\mapto \sum_{i=1}^m  y_i{\overline x}_i)$$
to identify
$$(\Lambda^m)^*~=~\Lambda^m ~.$$
The duality isomorphism
$$\eqalign{*~:~{\rm Hom}_\Lambda(\Lambda^m,\Lambda^n) ~&\isa~ {\rm Hom}_\Lambda((\Lambda^n)^*,(\Lambda^m)^*)~=~
 {\rm Hom}_\Lambda(\Lambda^n,\Lambda^m)~;\ccr
 &f \mapto f^*}$$
can be identified with the isomorphism defined by conjugate transposition of matrices
$$M_{m,n}(\Lambda) ~\isa~ M_{n,m}(\Lambda)~;~\alpha~=~(a_{ij}) \mapto 
 \alpha^*~=~(b_{ji})~,~b_{ji} = \overline{a}_{ij}~.\eqno{\square}$$

\noindent {\sc Example 1.13}
A $2\times 2$ matrix
$$\pmatrix{a&b \ccr c & d} \in M_{2,2}(\Lambda)$$
corresponds to the $\Lambda$-module morphism
$$f~=~\pmatrix{a & b \ccr c & d }~:~ \Lambda\oplus \Lambda \to \Lambda\oplus \Lambda~;~
(x,y) \mapto (xa+yb,xc+yd)~.$$
The conjugate transpose matrix
$$\pmatrix{\overline{a} & {\overline c} \ccr
{\overline b} & {\overline d}} \in M_{2,2}(\Lambda)$$
corresponds to the dual $\Lambda$-module morphism
$$\eqalign{f^*~=~\pmatrix{\overline{a} & {\overline c} \ccr
{\overline b} & {\overline d} }~:~
(\Lambda\oplus \Lambda)^*~=~& \Lambda\oplus \Lambda \to (\Lambda\oplus \Lambda)^*~=~\Lambda\oplus \Lambda~;\ccr
&(x,y)\mapto (x\overline{a}+y{\overline c},x{\overline b}+y{\overline d})~.}
\eqno{\lower18pt\hbox{$\square$}}$$

The dual of a chain complex of modules over a ring with involution $\Lambda$
$$C~:~\dots \xymatrix{\ar[r]&} C_{r+1} \xymatrix{\ar[r]^{\displaystyle{d}}&} C_r 
\xymatrix{\ar[r]^{\displaystyle{d}}&} C_{r-1} \xymatrix{\ar[r]&} \dots$$
is the cochain complex
$$C^*~:~\dots \xymatrix{\ar[r]&} C^{r-1} \xymatrix{\ar[r]^{\displaystyle{d^*}}&} C^r 
\xymatrix{\ar[r]^{\displaystyle{d^*}}&} C^{r+1} \xymatrix{\ar[r]&} \dots$$
with 
$$C^r~=~(C_r)^*~=~{\rm Hom}_{\Lambda}(C_r,\Lambda)~.$$

\noindent{\sc Definition 1.14} 
An {\it $m$-dimensional geometric Poincar\'e complex} (Wall [27])
is a finite $CW$ complex $X$ with an orientation character
$w(X):\pi_1(X)\to \bbz_2$ and a $w(X)$-twisted fundamental class $[X] \in H_m(X;\bbz^{w(X)})$ 
such that cap product defines $\bbz[\pi_1(X)]$-module isomorphisms
$$[X] \cap -~:~ H^*_{w(X)}(\widetilde X) ~\isa~ H_{m-*}(\widetilde X)$$
with $\widetilde{X}$ the universal cover of $X$. The $w(X)$-twisted
cohomology groups are given by
$$H^*_{w(X)}(\widetilde{X})~=~H^*(C(\widetilde{X})^*)$$
with $C(\widetilde{X})$ the cellular $\bbz[\pi_1(X)]$-module
chain complex, using the $w(X)$-twisted involution on $\bbz[\pi_1(X)]$ (1.5)
to define the left $\bbz[\pi_1(X)]$-module structure on the dual cochain
complex
$$C(\widetilde{X})^*~=~{\rm Hom}_{\bbz[\pi_1(X)]}
(C(\widetilde{X}),\bbz[\pi_1(X)])~.\eqno{\square}$$

The orientation character $w(X):\pi_1(X) \to \bbz_2$ sends a loop
$g:S^1 \to X$ to $w(g)=+1$ (resp.  $=-1$) if $g$ is
orientation-preserving (resp.  orientation-reversing).  

An oriented Poincar\'e complex $X$ (1.1) is just a Poincar\'e complex
(1.14) with $w(X)=+1$.

\noindent{\sc Example 1.15} A compact $m$-dimensional manifold is an
$m$-dimensional geometric Poincar\'e complex.\hfill$\square$

\noindent {\sectionfont \S2. Quadratic forms}

In the first instance suppose that the ground ring $\Lambda$
is commutative, with the identity involution $\bar a = a$ (1.4).
A symmetric form $(K,\lambda)$ over $\Lambda$ 
is a $\Lambda$-module $K$ together with a bilinear pairing
$$\lambda ~:~ K \times K \to \Lambda~;~(x,y) \mapto \lambda(x,y)$$
such that for all $x,y,z \in K$ and $a \in \Lambda$
$$\eqalign{&\lambda(x,ay) ~=~a \lambda(x,y) ~,\cr
&\lambda(x,y+z) ~=~ \lambda(x,y) + \lambda(x,z)~, \ccr 
&\lambda(x,y) ~=~ \lambda(y,x) \in \Lambda~.}$$
A quadratic form $(K,\lambda,\mu)$ over $\Lambda$ is a symmetric form
$(K,\lambda)$ together with a function
$$\mu ~:~ K \to Q_{+1}(\Lambda)~=~\Lambda~;~x\mapto \mu (x)$$
such that for all $x,y \in K$ and $a \in \Lambda$
$$\eqalign{&\mu(x+y)~=~\mu (x) + \mu (y) + \lambda(x,y)~,\cr
&\mu(ax)~=~a^2 \mu(x) \in  Q_{+1}(\Lambda) ~.}$$
In particular, for every $x \in K$
$$2 \mu(x)~ =~ \lambda(x,x) \in Q^{+1}(\Lambda) ~=~ \Lambda ~.$$
If $2 \in \Lambda$ is invertible (e.g. if $\Lambda$ is a field of characteristic
$\ne 2$, such as $\bbr,\bbc,\bbq$) there is no 
difference between symmetric and quadratic forms, with $\mu$ determined
by $\lambda$ according to $\mu (x) = \lambda(x,x)/2$. 

A symplectic form $(K,\lambda)$ over a commutative ring $\Lambda$ 
is a $\Lambda$-module $K$ together with a bilinear pairing
$\lambda : K \times K \to \Lambda$ such that for all $x,y,z \in K$ and $a \in \Lambda$
$$\eqalign{&\lambda(x,at) ~=~a \lambda(x,y) ~,\cr
&\lambda(x,y+z) ~=~ \lambda(x,y) + \lambda(x,z) ~,\cr
&\lambda(x,y) ~=~ - \lambda(y,x) \in \Lambda~.}$$
A $(-1)$-quadratic form $(K,\lambda,\mu)$ over $\Lambda$ is a symplectic form
$(K,\lambda)$ together with a function
$$\mu ~:~ K \to Q_{-1}(\Lambda)~=~\Lambda/\{2a \,\vert\, a \in \Lambda\}~;~x\mapto \mu (x)$$
such that for all $x,y \in K$ and $a \in \Lambda$
$$\eqalign{&\mu(x+y)~=~\mu (x) + \mu (y) + \lambda(x,y)~,\cr
&\mu(ax)~=~a^2 \mu(x) \in  Q_{-1}(\Lambda) ~.}$$
In particular, for every $x \in K$
$$2 \mu(x) ~=~ \lambda(x,x) \in Q^{-1}(\Lambda) ~=~ 
\{a \in \Lambda \,\vert\, 2a =0\} ~.$$
If $2 \in \Lambda$ is invertible then $Q_{-1}(\Lambda) = 0$ and there is no 
difference between symplectic and $(-1)$-quadratic forms, with $\mu = 0$. 

In the applications of forms to surgery theory it is
necessary to work with quadratic and $(-1)$-quadratic forms 
over noncommutative group rings with the involution
as in 1.5. \S2 develops the general theory of
forms over rings with involution, taking account of these differences. 

Let $X$ be an $m$-dimensional geometric Poincar\'e complex
with universal cover $\widetilde X$ and fundamental group ring
$\Lambda =\bbz [\pi_1(X)]$, with the $w(X)$-twisted involution. 
The Poincar\'e duality isomorphism
$$\phi ~=~ [X] \cap -~:~ H^{m-r}_{w(X)}(\widetilde X)~\isa~ H_r(\widetilde X)$$
and the evaluation pairing
$$H^r_{w(X)}(\widetilde X) \to H_r(\widetilde X)^* ~=~ 
{\rm Hom}_\Lambda(H_r(\widetilde X),\Lambda)~;
y \mapto (x \mapto \langle y,x\rangle)$$
can be combined to define a sesquilinear pairing
$$\lambda ~:~ H_r(\widetilde X) \times H_{m-r}(\widetilde X) \to \Lambda ~;~
(x,\phi (y)) \mapto \langle y,x\rangle $$
such that
$$\lambda(y,x) ~=~ (-1)^{r(m-r)} \overline {\lambda(x,y)}~$$
with $\Lambda \to \Lambda ; a \mapto \bar a$ the involution of 1.5. 

If $M$ is an $m$-dimensional manifold with fundamental group ring
$\Lambda =\bbz[\pi_1(M)]$ the pairing $\lambda :H_r(\widetilde M)
\times H_{m-r}(\widetilde M) \to \Lambda$ can be interpreted
geometrically using the geometric intersection numbers of cycles.  For
any two immersions $x:S^r \looparrowright \widetilde M$, 
$y : S^{m-r} \looparrowright \widetilde M$ in general position
$$\lambda(x,y)~=~ \sum_{g \in \pi_1(M)}n_g g \in \Lambda$$
with $n_g \in\bbz$ the algebraic number of intersections in
$\widetilde M$ of $x$ and $gy$. In particular, for $m=2n$ there is
defined a $(-1)^n$-symmetric pairing
$$\lambda ~:~ H_n(\widetilde M) \times H_n(\widetilde M) \to \Lambda$$
which is relevant to surgery in the middle dimension $n$. An element
$x \in \pi_n(M)$ can be killed by surgery if and only if it is
represented by an embedding $S^n \times D^n \hookrightarrow M^{2n}$.
The condition that the Hurewicz image $x\in H_n(\widetilde M)$ 
be such that $\lambda(x,x) = 0 \in \Lambda$ is necessary but not
sufficient to kill $x \in \pi_n(M)$ by surgery. 
The theory of forms developed in \S2 is required for an algebraic
formulation of the necessary and sufficient condition for an element in
the kernel $K_n(M)$ of an $n$-connected $2n$-dimensional normal map
$(f,b):M \to X$ to be killed by surgery, assuming $n \ge 3$.

As in \S1 let $\Lambda$ be a ring with involution, not necessarily a group
ring.

\noindent {\sc Definition 2.1}
A {\it sesquilinear pairing} $(K,L,\lambda)$ on $\Lambda$-modules $K,L$
is a function
$$\lambda~:~K\times L \to \Lambda~;~(x,y) \mapto \lambda(x,y)$$
such that for all $w,x\in K$, $y,z \in L$, $a,b \in \Lambda$
{\parindent=35pt
\parskip=2pt
\item{\rm (i)} $\lambda(w+x,y+z)~=~\lambda(w,y)+ \lambda(w,z)+
\lambda(x,y)+ \lambda(x,z)\in \Lambda$~,
\item{\rm (ii)} $\lambda(ax,by)~=~b\lambda(x,y)\overline{a} \in \Lambda~.$\par}
\noindent The {\it dual} (or {\it transpose}) sesquilinear pairing is
$$T\lambda~:~L \times K \to \Lambda~;~(y,x) \mapto T \lambda(y,x)~
=~\overline{\lambda(x,y)}~.\eqno{\square}$$

\noindent {\sc Definition 2.2}
Given $\Lambda$-modules $K,L$ let $S(K,L)$ be the additive group of
sesquilinear pairings $\lambda:K \times L \to \Lambda$. Transposition
defines an isomorphism
$$T~:~S(K,L)~\cong~S(L,K)$$
such that
$$T^2~=~{\rm id.}~:~S(K,L) ~\isa~ S(L,K) ~\isa~ S(K,L)~.\eqno{\square}$$

\noindent {\sc Proposition 2.3}
{\it For any $\Lambda$-modules $K,L$ there is a natural isomorphism of 
additive groups
$$\eqalign{&S(K,L)~\isa~{\rm Hom}_\Lambda(K,L^*)~;\ccr
&(\lambda\, :\,K\times L\to \Lambda)\mapto 
(\lambda\,:\,K \to L^*\,;\,x\mapto (\,y\mapto \lambda(x,y)\,))\,.}$$
For f.\ g. projective $K,L$ the transposition isomorphism 
$T:S(K,L)\cong\allowbreak S(L,K)$ corresponds to the duality isomorphism}
$$\eqalign{*\,:\,&{\rm Hom}_\Lambda (K,L^*) \,\isa\, {\rm Hom}_\Lambda(L,K^*)\,;\ccr
&(\lambda :K\to L^*) \mapto ({\lambda}^*\,:\,L\,\to\, K^*\,;\,
y \mapto (x\,\mapto \, {\overline{\lambda(y,x)}}))\,.}
\eqno{\lower12pt\hbox{$\square$}}$$

Use 2.3 to identify
$$\eqalign{&S(K,L)~=~{\rm Hom}_\Lambda(K,L^*)~,~S(K)~=~{\rm Hom}_\Lambda (K,K^*)~,\ccr
&Q^{\epsilon}(K)~=~{\rm ker}(1-T_{\epsilon}:{\rm Hom}_\Lambda (K,K^*)
\to {\rm Hom}_\Lambda (K,K^*))~,\ccr
&Q_{\epsilon}(K)~=~{\rm coker}(1-T_{\epsilon}:{\rm Hom}_\Lambda (K,K^*)
\to {\rm Hom}_\Lambda (K,K^*))}$$
for any f.\ g. projective $\Lambda$-modules $K,L$.

\noindent {\sc Remark 2.4}
For f.\ g. free $\Lambda$-modules $\Lambda^m$, $\Lambda^n$ it is possible to identify 
$S(\Lambda^m ,\Lambda^n)$ with the additive group $M_{m,n}(\Lambda)$ of $m\times n$
matrices $(a_{ij})$ with entries $a_{ij} \in \Lambda$, using the isomorphism
$$M_{m,n}(\Lambda) ~\isa~ S(\Lambda^m ,\Lambda^n)~;~(a_{ij}) \mapto \lambda $$
defined by
$$\lambda((x_1,x_2,\dots ,x_m),(y_1,y_2,\dots ,y_n))~=~\sum^m _{i=1}
 \sum^n_{j=1} y_ja_{ij}{\overline x}_i~.$$
The transposition isomorphism 
$T:S(\Lambda^m,\Lambda^n)\cong S(\Lambda^n,\Lambda^m)$ 
corresponds to the isomorphism defined by conjugate transposition of matrices
$$T~:~M_{m,n}(\Lambda) ~\isa~ M_{n,m}(\Lambda)~;~(a_{ij}) \mapto 
(b_{ji})~,~b_{ji}~=~\overline{a}_{ij}~.\eqno{\square}$$

The group $S(K,L)$ is particularly significant in the case $K=L$\t :
 
\noindent {\sc Definition 2.5} (i) Given a $\Lambda$-module $K$ let
$$S(K)~=~S(K,K)$$ 
be the abelian group of sesquilinear pairings $\lambda:K \times K \to \Lambda$.
\hfil\break
(ii) The $\epsilon$-{\it transposition involution} is given for 
$\epsilon=\pm 1$ by
$$T_{\epsilon}~:~S(K) ~\isa~ S(K)~;~\lambda \mapto 
T_{\epsilon} \lambda ~=~\epsilon (T\lambda)~,$$
such that
$$T_{\epsilon} \lambda(x,y) ~=~
\epsilon \overline {\lambda(y,x)} \in \Lambda~,~
(T_{\epsilon})^2~=~{\rm id}~:~ S(K) \to S(K)~.\eqno{\square}$$

\noindent {\sc Definition 2.6}
The $\epsilon$-{\it symmetric group}
of a $\Lambda$-module $K$ is the additive group
$$Q^{\epsilon}(K)~=~{\rm ker}(1-T_{\epsilon}:S(K)\to S(K))~.$$
The $\epsilon$-{\it quadratic group} of $K$ is the additive group
$$Q_{\epsilon}(K)~=~{\rm coker}(1-T_{\epsilon}:S(K)\to S(K))~.$$
The $\epsilon$-{\it symmetrization morphism} is given by
$$1+T_{\epsilon}~:~Q_{\epsilon}(K) \to Q^{\epsilon}(K)~;~
\psi\mapto \psi+ T_{\epsilon} \psi~.\eqno{\square}$$

For $\epsilon = +1$ it is customary to refer to $\epsilon$-symmetric
and $\epsilon$-quadratic objects as symmetric and quadratic, as in the
commutative case. 

For $K=\Lambda$ there is an isomorphism of additive groups with
involution
$$\Lambda ~\isa~ S(\Lambda)~;~a \mapto ( (x,y) \mapto ya{\overline x})$$
allowing the identifications
$$\eqalign{&Q^{\epsilon}(\Lambda)~=~\{a\in \Lambda\,\vert\, \epsilon \overline{a}=a\}~,\ccr
&Q_{\epsilon}(\Lambda)~=~\Lambda/\{a- \epsilon \overline{a}\,\vert\, a\in \Lambda\}~,\ccr
&1+T_{\epsilon}~:~Q_{\epsilon}(\Lambda) \to Q^{\epsilon}(\Lambda)~;~
a \mapto a+ \epsilon \overline{a}~.}$$

\noindent {\sc Example 2.7} Let $\Lambda=\bbz$. The
$\epsilon$-symmetric and $\epsilon$-quadratic groups of $K=\bbz$ are
given by
$$\eqalign{&
Q^{\epsilon} (\bbz)~=~\cases{\bbz&if~$\epsilon =+1$ \ccr 
0&if~$ \epsilon =-1$~,}\ccr
&Q_{\epsilon} (\bbz)~=~
\cases{\bbz&if $\epsilon =+1$ \ccr\bbz/2&if~$\epsilon =-1$ }}$$
with generators represented by $1\in\bbz$, and with
$$1+T_+~=~2~:~Q_{+1}(\bbz)~=~\bbz \to Q^{+1}(\bbz)~=~\bbz~.\eqno{\square}$$

\noindent {\sc Definition 2.8}
An $\epsilon$-{\it symmetric form} $(K,\lambda)$ over $\Lambda$ is a
$\Lambda$-module $K$ together with an element $\lambda \in Q^{\epsilon}(K)$.
Thus $\lambda $ is a sesquilinear pairing
$$\lambda ~:~K \times K \to \Lambda~;~(x,y) \mapto \lambda(x,y)$$
such that for all $x,y \in K$
$$\lambda(x,y) ~=~ \epsilon \overline{\lambda(y,x)} \in \Lambda~.$$
The {\it adjoint} of $(K,\lambda)$ is the $\Lambda$-module morphism
$$K \to K^* ~;~ x \mapto (y \mapto \lambda(x,y))$$
which is also denoted by $\lambda$. The form is {\it nonsingular}
if $\lambda : K \to K^*$ is an isomorphism.
\hfill$\square$

Unless specified otherwise, only forms $(K,\lambda)$ with $K$ a f.\t g. 
projective $\Lambda$-module will be considered.

\noindent {\sc Example 2.9}
 The symmetric form $(\Lambda,\lambda)$ defined by
$$\lambda~=~1 ~:~\Lambda \to \Lambda^*~;~a \mapto ( b \mapto b\overline{a})$$
is nonsingular. \hfill$\square$

\noindent {\sc Definition 2.10}
For any f.\ g. projective $\Lambda$-module $L$ define
the nonsingular {\it hyperbolic $\epsilon$-symmetric form}
$$H^{\epsilon}(L)~=~(L \oplus L^*, \lambda)$$
by
$$\eqalign{\lambda ~=~\pmatrix{0 & 1 \ccr \epsilon & 0}~
:~ L \oplus L^*& \to (L \oplus L^*)^* ~=~L^* \oplus L~; \ccr
&(x,f)\mapto ((y,g)\mapto f(y) + \epsilon \overline {g(x)})~.}
\eqno{\lower16pt\hbox{$\square$}}$$

\noindent {\sc Example 2.11} Let $X$ be an $m$-dimensional
geometric Poincar\'e complex, and let $\widetilde X$ be
a regular oriented covering of $X$ with group of
covering translations $\pi$ and orientation character $w:\pi \to \bbz_2$.
An element $g \in \pi$ has $w(g)=+1$ (resp. $-1$) if and only if the
covering translation $g:\widetilde{X} \to \widetilde{X}$ is
orientation-preserving (resp. reversing).\hfil\break
(i) Cap product with the fundamental class $[X] \in H_m(X;\bbz^w)$ 
defines the Poincar\'e duality $\bbz [\pi]$-module
isomorphisms
$$[X] \cap - ~:~ H^{m-*}_w(\widetilde X)~ \isa~ H_*(\widetilde X)~.$$ 
If $m=2n$ and $X$ is a manifold geometric intersection numbers define
a $(-1)^n$-symmetric form $(H_n(\widetilde X),\lambda)$ 
over $\bbz[\pi]$ with adjoint the composite
$$\lambda ~:~ H_n(\widetilde X)~ \raise4pt\hbox{$([X] \cap -)^{-1} \atop
\ra{4}$}~ H^n_w(\widetilde X)~ \raise4pt\hbox{${\rm evaluation}
\atop \ra{4}$}~H_n(\widetilde X)^*~.$$
(ii) In general $H_n(\widetilde X)$ is not a f.\ g. projective $\bbz[\pi]$-module. 
If $H_n(\widetilde X)$ is f.\ g. projective then
the evaluation map is an isomorphism, and $(H_n(\widetilde X),\lambda)$
is a nonsingular form.\hfill$\square$

\noindent{\sc Remark 2.12} 
(i) Let $M$ be a $2n$-dimensional manifold, with
universal cover $\widetilde M$ and intersection pairing 
$\lambda:H_n(\widetilde M) \times H_n(\widetilde M) \to \bbz[\pi_1(M)]$.
An element $x \in {\rm im}(\pi_n(M) \to H_n(\widetilde M))$
can be killed by surgery if and only if it can be represented by an
embedding $x : S^n \times D^n \hookrightarrow M$, in which case
the homology class $x \in H_n(\widetilde M)$ is such that $\lambda(x,x) = 0$. 
However, the condition $\lambda(x,x)=0$ given by the symmetric
structure alone is not sufficient for the existence of such an
embedding -- see (ii) below for an explicit example.\hfil\break
(ii) The intersection form over $\bbz$ for $M^{2n}=S^n \times S^n$ is the
hyperbolic form (2.10) 
$$(H_n(S^n \times S^n),\lambda)~=~H^{(-1)^n}(\bbz)~.$$ 
The  element $x =(1,1) \in H_n(S^n \times S^n)$ is such that
$$\lambda(x,x)~=~\chi(S^n)~=~1+(-1)^n\in\bbz~,$$
so that $\lambda(x,x)=0$ for odd $n$.
The diagonal embedding $\Delta:S^n \hookrightarrow S^n
\times S^n$ has normal bundle $\nu_{\Delta}=\tau_{S^n}:S^n \to BO(n)$,
which is non-trivial for $n \neq 1,3,7$, so that it is not possible to
kill $x=\Delta_*[S^n] \in H_n(S^n \times S^n)$ by surgery in these dimensions.
\hfill$\square$
 
\noindent {\sc Definition 2.13}
An $\epsilon$-{\it quadratic form} $(K,\lambda,\mu)$ over $\Lambda$
is an $\epsilon$-symmetric form $(K,\lambda)$ together with a function
$$\mu ~:~K \to Q_{\epsilon}(\Lambda)~;~x \mapto \mu (x)$$
such that for all $x,y \in K$, $a \in \Lambda$

{\parindent=50pt
\parskip=4pt
\item{\rm (i)} $\mu (x+y) - \mu (x) - \mu (y)~=~\lambda(x,y) \in Q_{\epsilon}(\Lambda)~,$
\item{\rm (ii)} $\mu (x) + \epsilon \overline{\mu (x)}~=~
\lambda(x,x) \in
{\rm im}(1+T_{\epsilon}:Q_{\epsilon}(\Lambda)\to Q^{\epsilon}(\Lambda))~,$ 
\item{(iii)} $\mu (ax)~=~a \mu (x)\overline{a} \in Q_{\epsilon}(\Lambda)~.$\par}
\line{\hfill$\square$}

\noindent {\sc Definition 2.14}
For any f.\ g. projective $\Lambda$-module $L$ define the nonsingular 
{\it hyperbolic $\epsilon$-quadratic form} over $\Lambda$ by
$$H_{\epsilon}(L)~=~(L \oplus L^*, \lambda , \mu)$$
with 
$$\eqalign{&\lambda~=~\pmatrix{0 & 1 \ccr \epsilon & 0}~
:~ L \oplus L^* \to (L \oplus L^*)^* ~=~L^* \oplus L~; \ccr
& \hskip100pt 
(x,f)\mapto ((y,g)\mapto f(y) + \epsilon \overline {g(x)})~,\ccr
&\mu ~:~ L \oplus L^* \to Q_{\epsilon}(\Lambda) ~;~(x,f)\mapto f(x)~.}$$
$(L\oplus L^*,\lambda)$ is the hyperbolic $\epsilon$-symmetric form
$H^{\epsilon}(L)$ of 2.10.\hfill$\square$

\noindent {\sc Example 2.15} (Wall [28,\t Chapter 5])
An $n$-connected normal map $(f,b):M^{2n}\to X$ from a $2n$-dimensional
manifold with boundary $(M,\partial M)$ to a geometric Poincar\'e pair
$(X,\partial X)$ with $\partial f =f\vert :\partial M \to \partial X$ a
homotopy equivalence determines a $(-1)^n$-quadratic form
$(K_n(M),\lambda ,\mu)$ over $\Lambda=\bbz[\pi_1(X)]$ with the 
$w(X)$-twisted involution (1.5), with
$$K_n(M)~=~\pi_{n+1}(f)~=~H_{n+1}(\widetilde f)~=~ 
{\rm ker}(\widetilde f_*:H_n(\widetilde M) \to H_n(\widetilde X))$$ 
the stably f.\ g. free kernel $\Lambda$-module, and $\widetilde f :
\widetilde M \to \widetilde X$ a $\pi_1(X)$-equivariant lift of $f$ to
the universal covers. Note that $K_n(M)=0$ if (and for $n \ge 2$ only if) $f:M \to X$ 
is a homotopy equivalence, by the theorem of J.H.C.\t Whitehead. 
\hfil\break
(i) The pairing $\lambda : K_n(M) \times K_n(M) \to
\Lambda$ is defined by geometric intersection numbers, as follows.
Every element $x \in K_n(M)$ is represented by an
$X$-nullhomotopic framed immersion $g:S^n \looparrowright M$ with
a choice of path in $g(S^n)\subset M$ from the base point 
$\ast \in M$ to $g(1)\in M$. Any two elements $x,y \in K_n(M)$ 
can be represented by such immersions $g,h:S^n \looparrowright M$ with
transverse intersections and self-intersections. 
The intersection of $g$ and $h$
$$D(g,h)~=~\{(a,b) \in S^n \times S^n\,\vert\,g(a)=h(b)\in M\}$$
is finite. For each intersection point $(a,b) \in D(g,h)$ let 
$$\gamma(a,b) \in \pi_1(M)~=~\pi_1(X)$$ 
be the homotopy class of the loop in $M$ obtained by joining the path 
in $g(S^n)\subset M$ from the base point $\ast \in M$ to $g(a)$
to the path in $h(S^n) \subset M$ from $h(b)$ back to the base point.
Choose an orientation for $\tau_*(M)$ and transport it to an
orientation for $\tau_{g(a)}(M)=\tau_{h(b)}(M)$ by the path for $g$, and let
$$\epsilon(a,b)~=~[\tau_a(S^n) \oplus \tau_b(S^n) :\tau_{g(a)}(M)] 
\in \{\pm 1\}$$
be $+1$ (resp. $-1$) if the isomorphism
$$(dg~dh)~:\tau_a(S^n) \oplus \tau_b(S^n)~\isa~\tau_{g(a)}(M)$$
is orientation-preserving (resp. reversing).
The {\it geometric intersection} of $x,y \in K_n(M)$ is given by
$$\lambda(x,y)~=~\sum\limits_{(a,b) \in D(g,h)}I(a,b) \in \Lambda$$ 
with
$$I(a,b)~=~\epsilon(a,b)\gamma(a,b) \in \Lambda~.$$
It follows from
$$\eqalign{
&\epsilon(b,a)~=~
[\tau_bS^n \oplus \tau_aS^n :\tau_aS^n \oplus \tau_bS^n]\epsilon(a,b)\cr
&\hphantom{\epsilon(b,a)~}=~
{\rm det}(\pmatrix{0 & 1 \cr 1 & 0}:\bbr^n \oplus \bbr^n \to
\bbr^n \oplus \bbr^n)\epsilon(a,b)\cr
&\hphantom{\epsilon(b,a)~}=~(-1)^n\epsilon(a,b) \in \{\pm 1\}~,\cr
&\gamma(b,a)~=~w(X)(\gamma(a,b))\gamma(a,b)^{-1} \in \pi_1(X)~,\cr
&I(b,a)=(-1)^n\overline{I(a,b)} \in \Lambda}$$
that 
$$\lambda(y,x)~=~(-1)^n\overline{\lambda(x,y)} \in \Lambda$$
(which also holds from purely homological considerations).\hfil\break
(ii) The quadratic function $\mu:K_n(M) \to Q_{(-1)^n}(\Lambda)$ is 
defined by geometric self-intersection numbers, as follows.
Represent $x \in K_n(M)$ by an immersion $g:S^n \looparrowright M$ 
as in (i), with transverse self-intersections. The double point set of $g$
$$\eqalign{D_2(g)~&=~D(g,g)\backslash \Delta(S^n)\cr
&=~\{(a,b) \in S^n \times S^n\,\vert\,a\neq b\in S^n,g(a)=g(b)\in M\}}$$
is finite, with a free $\bbz_2$-action $(a,b) \mapsto (b,a)$. 
For each $(a,b) \in D_2(g)$ let $\gamma(a,b)$ be the loop in
$M$ obtained by transporting to the base point the image under $g$ of a
path in $S^n$ from $a$ to $b$.  The {\it geometric self-intersection} 
of $x$ is defined by
$$\mu (x)~=~\sum\limits_{(a,b) \in D_2(g)/\bbz_2}I(a,b) \in Q_{(-1)^n}(\Lambda)~,$$ 
with $I(a,b)=\epsilon(a,b)\gamma(a,b)$ as in (i). Note that $\mu(x)$ is independent
of the choice of ordering of $(a,b)$ since 
$I(b,a)~=~(-1)^n\overline{I(a,b)} \in \Lambda$.\hfil\break
(iii) The kernel $(-1)^n$-quadratic form $(K_n(M),\lambda,\mu)$ is such 
that $\mu(x)=0$ if (and for $n \geq 3$ only if) $x \in K_n(M)$ can be killed 
by surgery on $S^n \subset M^{2n}$, i.e. represented by an embedding
$S^n \times D^n \hookrightarrow M$ with a nullhomotopy in $X$ --
the condition $\mu(x)=0$ allows the double points of a representative
framed immersion $g:S^n \looparrowright M$ to be matched in pairs, 
which for $n \geq 3 $ can be cancelled by the Whitney trick. 
The effect of the surgery is a bordant $(n-1)$-connected  normal map
$$(f',b')~:~M'^{2n}~=~{\rm cl.}(M\backslash S^n \times D^n) \cup D^{n+1} \times S^{n-1}
 \to X$$
with kernel $\Lambda$-modules
$$K_i(M')~=~\cases{
{\rm coker}(x^*\lambda:K_n(M)\to \Lambda^*)&if $i=n-1$\ccr
{\displaystyle{{\rm ker}(x^*\lambda:K_n(M)\to \Lambda^*)}\over 
\displaystyle{{\rm im}(x:\Lambda \to K_n(M))}}&if $i=n$\ccr
{\rm ker}(x:\Lambda \to K_n(M))&if $i=n+1$\ccr
0&otherwise~.}$$
Thus $(f',b')$ is $n$-connected if and only if $x$ generates a
direct summand $L=\langle x\rangle \subset K_n(M)$, in which case
$L$ is a sublagrangian of $(K_n(M),\lambda,\mu)$ in the terminology of \S5, 
with
$$\eqalign{
&L \subseteq L^{\perp}~=~\{y \in K_n(M)\,\vert\,\lambda(x,y)=0\}~,\cr
&(K_n(M'),\lambda',\mu')~=~(L^{\perp}/L,[\lambda],[\mu])~,\cr
&(K_n(M),\lambda,\mu)~\cong~(K_n(M'),\lambda',\mu')\oplus H_{(-1)^n}(\Lambda)~.}$$
(iv) The effect on $(f,b)$ of a surgery on an $X$-nullhomotopic 
embedding $S^{n-1} \times D^{n+1} \hookrightarrow M$ is an
$n$-connected bordant normal map
$$(f'',b'')~:~M''^{2n}~=~{\rm cl.}(M\backslash S^{n-1} \times D^{n+1}) 
\cup D^n \times S^n~=~M \# (S^n \times S^n)  \to X$$
with kernel $\Lambda$-modules
$$K_i(M'')~=~\cases{K_n(M) \oplus \Lambda \oplus \Lambda^*&if $i=n$\ccr
0&otherwise}$$
and kernel $(-1)^n$-quadratic form
$$(K_n(M''),\lambda'',\mu'')~=~(K_n(M),\lambda,\mu)\oplus H_{(-1)^n}(\Lambda)~.$$
(v) The main result of even-dimensional surgery obstruction theory is
that for $n \ge 3$ an $n$-connected $2n$-dimensional normal map
$(f,b):M^{2n} \to X$ is normal bordant to a homotopy equivalence if
and only if there exists an isomorphism of $(-1)^n$-quadratic forms
over $\Lambda=\bbz[\pi_1(X)]$ of the type
$$(K_n(M),\lambda,\mu) \oplus H_{(-1)^n}(\Lambda^k)~
\isa~ H_{(-1)^n}(\Lambda^{k'})$$
for some $k,k' \ge 0$.\hfill$\square$

\noindent {\sc Example 2.16}
There is also a relative version of 2.15.
An $n$-connected $2n$-dimensional normal map of pairs
$(f,b):(M^{2n},\partial M) \to (X, \partial X)$
has a kernel $(-1)^n$-quadratic form $(K_n(M),\lambda,\mu)$
over $\bbz[\pi_1(X)]$ is nonsingular if and only if
$\partial f:\partial M \to \partial X$ is a homotopy equivalence
(assuming $\pi_1(\partial X)\cong \pi_1(X)$).\hfill$\square$

\noindent {\sc Remark 2.17} 
(Realization of even-dimensional surgery obstructions, Wall [28,\t 5.8])\hfil\break 
(i) Let $X^{2n-1}$ be a $(2n-1)$-dimensional manifold, and suppose
given an embedding $e:S^{n-1} \times D^n \hookrightarrow X$, together
with a null-homotopy $\delta e$ of $e\vert : S^{n-1}\hookrightarrow X$
and a null-homotopy of the map $S^{n-1} \to O$ comparing the
(stable) trivializations of $\nu_{e\vert}:S^{n-1} \to BO(n)$
given by $e$ and $\delta e$. Then there is defined an $n$-connected
$2n$-dimensional normal map
$$(f,b)~:~(M;\partial_-M,\partial_+M) \to X \times ([ 0,1];\{0\},\{1\})$$
with 
$$\eqalign{&\partial_-f~=~{\rm id.}~:~\partial_-M~=~X \to X~,\cr
&M^{2n}~=~X \times [ 0,1] \cup_e D^n \times D^n~,\cr
&\partial_+M~=~{\rm cl.}(X \backslash e(S^{n-1} \times D^n)) \cup D^n \times S^{n-1}~.}$$
The kernel $(-1)^n$-quadratic form $(\Lambda,\lambda,\mu)$ over 
$\Lambda$ (2.16) is the (self-)\-intersection of the
framed immersion $S^{n-1} \times [ 0,1] \looparrowright X \times [ 0,1]$
defined by the track of a regular homotopy 
$e_0 \simeq e:S^{n-1}\times D^n \to X$
from a trivial unlinked embedding 
$$e_0~:~S^{n-1}\times D^n \hookrightarrow S^{2n-1}~=~S^{n-1} \times D^n \cup
D^n \times S^{n-1} \hookrightarrow X \# S^{2n-1}~=~X~.$$
Moreover, every form $(\Lambda,\lambda,\mu)$
arises in this way\t : starting with $e_0$ construct a regular homotopy
$e_0 \simeq e$ to a (self-)linked embedding $e$ such that the track has
(self-)intersection $(\lambda,\mu)$. \hfil\break
(ii) Let $(K,\lambda,\mu)$ be a $(-1)^n$-quadratic form over $\bbz[\pi]$, with
$\pi$ a finitely presented group and $K=\bbz[\pi]^k$ f.\ g. free.
Let $n \geq 3$, so that there exists a $(2n-1)$-dimensional manifold $X^{2n-1}$ 
with $\pi_1(X)=\pi$. For any such $n \geq 3$, $X$ there exists an $n$-connected 
$2n$-dimensional normal map 
$$(f,b)~:~(M^{2n};\partial_-M,\partial_+M) \to X^{2n-1} \times ([ 0,1];\{0\},\{1\})$$
with kernel form $(K,\lambda,\mu)$ and
$$\eqalign{&\partial_-f~=~{\rm id.}~:~\partial_-M~=~X\to X~,~K_n(M)~=~K~,\cr
&K_{n-1}(\partial_+M)~=~{\rm coker}(\lambda:K \to K^*)~,\cr
&K_n(\partial_+M)~=~{\rm ker}(\lambda:K \to K^*)~.}$$
The map $\partial_+f:\partial_+M \to X$ is a homotopy equivalence if and only
if the form $(K,\lambda,\mu)$ is nonsingular. Given $(K,\lambda,\mu)$, $X$ the
construction of $(f,b)$ proceeds as in (i).\hfill$\square$

\noindent {\sc Example 2.18} For $\pi_1(X)=\{1\}$ the realization of
even-dimensional surgery obstructions (2.17) is essentially the same as
the Milnor [11], [12] construction of $(n-1)$-connected $2n$-dimensional
manifolds by {\it plumbing} together $n$-plane bundles over $S^n$.  Let
$G$ be a finite connected graph without loops (= edges joining a vertex
to itself), with vertices $v_1,v_2,\dots,v_k$. Suppose given an oriented 
$n$-plane bundle over $S^n$ at each vertex 
$$\omega_1,\omega_2,\dots,\omega_k\in \pi_n(BSO(n))~=~\pi_{n-1}(SO(n))~,$$
regarded as a weight. Let $(\bbz^k,\lambda)$ be the $(-1)^n$-symmetric form 
over $\bbz$ defined by the $(-1)^n$-symmetrized adjacency matrix of $G$
and the Euler numbers $\chi(\omega_i)\in \bbz$, with
$$\eqalign{
&\lambda_{ij}~=~\cases{
\hbox{\rm no. of edges in $G$ joining $v_i$ to $v_j$}&if $i<j$\cr
(-1)^n\hbox{\rm (no. of edges in $G$ joining $v_i$ to $v_j$)}&if $i>j$\cr
\chi(\omega_i)&if $i=j$~,}\cr
&\lambda~:~\bbz^k \times \bbz^k \to \bbz~;~
((x_1,x_2,\dots,x_k),(y_1,y_2,\dots,y_k)) \mapsto
\sum\limits^k_{i=1}\sum\limits^k_{j=1}\lambda_{ij}x_iy_j~.}$$
The graph $G$ and the Euler numbers $\chi(\omega_i)$ 
determine and are determined by the form $(\bbz^k,\lambda)$.\hfil\break
(i) See Browder [1,\t Chapter V] for a detailed account of the
plumbing construction which uses $G$ to glue together the 
$(D^n,S^{n-1})$-bundles
$$(D^n,S^{n-1}) \to (E(\omega_i),S(\omega_i)) \to S^n~~(i=1,2,\dots,k)$$
to obtain a connected $2n$-dimensional manifold with boundary 
$$(P,\partial P)~=~(P(G,\omega),\partial P(G,\omega))$$
such that $P$ is an identification space
$$P~=~(\coprod\limits^k_{i=1}E(\omega_i))/\sim$$
with 1-skeleton homotopy equivalent to $G$, fundamental group
$$\pi_1(P)~=~\pi_1(G)~=~\ast_g\bbz$$
the free group on $g=1-\chi(G)$ generators, homology
$$H_r(P)~=~\cases{\bbz&if $r=0$\cr
\bbz^g&if $r=1$\cr
\bbz^k&if $r=n$\cr
0&otherwise~,}$$
and intersection form $(H_n(P),\lambda)$. Killing $\pi_1(P)$
by surgeries removing $g$ embeddings $S^1 \times D^{2n-1} \subset P$
representing the generators, there is obtained an
$(n-1)$-connected $2n$-dimensional manifold with boundary
$$(M,\partial M)~=~(M(G,\omega),\partial M(G,\omega))$$
such that
$$\eqalign{&H_r(M)~=~\cases{\bbz&if $r=0$\cr
\bbz^k&if $r=n$\cr
0&otherwise~,}\cr
&\lambda~:~H_n(M) \times H_n(M) \to \bbz~;\cr
&\hskip25pt ((x_1,x_2,\dots,x_k),(y_1,y_2,\dots,y_k)) \mapsto 
\sum\limits^k_{i=1}\sum\limits^k_{j=1}\lambda_{ij}x_iy_j~,\cr
&\tau_M~\simeq~\bigvee^k_{i=1}(\omega_i\oplus \epsilon^n)~:~
M~\simeq~ \bigvee^k_{i=1}S^n \to BSO(2n)~,\cr
&H_r(\partial M)~=~\cases{\bbz&if $r=0,2n-1$\cr
{\rm coker}(\lambda:\bbz^k \to \bbz^k)&if $r=n-1$\cr
{\rm ker}(\lambda:\bbz^k \to \bbz^k)&if $r=n$\cr
0&otherwise~.}}$$
If $G$ is a tree then $g=0$, $\pi_1(P(G,\omega))=\{1\}$, and
$$(M(G,\omega),\partial M(G,\omega))~=~(P(G,\omega),\partial P(G,\omega))~.$$
(ii) By Wall [26] for $n \geq 3$ an integral $(-1)^n$-symmetric matrix 
$(\lambda_{ij})_{1 \leq i,j \leq k}$ and 
elements $\omega_1,\omega_2,\dots,\omega_k\in \pi_n(BSO(n))$ with
$$\lambda_{ii}~=~\chi(\omega_i) \in \bbz~~(i=1,2,\dots,k)$$
determine an embedding
$$x~=~\bigcup_kx_i~:~\bigcup_kS^{n-1} \times D^n \hookrightarrow S^{2n-1}$$
such that\t :
{\parindent=23pt
\parskip=2pt
\item{(a)} for $1 \leq i < j \leq k$
$$\hbox{\rm linking number}(x_i(S^{n-1}\times 0) \cap x_j(S^{n-1}\times 0)
\hookrightarrow S^{2n-1})~=~\lambda_{ij} \in \bbz~,$$
\item{(b)} for $1 \leq i \leq k$ 
$x_i:S^{n-1}\times D^n \hookrightarrow S^{2n-1}$ is isotopic to the embedding}
$$\eqalign{&
e_{\omega_i}~:~S^{n-1} \times D^n \hookrightarrow S^{2n-1}~=~
S^{n-1} \times D^n\cup D^n\times S^{n-1}~;\cr
&\hskip100pt (s,t) \mapsto (s,\omega_i(s)(t))~.}$$
\noindent
Using $x$ to attach $k$ $n$-handles to $D^{2n}$ there is obtained
an oriented $(n-1)$-connected $2n$-dimensional manifold 
$$M(G,\omega)~=~D^{2n} \cup_x \bigcup_k \hbox{\rm $n$-handles}~ 
D^n \times D^n$$
with boundary an oriented $(n-2)$-connected $(2n-1)$-dimensional manifold 
$$\partial M(G,\omega)~=~{\rm cl.}(S^{2n-1}\backslash 
x(\bigcup_k S^{n-1} \times D^n)) \cup\bigcup_k D^n \times S^{n-1}~.$$
Moreover, every oriented $(n-1)$-connected $2n$-dimensional manifold 
with non-empty $(n-2)$-connected boundary is of the form 
$(M(G,\omega),\partial M(G,\omega))$, 
with $(\lambda_{ij},\omega_i)$ the complete set of diffeomorphism invariants.\hfil\break
(iii) Stably trivialized $n$-plane bundles over $S^n$ are classified
by $Q_{(-1)^n}(\bbz)$, with an isomorphism
$$Q_{(-1)^n}(\bbz)~\isa~\pi_{n+1}(BSO,BSO(n)) ~;~1 \mapsto 
(\delta\tau_{S^n},\tau_{S^n})$$
with 
$$\delta\tau_{S^n}~:~\tau_{S^n} \oplus \epsilon~ \cong~\epsilon^{n+1}$$
the stable trivialization given by the standard embedding $S^n \subset S^{n+1}$.
The map 
$$Q_{(-1)^n}(\bbz)~=~\pi_{n+1}(BSO,BSO(n)) \to \pi_n(BSO(n))~;~
1 \mapsto \tau_{S^n}$$
is an injection for $n \neq 1,3,7$. With $G$ as above, suppose now
that the vertices $v_1,v_2,\dots,v_k$ are weighted by elements
$$\mu_1,\mu_2,\dots,\mu_k\in \pi_{n+1}(BSO,BSO(n))~=~Q_{(-1)^n}(\bbz)~.$$
Define
$$\eqalign{\omega_i~=~[\mu_i] \in~ 
&{\rm im}(\pi_{n+1}(BSO,BSO(n)) \to \pi_n(BSO(n)))\cr
&=~{\rm ker}(\pi_n(BSO(n)) \to \pi_n(BSO))}$$
and let $(\bbz^k,\lambda,\mu)$ be the $(-1)^n$-quadratic form over $\bbz$ 
with $\lambda$ as before and
$$\mu~:~\bbz^k \to Q_{(-1)^n}\bbz~;~(x_1,x_2,\dots,x_k) \mapsto 
\sum\limits_{1 \leq i< j \leq k}\lambda_{ij}x_ix_j+
\sum\limits^k_{i=1}\mu_i(x_i)^2~,$$
such that
$$\lambda_{ii}~=~\chi(\omega_i)~=~(1+(-1)^n)\mu_i \in \bbz~.$$
The $(n-1)$-connected $2n$-dimensional manifold 
$$M(G,\mu_1,\mu_2,\dots,\mu_n)~=~M(G,\omega_1,\omega_2,\dots,\omega_n)$$
is stably parallelizable, with an $n$-connected normal map
$$(M(G,\mu_1,\mu_2,\dots,\mu_n),\partial M(G,\mu_1,\mu_2,\dots,\mu_n))
 \to (D^{2n},S^{2n-1})$$
with kernel form $(\bbz^k,\lambda,\mu)$.\hfil\break
(iv) For $n \geq 3$ the realization of a $(-1)^n$-quadratic form 
$(\bbz^k,\lambda,\mu)$ over $\bbz$ (2.17) 
is an $n$-connected $2n$-dimensional normal map
$$\eqalign{(f,b)~:~&(M^{2n};S^{2n-1},\partial_+M)\cr
&=~({\rm cl.}(M(G,\mu_1,\dots,\mu_k)\backslash D^{2n});S^{2n-1},\partial
M(G,\mu_1,\dots,\mu_k))\cr
& \hskip100pt \to S^{2n-1} \times ([ 0,1];\{0\},\{1\})~.}$$
with kernel form $(\bbz^k,\lambda,\mu)$. 
If $(\bbz^k,\lambda,\mu)$ is nonsingular then 
$$\partial_+f~:~\Sigma^{2n-1}~=~\partial M(G,\mu_1,\dots,\mu_k) \to S^{2n-1}$$ 
is a homotopy equivalence, and $\Sigma^{2n-1}$ is a homotopy 
sphere with a potentially exotic differentiable structure (Milnor [10],\t 
Kervaire and Milnor [7]) -- see 2.20, 3.6 and 3.7 below.\hfill$\square$

\noindent{\sc Example 2.19}
(i) Consider the special case $k=1$ of 2.18 (i).
Here $G=\{v_1\}$ is the graph with one vertex, and 
$$\omega \in \pi_n(BSO(n))~=~\pi_{n-1}(SO(n))$$
classifies an $n$-plane bundle over $S^n$. The plumbed
$(n-1)$-connected $2n$-dimensional manifold with boundary is
the $(D^n,S^{n-1})$-bundle over $S^n$ 
$$(M(G,\omega),\partial M(G,\omega))~=~(E(\omega),S(\omega))$$
with
$$\eqalign{E(\omega)~&=~S^{n-1} \times D^n \cup_{
(x,y) \sim (x,\omega(x)(y))}S^{n-1} \times D^n\cr
&=~D^{2n} \cup_{e_{\omega}}D^n \times D^n}$$
obtained from $D^{2n}$ by attaching an $n$-handle along the embedding
$$\eqalign{&e_{\omega}~:~S^{n-1} \times D^n \hookrightarrow S^{2n-1}~=~
S^{n-1} \times D^n\cup D^n\times S^{n-1}~;\cr
&\hskip100pt (x,y) \mapsto (x,\omega(x)(y))~.}$$
(ii) Consider the special case $k=1$ of 2.18 (iii), the realization of a
$(-1)^n$-quadratic form $(\bbz,\lambda,\mu)$ over $\bbz$,
with $G=\{v_1\}$ as in (i). An element
$$\eqalign{\mu~=~(\delta\omega,\omega) \in \pi_{n+1}(BSO,BSO(n))~
&=~Q_{(-1)^n}(\bbz)\cr
&=~\cases{\bbz&if $n \equiv 0(\bmod\,2)$\cr
\bbz_2&if $n \equiv 1(\bmod\,2)$}}$$
classifies an $n$-plane bundle $\omega:S^n \to BSO(n)$
with a stable trivialization 
$$\delta\omega~:~\omega \oplus \epsilon^{\infty}~\cong~\epsilon^{n+\infty}~,$$ 
and
$$(M(G,\mu),\partial M(G,\mu))~=~(E(\omega),S(\omega))~.$$
For $n \neq 1,3,7$ $\delta\omega$ is determined by $\omega$.
For even $n$ $\mu \in Q_{+1}(\bbz)=\bbz$ and
$$\omega~=~\mu^*\tau_{S^n}~:~
S^n \xymatrix@C+10pt{\ar[r]^-{\displaystyle{\mu}}&} S^n 
\xymatrix@C+10pt{\ar[r]^-{\displaystyle{\tau_{S^n}}}&} BSO(n)$$
is the unique stably trivial $n$-plane bundle over $S^n$ with
Euler number
$$\chi(\omega)~=~2\mu\in \bbz~.$$
For odd $n\neq 1,3,7$ $\mu \in Q_{-1}(\bbz)=\bbz_2$ and
$$\omega~=~\cases{\tau_{S^n}&if $\mu=1$\cr
\epsilon^n&if $\mu=0$~.}$$
For $n=1,3,7$ 
$$\omega~=~\tau_{S^n}~=~\epsilon^n~:~S^n \to BSO(n)$$
and $\delta\omega$ is the (stable) trivialization of $\omega$ with 
mod 2 Hopf invariant $\mu$.\hfil\break
The plumbed $(n-1)$-connected $2n$-dimensional manifold 
$$(M(G,\mu),\partial M(G,\mu))~=~(M(G,\omega),\partial M(G,\omega))~=~(E(\omega),S(\omega))$$
(as in (i)) is stably parallelizable. The trace of the surgery on the normal map
$$(f_-,b_-)~=~{\rm id.}~:~\partial_-M_{\omega}~=~S^{2n-1} \to S^{2n-1}$$ 
killing $e_{\omega}:S^{n-1} \times D^n \hookrightarrow S^{2n-1}$
is an $n$-connected $2n$-dimensional normal map
$$(f_{\omega},b_{\delta\omega})~:~
(M^{2n}_{\omega};\partial_-M_{\omega},\partial_+M_{\omega}) 
\to S^{2n-1} \times ([ 0,1];\{0\},\{1\})$$
with
$$\eqalign{
&M_{\omega}~=~{\rm cl.}(M(G,\mu)\backslash D^{2n})~=~
{\rm cl.}(E(\omega)\backslash D^{2n})\cr
&\hphantom{M_{\omega}~}=~
S^{2n-1}\times [ 0,1]\cup_{e_{\omega}}D^n \times D^n~,\cr
&\partial_+M_{\omega}~=~{\rm cl.}(S^{2n-1} \backslash 
e_{\omega}(S^{n-1} \times D^n)) \cup D^n\times S^{n-1}\cr
&\hphantom{\partial_+M_{\omega}~}=~
D^n\times S^{n-1} \cup_{\omega} D^n\times S^{n-1}~=~S(\omega)~,\cr
&K_n(M_{\omega})~=~\bbz}$$
and kernel form $(\bbz,\lambda,\mu)$.
If $\mu=0 \in Q_{(-1)^n}(\bbz)$ then 
$$\omega~=~\epsilon^n~:~S^n \to BSO(n)~~,~~
\partial_+M_{\omega}~=~S(\epsilon^n)~=~S^{n-1} \times S^n~,$$
If $\mu=1 \in Q_{(-1)^n}(\bbz)$ then 
$$\omega~=~\tau_{S^n}~:~S^n \to BSO(n)~~,~~
\partial_+M_{\omega}~=~S(\tau_{S^n})~=~O(n+1)/O(n-1)~.$$ 
(iii) Consider the special case $k=2$ of 2.18 (i), with $G=I$ the
graph with 1 edge and 2 vertices 
$$
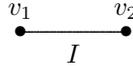
\beginpicture
\setcoordinatesystem units <4pt,4pt> 
\setlinear
\plot 0 0 10 0 /
\put{$\bullet$} at 0 0
\put{$\bullet$} at 10 0
\put{~} at 0 -4
\put{~} at 0 4
\put{$v_1$} at 0 2
\put{$v_2$} at 10 2
\put{$I$} at 5 -2
\endpicture$$
For any weights $\omega_1,\omega_2 \in \pi_n(BSO(n))$ there is obtained an
$(n-1)$-connected $2n$-dimensional manifold
$$M(I,\omega_1,\omega_2)~=~D^{2n} \cup_{e_{\omega_1} \cup e_{\omega_2}}
(D^n \times D^n \cup D^n \times D^n)$$
by plumbing as in Milnor [11],\t [12], with intersection form the
$(-1)^n$-symmetric form $(\bbz \oplus \bbz,\lambda)$ over $\bbz$ defined by
$$\eqalign{
\lambda~:~&\bbz \oplus \bbz \times \bbz \oplus \bbz \to \bbz~;\cr
&((x_1,x_2),(y_1,y_2)) \mapsto \chi(\omega_1)x_1y_1 +\chi(\omega_2)x_2y_2 +
x_1y_2 + (-1)^n x_2y_1~.}$$
(iv) Consider the special case $k=2$ of 2.18 (iii), with $G=I$ as in (iii).
For $\mu_1,\mu_2 \in Q_{(-1)^n}(\bbz)$ and
$$\omega_i~=~[\mu_i] \in {\rm im}(Q_{(-1)^n}(\bbz) \to \pi_n(BSO(n)))$$
the $(-1)^n$-quadratic form $(\bbz \oplus \bbz,\lambda,\mu)$ over $\bbz$
defined by
$$\mu~:~\bbz \oplus \bbz \to Q_{(-1)^n}(\bbz)~;~(x_1,x_2)
\mapsto \mu_1(x_1)^2 +\mu_2(x_2)^2+x_1x_2$$
is the kernel form of an $n$-connected $2n$-dimensional normal map
$$(f,b)~:~M(I,\mu_1,\mu_2)~=~M(I,\omega_1,\omega_2) \to D^{2n}~.$$
If $\mu_1=\mu_2=0$ then $(\bbz\oplus\bbz,\lambda,\mu)=H_{(-1)^n}(\bbz)$ 
is hyperbolic $(-1)^n$-quadratic form over $\bbz$, with
$$\eqalign{&\lambda~:~\bbz \oplus \bbz \times \bbz \oplus \bbz \to \bbz~;~
((x_1,x_2),(y_1,y_2)) \mapsto x_1y_2+(-1)^nx_2y_1~,\cr
&\mu~:~\bbz \oplus \bbz \to Q_{(-1)^n}(\bbz)~=~\bbz/\{1+(-1)^{n-1}\}~;~
(x_1,x_2) \mapsto x_1x_2~,}$$
and the plumbed manifold is a punctured torus
$$(M(I,0,0)^{2n},\partial M(I,0,0))~=~
({\rm cl.}(S^n \times S^n \backslash D^{2n}),S^{2n-1})~.$$
The hyperbolic form is the kernel of the $n$-connected $2n$-dimensional normal map
$$\eqalign{&(f,b)~:~(M;\partial_-M,\partial_+M)~=~
({\rm cl.}(M(I,0,0)\backslash D^{2n});S^{2n-1},S^{2n-1})\cr
&\hskip150pt \to S^{2n-1} \times ([ 0,1];\{0\},\{1\})}$$
defined by the trace of surgeries on the linked spheres
$$S^{n-1} \cup S^{n-1} \hookrightarrow S^{2n-1}~=~S^{n-1} \times D^n
\cup D^n \times S^{n-1}$$
with no self-linking. These are 
the attaching maps for the cores of the $n$-handles in the decomposition
$$M(I,0,0)~=~D^{2n} \cup D^n \times D^n \cup D^n \times D^n~,$$
using the standard framings of $S^{n-1} \subset S^{2n-1}$.
If $n$ is odd, say $n=2k+1$, and $\mu_0=\mu_1=1 \in Q_{-1}(\bbz)$
 the form in 2.18 (i) is just the Arf $(-1)$-quadratic form over $\bbz$
$(\bbz\oplus\bbz,\lambda,\mu')$ with
$$\mu'~:~\bbz \oplus \bbz \to Q_{(-1)}(\bbz)~=~\bbz_2~;~
(x,y) \mapto x^2+xy+y^2~.$$
The plumbed manifold 
$$M(I,1,1)^{4k+2}~=~D^{4k+2} \cup D^{2k+1} \times D^{2k+1} \cup D^{2k+1} \times D^{2k+1}$$
has the same attaching maps for the cores of the $(2k+1)$-handles as
$M(I,0,0)$, but now
using the framings of $S^{2k} \subset S^{4k+1}$ classified by
$$\tau_{S^{2k+1}} \in \pi_{2k+1}(BSO(2k+1))~=~\pi_{2k}(SO(2k+1))$$
(which is zero if and only if $k=0,1,3$).  
The Arf form is the kernel of the $(2k+1)$-connected $(4k+2)$-dimensional normal map
$$\eqalign{&(f',b')~:~(M';\partial_-M',\partial_+M')~=~
({\rm cl.}(M(I,1,1)\backslash D^{4k+2});S^{4k+1},\Sigma^{4k+1})\cr
&\hskip150pt \to S^{4k+1} \times ([ 0,1];\{0\},\{1\})}$$
defined by the trace of surgeries on the linked spheres
$$S^{2k} \cup S^{2k} \hookrightarrow S^{4k+1}~=~S^{2k} \times D^{2k+1}
\cup D^{2k+1} \times S^{2k}~,$$
with self-linking given by the non-standard framing.
(See 3.7 below for a brief account of the exotic sphere $\Sigma^{4k+1}$).
\hfill$\square$

\noindent{\sc Example 2.20} 
The sphere bundles $S(\omega)$ of certain oriented 4-plane bundles
$\omega$ over $S^4$ (the special case $n=4$ of 2.19 (i)) give explicit
exotic 7-spheres.  An oriented 4-plane bundle $\omega:S^4 \to BSO(4)$
is determined by the Euler number and first Pontrjagin class
$$\chi(\omega),p_1(\omega) \in H^4(S^4)~=~\bbz~,$$
which must be such that 
$$2\chi(\omega)\equiv p_1(\omega)(\bmod\,4)~,$$ 
with an isomorphism
$$\pi_4(BSO(4))~\isa~\bbz \oplus \bbz~;~\omega \mapto 
((2\chi(\omega)+p_1(\omega))/4,(2\chi(\omega)-p_1(\omega))/4)~.$$
If $\chi(\omega)=1$ then $S(\omega)$ is a homotopy 7-sphere, and
$$p_1(\omega)~=~2\ell\in H^4(S^4)~=~\bbz$$ 
for some odd integer $\ell$. The 7-dimensional differentiable manifold
$\Sigma_\ell^7=S(\omega)$ is homeomorphic to $S^7$ (by Smale's generalized
Poincar\'e conjecture, or by a direct Morse-theoretic argument). If
$\Sigma_\ell^7$ is diffeomorphic to $S^7$ then 
$$M^8~=~E(\omega) \cup_{\Sigma_\ell^7}D^8$$
is a closed 8-dimensional differentiable manifold with 
$$p_1(M)~=~p_1(\omega)~=~2\ell \in H^4(M)~=~\bbz~~,~~\sigma(M)~=~1 \in \bbz~.$$
By the Hirzebruch signature theorem
$$\eqalign{\sigma(M)~&=~\langle {\cal L}(M),[M] \rangle\cr
&=~(7p_2(M)-p_1(M)^2)/45\cr
&=~(7p_2(M)-4\ell^2)/45~=~1 \in \bbz~.}$$
If $\ell \not \equiv \pm 1(\bmod\,7)$ then 
$$p_2(M)~=~(45+4\ell^2)/7 \not \in H^8(M)~=~\bbz$$
so that there is no such diffeomorphism, and $\Sigma_\ell^7$
is an exotic 7-sphere (Milnor [10], Milnor and Stasheff [14,\t p.247]).
\hfill$\square$

\noindent {\sc Definition 2.21}
An {\it isomorphism} of $\epsilon$-symmetric forms 
$$f~:~(K,\lambda)~\isa~(K',\lambda')$$ 
is a $\Lambda$-module isomorphism $f:K\cong K'$ such that
$$\lambda'(f(x),f(y))~=~\lambda(x,y) \in \Lambda~.$$
An {\it isomorphism} of $\epsilon$-quadratic forms
$f:(K,\lambda ,\mu)\cong (K',\lambda',\mu')$ is an isomorphism of the
underlying $\epsilon$-symmetric forms 
$f:(K,\lambda)\cong (K',\lambda')$ such that 
$$\mu'(f(x)) ~=~\mu (x) \in Q_{\epsilon}(\Lambda)~.\eqno{\square}$$

\noindent {\sc Proposition 2.22}
{\it If there exists a central element $s\in \Lambda$ such that
$$s+{\overline s}=1 \in \Lambda$$ 
there is an identification of categories}
$$\{\epsilon {\it \hbox{\rm -}quadratic~forms~over}~\Lambda\}~=~
\{\epsilon {\it \hbox{\rm -}symmetric~forms~ over}~\Lambda\}~.$$
\noindent {\sc Proof}\t :
The $\epsilon$-symmetrization map 
$1+T_{\epsilon}:Q_{\epsilon}(K)\to Q^{\epsilon}(K)$ is an isomorphism
 for any $\Lambda$-module $K$, with inverse
$$Q^{\epsilon}(K) \to Q_{\epsilon}(K)~;~\lambda \mapto 
((x,y) \mapto s\lambda(x,y))~.$$
For any $\epsilon$-quadratic form $(K,\lambda,\mu)$ the 
$\epsilon$-quadratic function $\mu $ is determined by the 
$\epsilon$-symmetric pairing $\lambda $, with
$$\mu (x)~=~s \lambda(x,x) \in Q_{\epsilon}(\Lambda)~.\eqno{\square}$$

\noindent {\sc Example 2.23}
If $2\in \Lambda$ is invertible
then 2.22 applies with $s=1/2\in \Lambda$.\hfill $\square$

For any $\epsilon$-symmetric form $(K,\lambda)$ and $x \in K$
$$\lambda(x,x) \in Q^{\epsilon}(\Lambda)~.$$

\noindent {\sc Definition 2.24}
An $\epsilon$-symmetric form $(K,\lambda)$ is {\it even} if for all
$x\in K$
$$\lambda(x,x) \in {\rm im}(1+T_{\epsilon}:Q_{\epsilon}(\Lambda)\to 
Q^{\epsilon}(\Lambda))~.\eqno{\square}$$

\noindent {\sc Proposition 2.25}
{\it Let $\epsilon = 1$ or $-1$. If the $\epsilon$-symmetrization map 
$$1+T_{\epsilon}~:~Q_{\epsilon}(\Lambda)\to Q^{\epsilon}(\Lambda)$$
is an injection there is an identification of categories}
$$\{\epsilon {\it \hbox{\rm -}quadratic~forms~over}~\Lambda\}~=~
\{{\it even}~ \epsilon {\it \hbox{-}symmetric~forms~ over}~\Lambda\}~.$$
{\sc Proof}\t :
Given an even $\epsilon$-symmetric form $(K,\lambda)$ over $\Lambda$ 
there is a unique function $\mu : K \to Q_{\epsilon}(\Lambda)$ such that
for all $x \in K$
$$(1+T_{\epsilon})( \mu (x)) ~=~ \lambda(x,x) \in Q^{\epsilon}(\Lambda)~,$$
which then automatically satisfies the conditions of 2.13
for $(K,\lambda,\mu)$ to be an $\epsilon$-quadratic form.
\hfill$\square$

\noindent {\sc Example 2.26}
The symmetrization map 
$$1+T~=~2~:~Q_{+1}(\bbz)~=~\bbz \to Q^{+1} (\bbz)~=~\bbz$$ 
is an injection, so that quadratic forms over $\bbz$
coincide with the even symmetric forms.
\hfill$\square$

\noindent {\sc Example 2.27}
The $(-1)$-symmetrization map 
$$1+T_-~:~Q_{-1}(\bbz)~=~\bbz_2 \to Q^{-1} (\bbz)~=~0$$ 
is not an injection, so that $(-1)$-quadratic forms over $\bbz$
have a richer structure than even $(-1)$-symmetric forms.
The hyperbolic $(-1)$-symmetric form $(K, \lambda)=H^{-1}(\bbz)$ over $\bbz$ 
$$K~=~\bbz \oplus\bbz~~,~~
\lambda ~:~ K \times K \to\bbz ~;~ ((a,b),(c,d)) \mapto ad-bc$$
admits two distinct $(-1)$-quadratic refinements
$(K,\lambda,\mu)$, $(K,\lambda,\mu')$, with
$$\eqalign{
&\mu ~:~ K \to Q_{-1}(\bbz)~=~\bbz/2 ~;~ (x,y) \mapto xy~,\ccr
&\mu' ~:~ K \to Q_{-1}(\bbz)~=~\bbz /2 ~;~ (x,y) \mapto x^2+xy+y^2~.}$$
See \S3 below for the definition of the Arf invariant, which distinguishes 
the hyperbolic $(-1)$-quadratic form $(K,\lambda,\mu)=H_{-1}(\bbz)$ from
the Arf form $(K,\lambda,\mu')$ (which already appeared in 2.19 (iv)).
\hfill$\square$

\noindent {\sectionfont \S3. The even-dimensional $L$-groups}

The even-dimensional surgery obstruction groups $L_{2n}(\Lambda)$ will
now be defined, using the following preliminary result.
 
\noindent {\sc Lemma 3.1}
{\it For any nonsingular $\epsilon$-quadratic form $(K,\lambda,\mu)$
there is defined an isomorphism
$$(K,\lambda,\mu) \oplus (K,-\lambda,-\mu)~\cong~H_{\epsilon}(K)~,$$
with $H_{\epsilon}(K)$ the hyperbolic $\epsilon$-quadratic form 
{\rm (}2.14{\rm)}.}\hfil\break
\noindent {\sc Proof}\t :
Let $L$ be a f.\ g.  projective $\Lambda$-module such that $K \oplus L$
is f.\ g.  free, with basis elements $\{x_1,x_2,\dots ,x_k\}$ say.  Let
$$\lambda_{ij}~=~(\lambda \oplus 0)(x_j,x_i) \in \Lambda~~(1 \le i < j \le k)$$
and choose representatives $\mu_i \in \Lambda$ of $\mu (x_i) \in
Q_{\epsilon}(\Lambda)$ $(1\le i\le k)$.  Define the $\Lambda$-module
morphism
$$\eqalign{\psi_{K \oplus L} ~
&:~K\oplus L \to (K\oplus L)^*~;\ccr
&\sum^k_{i=1}a_ix_i\mapto (\sum^k_{j=1}b_jx_j\mapto 
\sum^k_{i=1}b_i \mu_i \overline{a}_i + \sum_{1\le i<j\le k}
b_j \lambda_{ij}{\overline a}_i)~.}$$
The $\Lambda$-module morphism defined by
$$\psi~:~K~\raise4pt\hbox{${\rm inclusion}\atop \ra{3.5}$}~ 
K\oplus L~ \raise4pt\hbox{$\psi_{K \oplus L} \atop \ra{2.6} $}~
(K\oplus L)^*~=~K^* \oplus L^*~ \raise4pt\hbox{${\rm projection}
\atop \ra{3.5}$} ~K^*$$ 
is such that 
$$\eqalign{&\lambda ~=~ \psi+ \epsilon \psi^* ~:~ K \to K^* ,\ccr
&\mu(x)~ =~\psi(x,x) \in Q_{\epsilon}(\Lambda)~(x \in K)~.}$$
As $(K,\lambda,\mu)$ is nonsingular $\psi+\epsilon \psi^*:K \to K^*$
is an isomorphism. The $\Lambda$-module morphism defined by
$$\widetilde \psi~=~ (\psi+ \epsilon \psi^*)^{-1}
\psi(\psi+ \epsilon \psi^*)^{-1} ~:~K^* \to K$$
is such that
$$(\psi+ \epsilon \psi^*)^{-1} ~=~ 
\widetilde \psi+ \epsilon \widetilde \psi^*~:~K^* \to K~.$$
Define an isomorphism of $\epsilon$-quadratic forms
$$f~:~H_{\epsilon}(K)~\isa~(K \oplus K, \lambda \oplus -\lambda,
\mu \oplus -\mu)$$
by 
$$f~=~\pmatrix{1 &-\epsilon \widetilde \psi^* \ccr 1 & \widetilde \psi}~:~
K \oplus K^*  \to K \oplus K ~.\eqno{\square}$$

\noindent {\sc Definition 3.2}
The {\it $2n$-dimensional $L$-group} $L_{2n}(\Lambda)$ is the group of
equivalence classes of nonsingular $(-1)^n$-quadratic forms
$(K,\lambda,\mu)$ on stably \break
f.\ g. free $\Lambda$-modules, subject to the equivalence relation
$$\eqalign{
&(K,\lambda , \mu) \sim (K',\lambda' , \mu')\ccr
&\hbox{\rm if there exists an isomorphism of $(-1)^n$-quadratic forms}\ccr
&(K,\lambda , \mu) \oplus H_{(-1)^n}(\Lambda^k )~\isa~ 
(K',\lambda', \mu') \oplus H_{(-1)^n}(\Lambda^{k'})\ccr
&\hbox{\rm for some f.\ g. free $\Lambda$-modules $\Lambda^k,\Lambda^{k'}$.}}$$
Addition and inverses in $L_{2n}(\Lambda)$ are given by
$$\eqalign{
(K_1,\lambda_1 , \mu_1) + (K_2,\lambda_2 , \mu_2)~&=~
(K_1 \oplus K_2,\lambda_1 \oplus \lambda_2, \mu_1 \oplus \mu_2) ~,\ccr
-(K,\lambda , \mu)~&=~(K,- \lambda ,- \mu) \in L_{2n}(\Lambda)~.}
\eqno{\lower8pt\hbox{$\square$}} $$

The groups $L_{2n}(\Lambda)$ only depend on the residue 
$n(\bmod\,2)$, so that only two $L$-groups have actually been defined, 
$L_0(\Lambda)$ and $L_2(\Lambda)$. Note that 3.2 uses Lemma 3.1 to justify 
$(K,\lambda , \mu) \oplus (K,-\lambda ,-\mu) \sim 0$.

\noindent {\sc Remark 3.3}
The surgery obstruction of an $n$-connected $2n$-dimen\-sional 
normal map $(f,b):M^{2n}\to X$ is an element
$$\sigma_*(f,b)~=~(K_n(M),\lambda,\mu) \in L_{2n}(\bbz[\pi_1(X)])$$
such that $\sigma_*(f,b)=0$ if (and for $n\ge 3$ only if) $(f,b)$
is normal bordant to a homotopy equivalence.
\hfill$\square$

\noindent {\sc Example 3.4}
Let $M=M^2_g$ be the orientable 2-manifold (=\ surface) of genus $g$,
with degree 1 map $f:M \to S^2$.
A choice of framing of the stable normal bundle of an embedding
$M \hookrightarrow \bbr^3$ determines a 1-connected 2-dimensional normal map
$(f,b): M \to S^2$.
For a standard choice of framing (i.e. one which extends to a 3-manifold $N$
with $\partial N=M$) the kernel form and the surgery obstruction are given by 
$$\sigma_*(f,b)~=~H_{-1}(\bbz^g)~=~0~ \in L_2(\bbz)~$$
and $(f,b)$ is normal bordant to a homotopy equivalence, i.e. 
$M$ is framed null-cobordant.
\hfill$\square$

\noindent {\sc Example 3.5}
The even-dimensional $L$-groups of the ring $\Lambda=$~$\bbr$ 
of real numbers with the identity involution are given by
$$L_{2n}(\bbr)~=~\cases{\bbz&if $n$ is even\ccr 0 &
if $n$ is odd~.}$$
Since $1/2 \in \bbr$ there is no difference between symmetric 
and quadratic forms over $\bbr$.\hfill$\square$

The {\it signature} (alias index) of a 
nonsingular symmetric form $(K,\lambda)$ over $\bbr$ is defined by
$$\eqalign{&\sigma (K,\lambda)~=~\hbox{\rm
no. of positive eigenvalues of $\lambda$}\ccr
&\hskip100pt -
\hbox{\rm no. of negative eigenvalues of } \lambda \in\bbz~.}$$
\noindent Here, the symmetric form $\lambda \in Q^{+1}(K)$ is identified
with the symmetric $k\times k$ matrix $(\lambda(x_i,x_j)_{1 \le i,j \le
k}) \in M_{k,k}(\bbr)$ determined by any choice of basis $x_1,x_2,\dots,x_k$ 
for $K$.  By Sylvester's law of inertia the rank $k$ and the
signature $\sigma(K,\lambda)$ define a complete set of invariants for
the isomorphism classification of nonsingular symmetric forms
$(K,\lambda)$ over $\bbr$, meaning that two forms are isomorphic if and
only if they have the same rank and signature.  A nonsingular quadratic
form $(K,\psi)$ over $\bbr$ is isomorphic to a hyperbolic form if and
only if it has signature $0$.  Two such forms $(K,\lambda)$,
$(K',\lambda')$ are related by an isomorphism
$$(K,\lambda) \oplus H_+ (\bbr^m )~\isa~ (K',\lambda') \oplus 
H_+ (\bbr^{m'})$$
if and only if they have the same signature
$$\sigma (K,\lambda)~=~\sigma (K',\lambda') \in\bbz~.$$
Moreover, every integer is the signature of a form, 
since $1 \in\bbz$ is the signature of the nonsingular symmetric form
($\bbr$,1) with
$$1~:~\bbr \to \bbr^*~;~x\mapto ( y\mapto xy)$$
and for any nonsingular
symmetric forms $(K,\lambda),(K',\lambda')$ over $\bbr$
$$\eqalign{&\sigma((K,\lambda) \oplus (K',\lambda'))~=~
\sigma(K,\lambda)+\sigma(K',\lambda')~,\ccr
&\sigma(K,-\lambda)~=~-\sigma(K,\lambda) \in\bbz~.}$$
The isomorphism of 3.5 in the case $n\equiv 0({\rm mod}\,2)$ is defined by
$$L_0 (\bbr)~\isa~\bbz~;~(K,\lambda)\mapto \sigma (K,\lambda)~.$$
$L_2(\bbr)=0$ because every nonsingular $(-1)$-symmetric 
(alias symplectic) form over $\bbr$ admits is isomorphic to a hyperbolic form.

It is not possible to obtain a complete isomorphism classification of
nonsingular symmetric and quadratic forms over $\bbz$ -- see Chapter II
of Milnor and Husemoller [13] for the state of the art in 1973. 
Fortunately, it is much easier to decide if two forms become isomorphic
after adding hyperbolics then whether they are actually isomorphic. 
Define the signature of a nonsingular symmetric form $(K,\lambda)$ over
$\bbz$ to be the signature of the induced nonsingular symmetric form
over $\bbr$ $$\sigma (K,\lambda)~=~\sigma (\bbr \otimes K,1\otimes
\lambda) \in\bbz~.$$ It is a non-trivial theorem that two nonsingular
even symmetric forms $(K,\lambda)$, $(K',\lambda')$ are related by an
isomorphism
$$(K,\lambda) \oplus H_+ (\bbz^m )~\isa~ (K',\lambda') \oplus 
H_+ (\bbz^{m'})$$
if and only if they have the same signature
$$\sigma (K,\lambda)~=~\sigma (K',\lambda') \in\bbz~.$$
Moreover, not every integer arises as the signature of 
an even symmetric form, only those divisible by 8.
The Dynkin diagram of the exceptional Lie group $E_8$ is a tree
$$
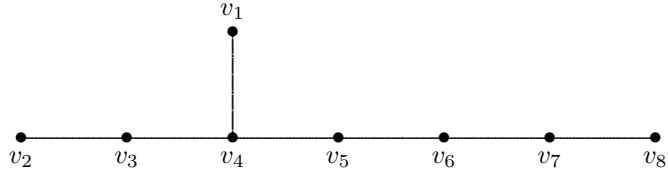
\beginpicture
\setcoordinatesystem units <4pt,4pt> 
\setlinear
\plot 0 0 10 0 20 0 30 0 40 0 50 0 60 0 /
\plot 20 0 20 10 /
\put{$\bullet$} at 0 0
\put{$\bullet$} at 10 0
\put{$\bullet$} at 20 0
\put{$\bullet$} at 30 0
\put{$\bullet$} at 40 0
\put{$\bullet$} at 50 0
\put{$\bullet$} at 60 0
\put{$\bullet$} at 20 10
\put{~} at 20 14
\put{~} at 20 -4
\put{$v_1$} at 20 12
\put{$v_2$} at 0 -2
\put{$v_3$} at 10 -2
\put{$v_4$} at 20 -2
\put{$v_5$} at 30 -2
\put{$v_6$} at 40 -2
\put{$v_7$} at 50 -2
\put{$v_8$} at 60 -2
\endpicture$$
\noindent
Weighing each vertex by $1 \in Q_{+1}(\bbz)=\bbz$
gives (by the method recalled in 2.18) a nonsingular quadratic form 
$(\bbz^8,\lambda_{E_8},\mu_{E_8})$ with signature
$$\sigma (\bbz^8,\lambda_{E_8})~=~8 \in\bbz~,$$
where
$$\lambda_{E_8}~=~\pmatrix{
2&0&0&1&0&0&0&0 \ccr
0&2&1&0&0&0&0&0 \ccr
0&1&2&1&0&0&0&0 \ccr
1&0&1&2&1&0&0&0 \ccr
0&0&0&1&2&1&0&0 \ccr
0&0&0&0&1&2&1&0 \ccr
0&0&0&0&0&1&2&1 \ccr
0&0&0&0&0&0&1&2 }~:~\bbz^8 \to (\bbz^8)^*$$
\smallskip
\noindent
and $\mu_{E_8}$ is determined by $\lambda_{E_8}$.

\noindent {\sc Example 3.6} 
(i) The signature divided by 8 defines an isomorphism
$$\sigma~:~L_{4k}(\bbz) \to \bbz~;~(K,\lambda,\mu) \mapsto \sigma(K,\lambda)/8$$
so that $(\bbz^8,\lambda_{E_8},\mu_{E_8}) \in L_{4k}(\bbz)$ represents a
generator.\hfil\break
(ii) See Kervaire and Milnor [7] and Levine [9] for the surgery
classification of exotic spheres in dimensions $n \geq 7$,
including the expression of the $h$-cobordism group $\Theta_n$ of 
$n$-dimensional exotic spheres as
$$\Theta_n~=~\pi_n(TOP/O)~=~\pi_n(PL/O)$$
and the exact sequence
$$\dots\to \pi_{n+1}(G/O) \to L_{n+1}(\bbz) \to \Theta_n \to \pi_n(G/O) \to \dots~.$$
(iii) In the original case $n=7$ (Milnor [10]) there is defined an isomorphism
$$\Theta_7 ~\isa~\bbz_{28}~;~\Sigma^7 \mapsto \sigma(W)/8$$
for any framed $8$-dimensional manifold $W$ with $\partial W =\Sigma^7$.
The realization (2.17) of $(\bbz^8,\lambda_{E_8},\mu_{E_8})$ as the kernel form 
of a 4-connected 8-dimensional normal bordism 
$$\eqalign{(f,b)~:~&(M^8,S^7,\partial_+M)~
=~({\rm cl.}(M(E_8,1,\dots,1)\backslash D^8);S^7,\partial M(E_8,1,\dots,1))\cr
& \hskip150pt
\to S^7 \times ([ 0,1];\{0\},\{1\})}$$
gives the exotic sphere
$$\Sigma^7~=~\partial_+M~=~\partial M(E_8,1,\dots,1)$$ 
generating $\Theta_7$\t : the framed 8-dimensional manifold 
$W=M(E_8,1,\dots,1)$ obtained by the $E_8$-plumbing of 8 copies of $\tau_{S^4}$ 
(2.18) has $\sigma(W)=8$. The 7-dimensional homotopy sphere 
$\Sigma^7_{\ell}$ defined for any odd integer $\ell$ in 2.20 is the
boundary of a framed 8-dimensional manifold $W_{\ell}$ with 
$\sigma(W_{\ell})=8(\ell^2-1)$.\break
\line{\hfill$\square$}

For any nonsingular $(-1)$-quadratic form $(K,\lambda,\mu)$ over $\bbz$ 
there exists a symplectic basis ${x_1,\dots ,x_{2m}}$ for $K$, such that
$$\lambda(x_i,x_j)~=~\cases{1&if $i-j=m$\cr
-1&if $j-i=m$\cr
0&otherwise~. }$$
The {\it Arf invariant} of $(K,\lambda,\mu)$ is defined using any such basis 
to be
$$c(K,\lambda,\mu)~=~\sum^m _{i=1}\mu (x_i)\mu (x_{i+m}) \in \bbz_2~.$$
\noindent{\sc Example 3.7} (i) The Arf invariant defines an isomorphism
$$c~:~L_{4k+2}(\bbz) \to \bbz_2~;~(K,\lambda,\mu) \mapsto c(K,\lambda,\mu)$$
The nonsingular $(-1)$-quadratic form $(\bbz \oplus\bbz,\lambda,\mu)$
over $\bbz$ defined by
$$\eqalign{&\lambda((x,y),(x',y'))~=~x'y-xy' \in \bbz~,\cr
&\mu(x,y)~=~x^2+xy+y^2 \in Q_{-1}(\bbz)~=~\bbz_2}$$
has Arf invariant $c(\bbz\oplus \bbz,\lambda,\mu)=1$, 
and so generates $L_{4k+2}(\bbz)$.\hfil\break
(ii) The realization (2.19 (iv)) of the Arf form
$(\bbz\oplus \bbz,\lambda,\mu)$ as the kernel form of a 5-connected
10-dimensional normal bordism 
$$\eqalign{(f,b)~:~&(M^{10},S^9,\partial_+M)~
=~({\rm cl.}(M(I,1,1)\backslash D^{10});S^9,\partial M(I,1,1))\cr
& \hskip100pt \to S^9 \times ([ 0,1];\{0\},\{1\})}$$
is obtained by plumbing together 2 copies of $\tau_{S^5}$ (2.18)
where $I$ is the tree with 1 edge and 2 vertices, and 
$\Sigma^9=\partial_+M=\partial M(I,1,1)$ is the exotic 9-sphere
generating $\Theta^9=\bbz_2$.
Coning off the boundary components gives the closed 10-dimensional 
$PL$ manifold $cS^9 \cup M \cup c \Sigma^9$ without differentiable 
structure of Kervaire [5].\hfill$\square$

\noindent {\sectionfont \S4. Split forms}

A ``split form'' on a $\Lambda$-module $K$ is an element 
$$\psi\in S(K)~=~{\rm Hom}_\Lambda(K,K^*)~,$$ 
which can be regarded as a sesquilinear pairing
$$\psi~:~K \times K \to \Lambda~;~(x,y) \mapsto \psi(x,y)~.$$
\indent
Split forms are more convenient to deal with than $\epsilon$-quadratic
forms in describing the algebraic effects of even-dimensional surgery
(in \S5 below), and are closer to the geometric applications such as
knot theory.  

The main result of \S4 is that the $\epsilon$-quadratic structures
$(\lambda,\mu)$ on a f.g. projective $\Lambda$-module $K$ correspond
to the elements of the $\epsilon$-quadratic group of 2.6
$$Q_{\epsilon}(K)~=~{\rm coker}(1-T_{\epsilon}:S(K) \to S(K))~.$$
The pair of functions 
$(\lambda ,\mu)$ used to define an $\epsilon$-quadratic form 
$(K,\lambda,\mu)$ can thus be replaced by an equivalence class of $\Lambda$-module 
morphisms $\psi:K\to K^*$ such that
$$\eqalign{
&\lambda(x,y) ~=~\psi(x,y) + \epsilon \overline{\psi(y,x)} \in \Lambda~,\ccr
&\mu (x)~=~ \psi(x,x) \in Q_{\epsilon}(\Lambda)}$$
i.e. by an equivalence class of split forms.

\noindent {\sc Definition 4.1}
(i) A {\it split form} $(K,\psi)$ over $\Lambda$ is a
f.\ g. projective $\Lambda$-module $K$ together with an element
$\psi\in S(K)$.\hfil\break 
(ii) A {\it morphism} (resp. {\it isomorphism}) of split forms over $\Lambda$
$$f~:~(K,\psi) \to (K',\psi')$$
is a $\Lambda$-module morphism (resp. isomorphism) $f:K \to K'$ such that
$$f^* \psi' f~=~\psi~:~K \to K^*~.$$
(iii) An {\it $\epsilon$-quadratic morphism} (resp. {\it isomorphism})
of split forms over $\Lambda$
$$(f,\chi)~:~(K,\psi) \to (K',\psi')$$
is a $\Lambda$-module morphism (resp. isomorphism) $f:K \to K'$ together 
with an element $\chi \in Q_{-\epsilon}(K)$ such that
$$f^* \psi' f - \psi~=~\chi-\epsilon\chi^*~:~K \to K^*~.$$
(iv) A split form $(K,\psi)$ is {\it $\epsilon$-nonsingular} 
if $\psi+\epsilon \psi^*:K \to K^*$ is a $\Lambda$-module isomorphism.\hfill$\square$

\noindent {\sc Proposition 4.2} {\it 
{\rm (i)} A split form $(K,\psi)$ determines an $\epsilon$-quadratic form 
$(K,\lambda,\mu)$ by
$$\eqalign{
&\lambda~=~\psi+\epsilon\psi^*~: ~K \to K^*~;~
x\mapto (y\mapto \psi(x,y) +\epsilon \overline {\psi(y,x)})~, \ccr
&\mu~:~K \to Q_{\epsilon}(\Lambda)~; ~x\mapto \psi(x,x)~.}$$
{\rm (ii)} Every $\epsilon$-quadratic form $(K,\lambda,\mu)$ is
determined by a split form $(K,\psi)$, which is unique up to 
$$\psi ~\sim~\psi'~\hbox{\it if}~~\psi'-\psi~=~\chi-\epsilon \chi^*~
\hbox{\it for some}~\chi:K \to K^*~.$$
{\rm (iii)} 
The isomorphism classes of (nonsingular) $\epsilon$-quadratic forms 
$(K,\lambda,\mu)$ over $\Lambda$ are in one-one correspondence
with the $\epsilon$-quadratic isomorphism classes 
of ($\epsilon$-nonsingular) split forms $(K,\psi)$ over $\Lambda$.} \hfil\break
{\sc Proof}\t :
(i) By construction.\hfil\break
(ii) There is no loss of generality in taking $K$ to be f.g. free, $K=\Lambda^k$.
An $\epsilon$-quadratic form $(\Lambda^k,\lambda,\mu)$ over $\Lambda$
is determined by a $k \times k$-matrix 
$\lambda=\{\lambda_{ij}\in \Lambda\,\vert\, 1 \le i,j \le k\}$ such that
$$\overline {\lambda_{ij}} ~=~ \epsilon \lambda_{ji} \in \Lambda~$$
and a collection of elements 
$\mu=\{\mu_i \in Q_{\epsilon}(\Lambda)\,\vert\,1 \le i \le k\}$ such that
$$\mu_i + \epsilon \overline{\mu_i}~=~\lambda_{ii} \in Q^{\epsilon}(\Lambda)~.$$
Choosing any representatives $\mu_i \in \Lambda$ of
$\mu_i \in Q_{\epsilon}(\Lambda^k)$ there is defined a split 
form $(\Lambda^k,\psi)$ with $\psi=\{\psi_{ij} \in 
\Lambda \,\vert\, 1 \le i,j \le k \}$
the $k \times k$ matrix defined by
$$\psi_{ij}~=~\cases{\lambda_{ij} &if $i<j$ \ccr \mu_i &if $i=j$
\ccr 0 &otherwise~.}$$
(iii) An $\epsilon$-quadratic (iso)morphism $(f,\chi):(K,\psi) \to (K',\psi')$
of split forms determines an (iso)morphism $f:(K,\lambda,\mu) \to (K',\lambda',\mu')$
of $\epsilon$-quadratic forms. 
Conversely, an $\epsilon$-quadratic form 
$(K,\lambda,\mu)$ determines an $\epsilon$-quadratic isomorphism class
of split forms $(K,\psi)$ as in 3.1, and every (iso)morphism of $\epsilon$-quadratic forms lifts to
an $\epsilon$-quadratic (iso)morphism of split forms.\hfill$\square$

Thus $Q_{\epsilon}(K)$ is both the group of isomorphism classes of
$\epsilon$-quadratic forms and the group of $\epsilon$-quadratic isomorphism
classes of split forms on a f.\ g. projective $\Lambda$-module $K$.

The following algebraic result will be used in 4.6 below to obtain a
homological split form $\psi$ on the kernel $\bbz[\pi_1(X)]$-module
$K_n(M)$ of an $n$-connected $2n$-dimensional normal map $(f,b):M \to
X$ with some extra structure, which determines the kernel
$(-1)^n$-quadratic form $(K_n(M),\lambda,\mu)$ as in 4.2 (i).

\noindent{\sc Lemma 4.3} {\it Let $(K,\lambda,\mu)$ be an
$\epsilon$-quadratic form over $\Lambda$.\hfil\break 
{\rm (i)} If $s:K \to K$ is an endomorphism such that
$$\pmatrix{s \cr 1-s}~:~(K,0,0) \to (K,\lambda,\mu) \oplus (K,-\lambda,-\mu)$$
defines a morphism of $\epsilon$-quadratic forms then $(K,\lambda s)$ is a split form 
which determines the $\epsilon$-quadratic form $(K,\lambda,\mu)$.\hfil\break
{\rm (ii)} If $(K,\lambda,\mu)$ is nonsingular and 
$(K,\psi)$ is a split form which determines $(K,\lambda,\mu)$ then 
$$s~=~\lambda^{-1}\psi~:~K \to K$$ 
is an endomorphism such that
$$\pmatrix{s \cr 1-s}~:~(K,0,0) \to (K,\lambda,\mu) \oplus (K,-\lambda,-\mu)$$
defines a morphism of $\epsilon$-quadratic forms.}\hfil\break
\noindent{\sc Proof}\t : (i) By 4.2 there exist a split form $(K,\psi)$ which
determines $(K,\lambda,\mu)$ and an $\epsilon$-quadratic morphism of split forms 
$$\big(\pmatrix{s \cr 1-s},\,\chi\,\big)~:~(K,0) \to (K,\psi) \oplus (K,-\psi)~.$$
It follows from
$$\eqalign{&\lambda~=~\psi+\epsilon \psi^* ~:~ K \to K^*~,\cr
&\pmatrix{s \cr 1-s}^*\pmatrix{ \psi & 0 \cr 0 &-\psi}\pmatrix{s \cr 1-s}~=~
\chi-\epsilon\chi^*~:~K \to K^*}$$
that 
$$\lambda s - \psi~=~\chi'-\epsilon \chi'^*~:~K \to K^*$$
with
$$\chi'~=~\chi-s^*\psi~:~K \to K^*~.$$
(ii) From the definitions.\hfill$\square$

In the terminology of \S5 the morphism of 4.3 (ii)
$$\pmatrix{s \cr 1-s}~:~(K,0,0) \to (K,\lambda,\mu) \oplus (K,-\lambda,-\mu)$$
is the inclusion of a lagrangian
$$\eqalign{L~&=~{\rm im}(\pmatrix{s \cr 1-s}:K \to K \oplus K)\cr
&=~{\rm ker}(\pmatrix{(-1)^{n-1}\psi^*&\psi}:K\oplus K \to K^*)~.}$$

\noindent {\sc Example 4.4} 
A {\it $(2n-1)$-knot} is an embedding of a homotopy $(2n-1)$-sphere in
a standard $(2n+1)$-sphere
$$\ell~:~\Sigma^{2n-1} \hookrightarrow S^{2n+1}~.$$
For $n=1$ this is just a classical knot
$\ell:\Sigma^1=S^1 \hookrightarrow S^3$; for $n \geq 3$ $\Sigma^{2n-1}$ is
homeomorphic to $S^{2n-1}$, by the generalized Poincar\'e conjecture,
but may have an exotic differentiable structure.  Split forms
$(K,\psi)$ first appeared as the Seifert forms over $\bbz$ of
$(2n-1)$-knots, originally for $n=1$.
See Ranicki [21,\t 7.8],\t [24] for a surgery treatment of
high-dimensional knot theory. In particular, a Seifert form 
is an integral refinement of an even-dimensional surgery kernel form, 
as follows.\hfil\break
(i) A $(2n-1)$-knot $\ell:\Sigma^{2n-1} \hookrightarrow S^{2n+1}$ is
{\it simple} if 
$$\pi_r(S^{2n+1}\backslash \ell(\Sigma^{2n-1}))~=~\pi_r(S^1)~~
(1 \leq r \leq n-1)~.$$
(Every $1$-knot is simple).
A simple $(2n-1)$-knot $\ell$ has a simple Seifert surface, that is an
$(n-1)$-connected framed codimension 1 submanifold $M^{2n} \subset
S^{2n+1}$ with boundary
$$\partial M ~=~ \ell(\Sigma^{2n-1}) \subset S^{2n+1}~.$$ 
The kernel of the $n$-connected normal map
$$(f,b)~=~{\rm inclusion}~:~(M,\partial M) \to (X,\partial X)~=~
(D^{2n+2},\ell(\Sigma^{2n-1}))$$
is a nonsingular $(-1)^n$-quadratic form $(H_n(M),\lambda,\mu)$ over $\bbz$.
The {\it Seifert form} of $\ell$ with respect to $M$ is the refinement
of $(H_n(M),\lambda,\mu)$ to a $(-1)^n$-nonsingular split form 
$(H_n(M),\psi)$ over $\bbz$ which is defined using
Alexander duality and the universal coefficient theorem
$$\psi~=~i_*~:~H_n(M) \to H_n(S^{2n+1}\backslash M)~\cong~H^n(M)~\cong~
H_n(M)^*$$
with $i:M \to S^{2n+1}\backslash M$ the map pushing $M$ off itself
along a normal direction in $S^{2n+1}$. 
If $i':M \to S^{2n+1}\backslash M$
pushes $M$ off itself in the opposite direction
$$i'_*~=~(-1)^{n+1}\psi^*~:~H_n(M) \to H_n(S^{2n+1}\backslash M)~
\cong~H^n(M)~\cong~H_n(M)^*$$
with
$$\eqalign{i_*-i'_*~&=~\psi+(-1)^n\psi^*~=~\lambda\cr
&=~([M]\cap -)^{-1}~:~H_n(M) \to H^n(M)~\cong~H_n(M)^*}$$
the Poincar\'e duality isomorphism. 
If $x_1,x_2,\dots,x_k \in H_n(M)$ is a basis then $(\psi(x_j,x_{j'}))$
is a {\it Seifert matrix} for the $(2n-1)$-knot $\ell$.
For any embeddings $x,y:S^n \hookrightarrow M$ 
$$\eqalign{\psi(x,y)~
&=~\hbox{\rm linking number($i x(S^n)\cup y(S^n) \subset S^{2n+1}$)}\ccr
&=~\hbox{\rm degree}(y^*ix:S^n \to S^n) \in \bbz}$$
with
$$y^*ix~:~S^n \xymatrix@C-5pt{\ar[r]^-{\displaystyle{x}}&}M 
\xymatrix@C-5pt{\ar[r]^-{\displaystyle{i}}&} S^{2n+1}\backslash M 
\xymatrix@C-5pt{\ar[r]^-{\displaystyle{y^*}}&} 
S^{2n+1}\backslash y(S^n) \simeq S^n~.$$
For $n \geq 3$ every element $x \in H_n(M)$ is represented by an
embedding $e:S^n \hookrightarrow M$, using the Whitney embedding theorem,
and $\pi_1(M)=\{1\}$. 
Moreover, for any embedding $x:S^n \hookrightarrow M$ the framed embedding
$M \hookrightarrow S^{2n+1}$ determines a stable trivialization
of the normal bundle $\nu_x:S^n \to BSO(n)$
$$\delta\nu_x~:~\nu_x \oplus \epsilon~\cong~\epsilon^{n+1}$$
such that
$$\psi(x,x)~=~(\delta\nu_x,\nu_x) \in \pi_{n+1}(BSO(n+1),BSO(n))~=~\bbz~.$$ 
Every element $x \in H_n(M)$ is represented by an embedding
$$e_1\times e_2~:~S^n \hookrightarrow M \times \bbr$$ 
with $e_1:S^n \looparrowright M$ a framed immersion such that the composite
$$S^n \xymatrix@C+10pt{\ar[r]^-{\displaystyle{e_1 \times e_2}}&} M \times \bbr 
\hookrightarrow S^{2n+1}$$
is isotopic to the standard framed embedding $S^n \hookrightarrow S^{2n+1}$.
Then
$$\psi(x,x)~=~\sum\limits_{(a,b) \in D_2(e_1),e_2(a)<e_2(b)}I(a,b) \in \bbz$$
is an integral lift of the geometric self-intersection (2.15 (ii))
$$\mu(x)~=~\sum\limits_{(a,b) \in D_2(e_1)/\bbz_2}I(a,b) \in Q_{(-1)^n}(\bbz)~,$$
with 
$$D_2(e_1)~=~\{(a,b) \in S^n \times S^n\,\vert\,a\neq b\in S^n,
e_1(a)=e_1(b)\in M\}$$ 
the double point set. 
For even $n$ $\psi(x,x)=\mu(x) \in Q_{+1}(\bbz)=\bbz$, while for
odd $n$ $\psi(x,x) \in \bbz$ is a lift of $\mu(x) \in Q_{-1}(\bbz)=\bbz_2$.
The Seifert form $(H_n(M),\psi)$ is such that $\psi(x,x)=0$ if (and
for $n \geq 3$ only if) $x \in H_n(M)$ can be killed by an ambient
surgery on $M^{2n} \subset S^{2n+1}$, i.e. represented by a framed
embedding of pairs 
$$x~:~(D^{n+1}\times D^n,S^n \times D^n) \hookrightarrow
(S^{2n+1}\times [ 0,1],M \times \{0\})$$
so that the effect of the surgery on $M$ is another Seifert surface 
for the $(2n-1)$-knot $\ell$
$$M'~=~{\rm cl.}(M \times x(S^n \times D^n)) \cup D^{n+1} \times S^{n-1}
\subset S^{2n+1}~.$$
If $x \in H_n(M)$ generates a direct summand 
$L=\langle x \rangle \subset H_n(M)$ then $M'$ is also $(n-1)$-connected, 
with Seifert form
$$(H_n(M'),\psi')~=~(L^{\perp}/L,[\psi])~,$$
where
$$L^{\perp}~=~\{y \in H_n(M)\,\vert\, (\psi+(-1)^n\psi^*)(x)(y)=0~{\rm for}~
x \in L\} \subseteq H_n(M)~.$$
(ii) Every $(-1)^n$-nonsingular split form $(K,\psi)$ over $\bbz$ 
is realized as the Seifert form of a simple $(2n-1)$-knot
$\ell:\Sigma^{2n-1} \hookrightarrow S^{2n+1}$ (Kervaire [6]). 
From the algebraic surgery point of view the realization proceeds as follows.
By 2.17 the nonsingular $(-1)^n$-quadratic form $(K,\lambda,\mu)$
determined by $(K,\psi)$ (4.2 (i)) is the kernel form of an $n$-connected
$2n$-dimensional normal map
$$(f,b)~:~(M^{2n},\Sigma^{2n-1}) \to (D^{2n},S^{2n-1})$$
with $f\vert:\Sigma^{2n-1} \to S^{2n-1}$ a homotopy equivalence.
The double of $(f,b)$ defines an $n$-connected $2n$-dimensional normal map
$$(g,c)~=~(f,b)\cup -(f,b)~:~N^{2n}~=~M\cup_{\Sigma^{2n-1}}-M \to
D^{2n} \cup_{S^{2n-1}}-D^{2n}~=~S^{2n}$$
with kernel form $(K \oplus K,\lambda\oplus -\lambda,\mu\oplus -\mu)$.
The direct summand
$$L~=~{\rm ker}(\pmatrix{(-1)^{n-1}\psi^*&\psi}:K\oplus K \to K^*) 
\subset K \oplus K$$
is such that for any $(x,y)\in L$
$$\mu(x)-\mu(y)~=~(1+(-1)^{n-1})\psi(x,x)~=~0 \in Q_{(-1)^n}(\bbz)~.$$
Let $k={\rm rank}_{\bbz}(K)$.
The trace of the $k$ surgeries on $(g,c)$ killing a basis 
$(x_j,y_j) \in K\oplus K$ ($j=1,2,\dots,k$) for $L$
is an $n$-connected $(2n+1)$-dimensional normal map
$$(W^{2n+1};N,S^{2n}) \to S^{2n} \times ([ 0,1];\{0\},\{1\})$$
such that
$$\ell~:~\Sigma^{2n-1} \hookrightarrow (\Sigma^{2n-1} \times D^2) \cup 
(W\cup D^{2n+1}) \cup (M \times [ 0,1])~\cong~S^{2n+1}$$
is a simple $(2n-1)$-knot with Seifert surface $M$ and Seifert form $(K,\psi)$.
Note that $M$ itself is entirely determined by the $(-1)^n$-quadratic form
$(K,\lambda,\mu)$, with ${\rm cl.}(M \backslash D^{2n})$ the trace of $k$ 
surgeries on $S^{2n-1}$ removing  
$$\bigcup_k S^{n-1} \times D^n \hookrightarrow S^{2n-1}$$
with (self-)linking numbers $(\lambda,\mu)$.  The embedding $M
\hookrightarrow S^{2n+1}$ is determined by the choice of split
structure $\psi$ for $(\lambda,\mu)$.  \hfil\break
(iii) In particular, (ii) gives a knot version
of the plumbing construction (2.18): let $G$ be a finite graph with
vertices $v_1,v_2,\dots,v_k$, weighted by $\mu_1,\mu_2,\dots,\mu_k \in
Q_{(-1)^n}(\bbz)$, so that there are defined a $(-1)^n$-quadratic form
$(\bbz^k,\lambda,\mu)$ and a plumbed stably parallelizable
$(n-1)$-connected $2n$-dimensional manifold with boundary
$$M^{2n}~=~M(G,\mu_1,\mu_2,\dots,\mu_k)~,$$
killing $H_1(G)$ by surgery if $G$ is not a tree.
A choice of split form $\psi$ for $(\lambda,\mu)$ determines a
compression of a framed embedding $M \hookrightarrow S^{2n+j}$ ($j$ large) 
to a framed embedding $M \hookrightarrow S^{2n+1}$, so that 
$\partial M \hookrightarrow S^{2n+1}$ is a codimension 2 framed embedding.  
The form $(\bbz^k,\lambda,\mu)$ is nonsingular if and only if
$\Sigma^{2n-1}=\partial M$ is a homotopy $(2n-1)$-sphere, in which case
$\Sigma^{2n-1} \hookrightarrow S^{2n+1}$ 
is a simple $(2n-1)$-knot with simple Seifert surface $M$.\hfil\break
(iv) Given a simple $(2n-1)$-knot $\ell:\Sigma^{2n-1} \hookrightarrow S^{2n+1}$
and a simple Seifert surface $M^{2n}\hookrightarrow S^{2n+1}$ there is
defined an $n$-connected $2n$-dimensional normal map
$$(f,b)~=~{\rm inclusion}~:~(M,\partial M) \to (X,\partial X)~=~
(D^{2n+2},\ell(\Sigma^{2n-1}))$$
as in (i). The knot complement is a $(2n+1)$-dimensional manifold with boundary 
$$(W,\partial W)~=~({\rm cl.}(S^{2n+1} \backslash (\ell(\Sigma^{2n-1}) \times D^2)),
\ell(\Sigma^{2n-1})\times S^1)$$
with a $\bbz$-homology equivalence $p:(W,\partial W) \to S^1$ such that
$$\eqalign{
&p\vert~=~{\rm projection}~:~
\partial W~=~\Sigma^{2n-1} \times S^1 \to S^1~,\cr
&p^{-1}({\rm pt.})~=~M \subset W~.}$$
Cutting $W$ along $M \subset W$ there is obtained a cobordism
$(N;M,M')$ with $M'$ a copy of $M$, and $N$ a deformation retract of
$S^{2n+1}\backslash M$, such that $(f,b)$ extends to an $n$-connected
normal map
$$(g,c)~:~(N;M,M') \to X \times ([ 0,1];\{0\},\{1\})$$
with $(g,c)\vert=(f',b'):M' \to X$ a copy of $(f,b)$. The 
$n$-connected $(2n+1)$-dimensional normal map
$$\eqalign{&(h,d)~=~(g,c)/((f,b)=(f',b'))~:\cr
&(W,\partial W)~=~(N;M,M')/(M=M') \to (X,\partial X) \times S^1}$$
is a $\bbz$-homology equivalence which is the identity on $\partial W$, 
and such that 
$$(f,b)~=~(h,d)\vert~:~
(M,\partial M)~=~h^{-1}((X,\partial X)\times \{{\rm pt.}\})
\to (X,\partial X)~.\eqno{\hbox{$\square$}}$$

\noindent {\sc Example 4.5} 
(i) Split forms over group rings arise in the following
geometric situation, generalizing 4.4 (iv).\hfil\break
Let $X$ be a $2n$-dimensional Poincar\'e complex, and let $(h,d):W \to
X \times S^1$ be an $n$-connected $(2n+1)$-dimensional normal map which is a
$\bbz[\pi_1(X)]$-homology equivalence. Cut $(h,d)$ along $X \times
\{{\rm pt.}\} \subset X \times S^1$ to obtain an $n$-connected
$2n$-dimensional normal map 
$$(f,b)~=~(h,d)\vert~:~M~=~h^{-1}(\{{\rm pt.}\}) \to X$$ 
and an $n$-connected normal bordism 
$$(g,c)~:~(N;M,M') \to X \times ([ 0,1];\{0\},\{1\})$$
with $N$ a deformation retract of $W\backslash M$,
such that $(g,c)\vert = (f,b):M \to X$, and such that $(g,c)\vert =
(f',b'):M' \to X$ is a copy of $(f,b)$. 
The inclusions $i:M \hookrightarrow N$, 
$i':M' \hookrightarrow N$ induce $\bbz[\pi_1(X)]$-module morphisms 
$$i_*~:~K_n(M) \to K_n(N)~~,~~i'_*~:~K_n(M')~=~K_n(M) \to K_n(N)$$ 
which fit into an exact sequence
$$\xymatrix@C+10pt{K_{n+1}(W)=0 \ar[r]& K_n(M) 
\ar[r]^-{\displaystyle{i_*-i'_*}}& K_n(N) \ar[r]&K_n(W)=0~,}$$
so that $i_*-i'_*:K_n(M) \to K_n(N)$ is an isomorphism. Let 
$(K_n(M),\lambda,\mu)$ be the kernel $(-1)^n$-quadratic form of $(f,b)$. 
The endomorphism
$$s~=~(i_*-i'_*)^{-1}i_*~:~K_n(M) \to K_n(M)$$ 
is such that 
$$\pmatrix{s \cr 1-s}~:~
(K_n(M),0,0) \to (K_n(M),\lambda,\mu) \oplus (K_n(M),-\lambda,-\mu)$$
defines a morphism of $(-1)^n$-quadratic forms, 
so that by 4.3 the split form $(K_n(M),\psi)$ with
$$\psi~=~\lambda s~:~K_n(M) \to K_n(M)^*$$
determines $(K_n(M),\lambda,\mu)$. Every element $x \in K_n(M)$ can 
be represented by a framed immersion 
$x:S^n \looparrowright M$ with a null-homotopy $fx \simeq *:S^n
\to X$.  Use the null-homotopy and the normal $\bbz[\pi_1(X)]$-homology 
equivalence $(h,d):W \to X \times S^1$ to extend $x$ to a framed immersion
$\delta x:D^{n+1} \looparrowright W$.  If $x_1,x_2,\dots,x_k \in
K_n(M)$ is a basis for the kernel f.g.  free $\bbz[\pi_1(X)]$-module
then
$$s(x_j)~=~\sum\limits^k_{j'=1}s_{jj'}x_{j'} \in K_n(M)$$
with 
$$\eqalign{
s_{jj'}~&=~\hbox{\rm linking number}
(ix_j(S^n) \cup x_{j'}(S^n) \subset W)\cr
&=~\hbox{\rm intersection number}
(ix_j(S^n) \cap \delta x_{j'}(D^{n+1}) \subset W)\in \bbz[\pi_1(X)]~.}$$
The split form $(K_n(M),\psi)$ is thus a (non-simply connected) Seifert
form.\hfil\break
(ii) Suppose given an $n$-connected $2n$-dimensional normal map
$(f,b):(M,\partial M) \to (X,\partial X)$, with kernel
$(-1)^n$-quadratic form $(K_n(M),\lambda,\mu)$ over $\bbz [\pi_1(X)]$. 
A choice of split form $\psi$ for $(\lambda,\mu)$ can be realized
by an $(n+1)$-connected $(2n+2)$-dimensional normal map
$$(g,c)~:~(L,\partial L) \to (X\times D^2,X \times S^1 \cup \partial X \times D^2)$$
which is a $\bbz[\pi_1(X)]$-homology equivalence with
$$\eqalign{&(f,b)~=~(g,c)\vert ~:~(M,\partial M)~=~
g^{-1}((X,\partial X) \times \{0\}) \to (X,\partial X)~,\ccr
&H_{*+1}(\widetilde{L},\widetilde{M})~=~K_*(M)~(=~0~{\rm for}~* \neq n)}$$
as follows. The inclusion $\partial M \hookrightarrow \partial L$ is a 
codimension 2 embedding with Seifert surface $M \hookrightarrow \partial L$ and 
Seifert form $(K_n(M),\psi)$ as in the relative version of (i), with
$$(h,d)~=~(g,c)\vert~:~W~=~\partial L \to 
X \times S^1 \cup \partial X \times D^2~.$$
The choice of split form $\psi$ for $(\lambda,\mu)$ determines a
sequence of surgeries on the $n$-connected $(2n+1)$-dimensional normal
map
$$(f,b) \times 1_{[ 0,1]}~ :~ M \times ([ 0,1];\{0\},\{1\}) \to
X \times ([ 0,1];\{0\},\{1\})$$ 
killing the (stably) f.\ g. free $\bbz[\pi_1(X)]$-module 
$$K_n(M \times [ 0,1])~=~K_n(M)~,$$
obtaining a rel $\partial$ normal bordant map
$$(f_N,b_N)~:~(N;M,M') \to X \times ([ 0,1];\{0\},\{1\})$$ 
with $K_i(N)=0$ for $i \neq n$.  The $\bbz[\pi_1(X)]$-module morphisms
induced by the inclusions $i:M \hookrightarrow N$, $i':M'
\hookrightarrow N$
$$\eqalign{&i_* ~:~ K_n(M) \to K_n(N)~~,~~i'_*~:~ K_n(M)~=~K_n(M') \to K_n(N)}$$
are such that $i_*-i'_*:K_n(M) \to K_n(N)$ is a 
$\bbz[\pi_1(X)]$-module isomorphism, with
$$\psi~ :~ K_n(M)~ \raise4pt\hbox{$(i_*-i'_*)^{-1}i_* \atop \ra{4}$}~ K_n(M)~
\raise4pt\hbox{${\rm adjoint}(\lambda) \atop \ra{4}$}~ K_n(M)^*~.$$
Thus it is possible to identify
$$\eqalign{&i_*~=~\psi~:~ K_n(M) \to K_n(N)~\cong~K_n(M)^*~,\ccr
&i'_*~=~(-1)^{n+1}\psi^*~:~ K_n(M)~=~K_n(M') \to K_n(N)~\cong~K_n(M)^*}$$
with
$$i_*-i'_*~=~\psi+(-1)^n\psi^*~=~{\rm adjoint}(\lambda)~:~K_n(M)~\isa~K_n(M)^*~.$$
The $(2n+1)$-dimensional manifold with boundary defined by
$$(V,\partial V)~=~(N/(M=M'),\partial M \times S^1)$$
is equipped with a normal map
$$(f_V,b_V)~:~(V,\partial V) \to (X \times S^1,\partial X \times S^1)$$
which is an $n$-connected $\bbz[\pi_1(X)]$-homology equivalence, with
$K_j(V)=0$ for $j \neq n+1$ and
$$K_{n+1}(V)~=~\hbox{\rm coker}(z\psi+(-1)^n\psi^*:K_n(M)[z,z^{-1}]
\to K_n(M)^*[z,z^{-1}])$$
identifying 
$$\bbz[\pi_1(X\times S^1)]~=~\bbz[\pi_1(X)][z,z^{-1}]~~(\overline{z}=z^{-1})~.$$ 
The trace of the surgeries on 
$(f,b)\times 1_{[ 0,1]}$ gives an extension of $(f_V,b_V)$ 
to an $(n+1)$-connected $(2n+2)$-dimensional normal bordism
$$(f_U,b_U)~:~(U;V,M \times S^1) \to X \times S^1 \times ([ 0,1];\{0\},\{1\})$$
with $K_i(U)=0$ for $i \neq n+1$ and (singular) kernel
$(-1)^{n+1}$-quadratic form over $\bbz[\pi_1(X)][z,z^{-1}]$
$$\eqalign{&(K_{n+1}(U),\lambda_U,\mu_U)\ccr
&=~(K_n(M)[z,z^{-1}],(1-z)\psi+(-1)^{n+1}(1-z^{-1})\psi^*,(1-z)\psi)~.}$$
The $(2n+2)$-dimensional manifold with boundary defined by
$$(W,\partial W)~=~(M \times D^2 \cup U,\partial M \times D^2 \cup V)$$
is such that $(f,b)$ extends to an $(n+1)$-connected normal map
$$(g,c)~=~(f,b) \times 1_{D^2} \cup (f_U,b_U)~:~(W,\partial W) \to 
(X \times D^2,\partial X \times D^2 \cup X \times S^1)$$
which is a $\bbz[\pi_1(X)]$-homology equivalence, with 
$H_{n+1}(\widetilde{W},\widetilde{M})=K_n(M)$. 
See Example 27.9 of Ranicki [24] for further details (noting
that the split form $\psi$ here corresponds to the asymmetric form 
$\lambda$ there). \hfil\break
(iii) Given a simple knot $\ell:\Sigma^{2n-1} \hookrightarrow S^{2n+1}$ 
and a simple Seifert surface $M^{2n} \subset S^{2n+1}$ there is defined an 
$n$-connected normal map
$$(f,b)~=~{\rm inclusion}~:~
(M^{2n},\partial M) \to (X,\partial X)~=~(D^{2n+2},\ell(\Sigma^{2n-1}))$$
with a Seifert form $\psi$ on $K_n(M)=H_n(M)$, as in 4.4.
For $n \geq 2$ the surgery construction of (i) applied to $(f,b)$, $\psi$
recovers the knot 
$$\ell~:~\Sigma^{2n-1}~=~\partial M \hookrightarrow \partial W~=~S^{2n+1}$$
with $M^{2n} \subset W=D^{2n+2}$ the Seifert surface pushed into the 
interior of $D^{2n+2}$. The knot complement 
$$(V^{2n+1},\partial V) ~=~
({\rm cl.}\big(S^{2n+1}\backslash(\ell(\Sigma^{2n-1}) \times D^2)\big),
\ell(\Sigma^{2n-1})\times S^1)$$ 
is such that there is defined an $n$-connected $(2n+1)$-dimensional normal
map
$$(f_V,b_V) ~:~ (V,\partial V) \to (X,\partial X)\times S^1$$
which is a homology equivalence, with
$$(f_V,b_V)\vert~=~(f,b)~:~(M,\partial M)~=~(f_V)^{-1}
((X,\partial X)\times \{*\}) \to (X,\partial X)~.$$
Cutting $(f_V,b_V)$ along $(f,b)$ results in a normal map as in (i)
$$(f_N,b_N)~:~(N^{2n+1};M^{2n},M'^{2n}) \to X \times 
([ 0,1];\{0\},\{1\}) ~.\eqno{\square}$$
\eject

\noindent {\sectionfont \S5. Surgery on forms}

\S5 develops algebraic surgery on forms.  The effect of a geometric
surgery on an $n$-connected $2n$-dimensional normal map is an algebraic
surgery on the kernel $(-1)^n$-quadratic form.  Moreover, geometric
surgery is possible if and only if algebraic surgery is possible.

Given an $\epsilon$-quadratic form $(K,\lambda,\mu)$ over $\Lambda$
it is possible to kill an element $x \in K$ by algebraic surgery 
if and only if $\mu(x)=0 \in Q_{\epsilon}(\Lambda)$ and $x$
generates a direct summand $\langle x\rangle =\Lambda x \subset K$.
The effect of the surgery is the $\epsilon$-quadratic form 
$(K',\lambda',\mu')$ defined on the subquotient 
$K'=\langle x\rangle^{\bot}/\langle x\rangle $ of $K$, 
with $\langle x\rangle^{\bot}=\{y \in K \,\vert\, \lambda(x,y)=0\in \Lambda\}$. 

\noindent {\sc Definition 5.1}
(i) Given an $\epsilon$-symmetric form $(K,\lambda)$ 
and a submodule $L\subseteq K$ define the {\it orthogonal} submodule
$$\eqalign{L^{\bot}~
&=~\{x\in K\,\vert\,\lambda(x,y)=0\in \Lambda~\hbox{\rm for all}~y\in L\}\ccr
&=~{\rm ker}(i^* \lambda : K \to L^*)}$$
with $i:L \to K$ the inclusion. If $(K,\lambda)$ is nonsingular and
$L$ is a direct summand of $K$ then so is $L^{\bot}$. \hfil\break
(ii) A {\it sublagrangian}
of a nonsingular $\epsilon$-quadratic form $(K,\lambda,\mu)$ over $\Lambda$
is a direct summand $L \subseteq K$ such that
$$\mu (L)~=~\{0\} \subseteq Q_{\epsilon}(\Lambda)~,$$
and
$$\lambda(L)(L)~=~\{0\}~~,~~L\subseteq L^{\bot}~.$$
(iii) A {\it lagrangian}
of $(K,\lambda,\mu)$ is a sublagrangian $L$ such that $L^{\bot}=L$.
\hfill$\square$

The main result of \S5 is that the inclusion of a sublagrangian is a
morphism of $\epsilon$-quadratic forms
$$i~:~(L,0,0) \to (K,\lambda,\mu)$$
which extends to an isomorphism
$$f~:~H_{\epsilon}(L) \oplus (L^{\bot}/L,[\lambda],[\mu])~ 
\isa ~(K,\lambda,\mu)$$
with $H_{\epsilon}(L)$ the hyperbolic $\epsilon$-quadratic form (2.14).

\noindent {\sc Example 5.2}
Let $(f,b):M^{2n} \to X$ be an $n$-connected $2n$-dimensional normal
map with kernel $(-1)^n$-quadratic form $(K_n(M),\lambda,\mu)$ over
$\bbz[\pi_1(X)]$, and $n \ge 3$.  An element $x \in K_n(M)$ generates a
sublagrangian $L=\langle x\rangle \subset K_n(M)$ if and only if it can
be killed by surgery on $S^n \times D^n \hookrightarrow M$ with trace
an $n$-connected normal bordism
$$((g,c);(f,b),(f',b'))~:~(W^{2n+1};M^{2n},M'^{2n}) \to 
X \times ([ 0,1];\{0\},\{1\})$$
such that $K_{n+1}(W,M')=0$. The kernel form of the effect of such a surgery
$$(f',b')~:~M'~=~{\rm cl.}(M \backslash S^n \times D^n) 
\cup D^{n+1} \times S^{n-1} \to X$$
is given by
$$(K_n(M'),\lambda',\mu') ~=~ (L^{\bot}/L,[\lambda],[\mu])~.$$
There exists an $n$-connected normal bordism $(g,c)$ of $(f,b)$ to a
homotopy equivalence $(f',b')$ with $K_{n+1}(W,M')=0$ if and only if
$(K_n(M),\lambda,\mu)$ admits a lagrangian. \hfill$\square$

\noindent {\sc Remark 5.3}
There are other terminologies. In the classical theory of quadratic 
forms over fields a lagrangian is a ``maximal isotropic subspace''.
Wall called hyperbolic forms ``kernels'' and the lagrangians ``subkernels''.
Novikov called hyperbolic forms ``hamiltonian'', and introduced the name 
``lagrangian'', because of the analogy with the hamiltonian formulation of 
physics.
\hfill$\square$

\noindent {\sc Example 5.4}
An $n$-connected $(2n+1)$-dimensional normal bordism
$$((g,c);(f,b),(f',b'))~:~(W^{2n+1};M^{2n},M'^{2n})\to 
X \times ([ 0,1];\{0\},\{1\})$$
with $K_{n+1}(W,M')=0$ determines a sublagrangian
$$L~=~{\rm im}(K_{n+1}(W,M) \to K_n(M)) \subset K_n(M)$$
of the kernel $(-1)^n$-quadratic form $(K,\lambda,\mu)$ of $(f,b)$,
with $K=K_n(M)$. The sublagrangian
$L$ is a lagrangian if and only if $(f',b')$ is a homotopy equivalence.
$W$ has a handle decomposition on $M$ of the type
$$W~=~M \times I \cup \bigcup_k (n+1)\hbox{\rm -handles}~D^{n+1} \times D^n~,$$
and $L \cong K_{n+1}(W,M)\cong\bbz[\pi_1(X)]^k$ is a f.\ g. free 
$\bbz[\pi_1(X)]$-module with rank the number $k$ of $(n+1)$-handles.
The exact sequences of stably f.\ g. free $\bbz[\pi_1(X)]$-modules
$$\eqalign{&0 \to K_{n+1}(W,M) \to K_n(M) \to K_n(W) \to 0~,\ccr
&0 \to K_n(M') \to K_n(W) \to K_n(W,M') \to 0}$$
are isomorphic to 
$$\eqalign{&0\to ~L~ \raise4pt\hbox{$i \atop \ra{1.5} $} ~
K \to K/L \to 0~,\ccr 
&0 \to L^{\bot}/L \to K/L~ \raise4pt\hbox
{$[i^* \lambda] \atop \ra{2}$} ~L^* \to 0~.}\eqno{\lower12pt\hbox{$\square$}}$$

\noindent {\sc Definition 5.5}
(i) A {\it sublagrangian} of an $\epsilon$-nonsingular split 
form $(K,\psi)$ is an $\epsilon$-quadratic morphism of split forms
$$(i,\theta) ~:~ (L,0) \to (K,\psi)$$
with $i : L \to K$ a split injection. \hfil\break
(ii) A {\it lagrangian} of $(K,\psi)$ is a sublagrangian such that the
sequence
$$0 \to L ~ \raise4pt\hbox{$i \atop \ra{1.5}$} ~ K~ 
\raise4pt\hbox{$i^*(\psi+ \epsilon \psi^*) \atop \ra{3.6}$}~
L^* \to 0$$
is exact.
\hfill$\square$

An $\epsilon$-nonsingular split form $(K,\psi)$ admits a
(sub)lagrang\-ian if and only if the associated $\epsilon$-quadratic
form $(K,\lambda,\mu)$ admits a (sub)lagrang\-ian. (Sub)lagra\-ngians
in split $\epsilon$-quadratic forms are thus (sub)lagrang\-ians in
$\epsilon$-quadratic forms with the $(-\epsilon)$-quadratic structure
$\theta$, which (following Novikov) is sometimes called the ``hessian''
form.

\noindent {\sc Definition 5.6}
The $\epsilon$-nonsingular {\it hyperbolic} split form
$H_{\epsilon}(L)$ is given for any f.\ g. projective $\Lambda$-module
$L$ by
$$H_{\epsilon}(L)~=~(L\oplus L^*, \pmatrix{0 & 1 \ccr 0 & 0 \ccr}~:~
L\oplus L^* \to (L\oplus L^*)^*~=~L^* \oplus L)~,$$
with lagrang\-ian $(i=\pmatrix{1 \ccr 0},0):(L,0) \to H_{\epsilon}(L)$.
\hfill$\square$

\noindent {\sc Theorem 5.7}
{\it An $\epsilon$-nonsingular split form $(K,\psi)$ admits a
lagrangian if and only if it is $\epsilon$-quadratic isomorphic
isomorphic to the hyperbolic form $H_{\epsilon}(L)$. Moreover, the
inclusion $(i,\theta):(L,0)\to (K,\psi)$ of a lagrangian extends to an
$\epsilon$-quadratic isomorphism of split forms
$(f,\chi):H_{\epsilon}(L) \cong (K,\psi)$.}\hfil\break
{\sc Proof}\t :
An isomorphism of forms sends lagrangians to lagrangians, so any form 
isomorphic to a hyperbolic has at least one lagrangian.
Conversely suppose that $(K,\psi)$ has a lagrangian 
$(i,\theta):(L,0) \to (K,\psi)$.
An extension of $(i,\theta)$ to an $\epsilon$-quadratic isomorphism 
$(f,\chi):H_{\epsilon}(L)\cong (K,\psi)$
determines a lagrangian $f(L^*) \subset K$ complementary to $L$.
Construct an isomorphism $f$ by choosing a complementary lagrangian 
to $L$ in $(K,\psi)$.
Let $i\in {\rm Hom}_\Lambda (L,K)$ be the inclusion, and choose a splitting 
$j'\in {\rm Hom}_\Lambda (L^*,K)$ of $i^* (\psi+ \epsilon \psi^*) \in
{\rm Hom}_\Lambda (K,L^*)$, so that
$$i^* (\psi+ \epsilon \psi^*)j'~=~ 1 \in {\rm Hom}_\Lambda (L^*,L^*)~.$$
In general, $j':L^*\to K$ is not the inclusion of a lagrangian, with
$j'^* \psi j' \ne 0\in Q_{\epsilon}(L^*)$.
Given any $k\in {\rm Hom}_\Lambda (L^*,L)$ there is defined another splitting
$$j~=~j' + ik~:~L^* \to K$$
such that
$$\eqalign{j^* \psi j~&=~j^{\prime *} \psi j'+ k^*i^* \psi ik + k^*i^* \psi
j' + j^{\prime *} \psi ik\cr
&=~j^{\prime *} \psi j' + k \in Q_{\epsilon}(L^*)~.}$$ 
Choosing a representative $\psi\in {\rm Hom}_\Lambda(K,K^*)$
of $\psi\in Q_{\epsilon}(K)$ and setting
$$k~=~-j'^*\psi j'~:~L^* \to L^*$$
there is obtained a splitting $j:L^*\to K$ 
which is the inclusion of a lagrangian
$$(j,\nu)~:~(L^*,0)\to (K,\psi)~.$$
The isomorphism of $\epsilon$-quadratic forms
$$(i~ j)~=~\pmatrix{\theta & 0 \ccr
j^*\psi i & \psi}~:~H_{\epsilon}(L) ~\isa~ (K,\psi)$$
is an $\epsilon$-quadratic isomorphism of split forms.
\hfill$\square$

\noindent {\sc Remark 5.8}
Theorem 5.7 is a generalization of Witt's theorem on the extension 
to isomorphism of an isometry of quadratic forms over fields.
The procedure for modifying the choice of complement to be a lagrangian 
is a generalization of the Gram-Schmidt method of constructing orthonormal
bases in an inner product space. Ignoring the split structure 5.7 shows
that a nonsingular $\epsilon$-quadratic form admits a
lagrangian (in the sense of 5.1 (iii)) if and only if it is isomorphic to a
hyperbolic form.
\hfill$\square$

\noindent {\sc Corollary 5.9}
{\it For any $\epsilon$-nonsingular split form $(K,\psi)$ 
the diagonal inclusion
$$\Delta ~:~K \to K\oplus K~;~x \mapto (x,x)$$
extends to an $\epsilon$-quadratic isomorphism of split forms}
$$H_{\epsilon}(K) ~\isa~ (K,\psi) \oplus (K,-\psi)~.$$
{\sc Proof}\t :
Apply 5.7 to the inclusion of the lagrangian
$$(\Delta,0) ~:~ (K,0) \to (K\oplus K ,\psi\oplus -\psi)~.$$
(This result has already been used in 3.1).
\hfill$\square$

\noindent {\sc Proposition 5.10}
{\it The inclusion $(i,\theta):(L,0) \to (K,\psi)$ of a sublagrangian
in an $\epsilon$-nonsingular split
form $(K,\psi)$ extends to an isomorphism of forms}
$$(f,\chi)~:~ H_{\epsilon}(L) \oplus (L^{\bot}/L,[\psi])~\isa~(K,\psi)~.$$
{\sc Proof}\t : 
For any direct complement $L_1$ to $L^{\bot}$ in $K$ there is defined
a $\Lambda$-module isomorphism
$$e~:~L_1 ~\isa~ L^* ~;~ x \mapto (y \mapto (1+T_{\epsilon})\psi(x,y))~.$$
Define a $\Lambda$-module morphism
$$j~:~ L^*~ \raise4pt\hbox{$e^{-1} \atop \ra{2.5} $}~ L_1 ~
\raise4pt\hbox{${\rm inclusion} \atop \ra{3}$}~ K~.$$
The $\epsilon$-nonsingular split form defined by
$$(H, \phi) ~=~(L \oplus L^* ,
 \pmatrix{ 0 & 1 \ccr 0 & j^* \psi j})$$
has lagrangian $L$, so that it is isomorphic to
the hyperbolic form $H_{\epsilon}(L)$ by 5.7.
Also, there is defined an $\epsilon$-quadratic morphism of split forms
$$(g~=~(i~j)~,\pmatrix{\theta & i^* \psi j \ccr
0 & 0})~:~(H,\phi) \to (K,\psi)$$
with $g:H \to K$ an injection split by
$$h~=~((1+T_{\epsilon})\phi)^{-1}g^*(1+T_{\epsilon})\psi~:~K \to H~.$$
The direct summand of $K$ defined by
$$\eqalign{
H^{\bot}~&
=~\{x \in K \,\vert\, (1+T_{\epsilon})\psi(x,gy)~=~0~{\rm for~all~} y \in H\}\ccr
&=~{\rm ker}(g^*(1+T_{\epsilon})\psi:K \to H^*)~=~{\rm ker}(h:K \to H)}$$
is such that
$$K~=~g(H) \oplus H^{\bot}~.$$
It follows from the factorization
$$i^*(1+T_{\epsilon})\psi~:~ K \raise4pt\hbox{$h \atop \ra{1.5}$}~
H~=~L \oplus L^* ~ \raise4pt\hbox{${\rm projection} \atop \ra{3.5}$}~L^*$$
that
$$L^{\bot}~=~{\rm ker}(i^*(1+T_{\epsilon})\psi:K\to L^*)~=~L\oplus H^{\bot}~.$$
The restriction of $\psi\in S(K)$ to $H^{\bot}$ defines
an $\epsilon$-nonsingular split form $(H^{\bot},\phi^ {\bot})$.
The injection $g$ and the inclusion $g^{\bot}:H^{\bot} \to K$ are the
components of a $\Lambda$-module isomorphism 
$$f~=~(g~g^{\bot})~:~H\oplus H^{\bot} \to K$$
which defines an $\epsilon$-quadratic isomorphism of split forms
$$(f,\chi)~:~(H,\phi)\oplus (H^{\bot},\phi^{\bot})~\isa~(K,\psi)$$
with
$$(H^{\bot}, \phi^{\bot})~\cong~(L^{\bot}/L,[\psi])~.\eqno{\hbox{$\square$}}$$

\noindent {\sc Example 5.11}
An $n$-connected $(2n+1)$-dimensional normal bordism
$$((g,c);(f,b),(f',b'))~:~
(W^{2n+1};M^{2n},M'^{2n})\to X\times ([ 0,1];\{0\},\{1\})$$
is such that $W$ has a handle decomposition on $M$ of the type
$$W~=~M \times I \cup \bigcup_k n\hbox{\rm -handles}~D^n \times D^{n+1}
\cup \bigcup_{k'}(n+1)\hbox{\rm -handles}~ D^{n+1} \times D^n~.$$
Let
$$(W;M,M')~=~(W';M,M'') \cup_{M''} (W'';M'',M')$$ 
with
$$\eqalign{
&W'~=~M \times [ 0,1] \cup \bigcup_k ~n\hbox{\rm -handles}~ D^n \times D^{n+1}~,\ccr
&M''~=~{\rm cl.}(\partial W' \backslash M)~,\ccr
&W''~=~M'' \times [ 0,1] \cup \bigcup_{k'} ~(n+1)\hbox{\rm -handles}~ D^{n+1}
 \times D^n~.}$$
The restriction of $(g,c)$ to $M''$ is an $n$-connected $2n$-dimensional
normal map
$$(f'',b'')~:~M''~\cong~M \# (\#_k S^n \times S^n)~\cong~M'
\# (\#_{k'} S^n \times S^n) \to X$$
with kernel $(-1)^n$-quadratic form
$$\eqalign{(K_n(M''),\lambda'' , \mu'')~ & \cong ~
(K_n(M),\lambda , \mu) \oplus H_{(-1)^n}(\bbz[\pi_1(X)]^k) \ccr
& \cong~(K_n(M'),\lambda' , \mu')
\oplus H_{(-1)^n}(\bbz[\pi_1(X)]^{k'})~.}$$
Thus $(K_n(M''),\lambda'', \mu'')$ has sublagrangians
$$\eqalign{&L~=~{\rm im}(K_{n+1}(W',M'') \to K_n(M''))~\cong~\bbz[\pi_1(X)]^k~,\ccr
&L'~=~{\rm im}(K_{n+1}(W'',M'') \to K_n(M''))~\cong~\bbz[\pi_1(X)]^{k'}}$$
such that
$$\eqalign{&(L^{\bot}/L,[\lambda''],[\mu''])~\cong ~(K_n(M),\lambda , \mu)~,\ccr
&(L'^{\bot}/L',[\lambda'']' , [\mu'']')~\cong ~(K_n(M'),\lambda' , \mu')~.}$$
Note that $L$ is a lagrangian if and only if $(f,b):M \to X$ is a homotopy
equivalence. Similarly for $L'$ and $(f',b'):M' \to X$.
\hfill$\square$

\noindent {\sectionfont \S6. Short odd complexes}

A ``$(2n+1)$-complex'' is the algebraic structure best suited to
describing the surgery obstruction of an $n$-connected
$(2n+1)$-dimensional normal map. In essence it is a 1-dimensional 
chain complex with $(-1)^n$-quadratic Poincar\'e duality. 
 
As before, let $\Lambda$ be a ring with involution.

\noindent {\sc Definition 6.1}
A {\it $(2n+1)$-complex over $\Lambda$} $(C,\psi)$ is a f.\ g.
free $\Lambda$-module chain complex of the type
$$C~:~\dots \to 0 \to ~C_{n+1}~ 
\raise4pt\hbox{$d \atop \ra{1.5} $}~ C_n \to 0\to \dots $$
together with two $\Lambda$-module morphisms
$$\psi_0~:~C^n~=~(C_n)^* \to C_{n+1}~,~\psi_1~:~C^n \to C_n$$
such that
$$d\psi_0 + \psi_1 + (-1)^{n+1}{\psi^*_1}~=~ 0~:~C^n \to C_n~,$$
and such that the chain map
$$(1+T)\psi_0~:~C^{2n+1-*} \to C$$
defined by
$$\eqalign{&d_{C^{2n+1-*}}~=~(-1)^{n+1}d^*~:\ccr
& \indent\indent (C^{2n+1-*})_{n+1}~=~C^n \to (C^{2n+1-*})_n~=~C^{n+1}~,\ccr
&(1+T)\psi_0~= ~\cases{\psi_0~:~(C^{2n+1-*})_{n+1}~=~C^n \to C_{n+1}\ccr
\psi_0^*~:~(C^{2n+1-*})_n~=~C^{n+1} \to C_n~,} \ccr 
&(C^{2n+1-*})_r~=~C^{2n+1-r}~=~0~{\rm for}~r \ne n,n+1}$$
is a chain equivalence
$$\xymatrix@C+4pt{
C^{2n+1-*} : \ar[d]_{\displaystyle{(1+T)\psi_0}} \dots \ar[r] &
0 \ar[r] \ar[d] &
C^n \ar[r]^{\displaystyle{(-1)^{n+1}d^*}} \ar[d]_{\displaystyle{\psi_0}} & 
C^{n+1} \ar[r] \ar[d]_{\displaystyle{\psi_0^*}} & 0 \ar[d] \ar[r] & \dots \\
C~: \dots \ar[r] &0 \ar[r] &
C_{n+1} \ar[r]^{\displaystyle{d}} & C_n \ar[r] & 0 \ar[r] & \dots}
\eqno{\lower20pt\hbox{$\square$}}$$

\noindent {\sc Remark 6.2}
A $(2n+1)$-complex is essentially the inclusion of a lagrangian
in a hyperbolic split $(-1)^n$-quadratic form
$$(\pmatrix{\psi_0 \ccr d^* },-\psi_1)~:~(C^n,0) \to 
 H_{(-1)^n}(C_{n+1})~.$$
The chain map $(1+T)\psi_0 : C^{2n+1-*} \to C$
is a chain equivalence if and only if the algebraic mapping cone
$$0 \to C^n~
\raise8pt\hbox{$\pmatrix{\psi_0 \ccr d^* } \atop \ra{4}$}~
C_{n+1} \oplus C^{n+1}~ 
\raise4pt\hbox{$\pmatrix{d & (-1)^n{\psi^*_0}} \atop \ra{5.5} $}~
C_n \to 0$$
is contractible, which is just the lagrangian condition. The triple
$$\hbox{\rm (\ form\ ;\ lagrangian\ ,\ lagrangian\ )}~=~
(H_{(-1)^n}(C_{n+1});C_{n+1},{\rm im}\pmatrix{\psi_0 \ccr d^*})$$
is an example of a ``$(-1)^n$-quadratic formation''. 
Formations will be studied in greater detail in \S9 below. 
\hfill$\square$

\noindent {\sc Example 6.3}
Define a {\it presentation} of an $n$-connected $(2n+1)$-dimensional normal map 
$(f,b):M^{2n+1}\to X$ to be a normal bordism
$$((g,c);(f,b),(f',b'))~:~
(W^{2n+2};M^{2n+1},M^{\prime 2n+1}) \to X \times ([ 0,1];\{0\},\{1\})$$
such that $W\to X\times [ 0,1]$ is $n$-connected, with
$$K_r(W)~=~0~{\rm for}~r \ne n+1 ~.$$
Then $K_{n+1}(W)$ a f.\ g. free $\bbz[\pi_1(X)]$-module and
$W$ has a handle decomposition on $M$ of the type
$$W~=~M \times I \cup \bigcup_k (n+1)\hbox{\rm -handles}~ D^{n+1} \times D^{n+1}~,$$
and $K_{n+1}(W,M)\cong\bbz[\pi_1(X)]^k$ is a f.\ g. free 
$\bbz[\pi_1(X)]$-module with rank the number $k$ of $(n+1)$-handles.
Thus $(W;M,M')$ is the trace of surgeries on $k$ disjoint 
embeddings $S^n \times D^{n+1}\hookrightarrow M^{2n+1}$ with null-homotopy
in $X$ representing a set of $\bbz[\pi_1(X)]$-module generators of $K_n(M)$~.
For every $n$-connected $(2n+1)$-dimensional normal map
$(f,b):M^{2n+1}\to X$ the kernel $\bbz[\pi_1(X)]$-module $K_n(M)$ is f.\ g.,
so that there exists a presentation 
$(g,c):(W;M,M')\to X \times ([ 0,1];\{0\},\{1\})$.
Poincar\'e duality and the universal coefficient theorem give natural 
identifications of f.\ g. free $\bbz[\pi_1(X)]$-modules
$$\eqalign{
&K_{n+1}(W)~=~K^{n+1}(W,\partial W)~=~K_{n+1}(W,\partial W)^*~~
(\partial W=M \cup M')~,\ccr
&K_{n+1}(W,M) ~=~ K^{n+1}(W,M') ~=~ K_{n+1}(W,M')^*~.}$$
The presentation determines a $(2n+1)$-complex $(C,\psi)$ such that
$$H_*(C) ~=~ K_*(M)~,$$
with
$$\eqalign{&d~=~({\rm inclusion})_*~:~C_{n+1}~=~K_{n+1}(W,M')~
\to C_n~=~K_{n+1}(W,\partial W)~,\ccr
&\psi_0~=~({\rm inclusion})_*~:~C^n~=~K_{n+1}(W) \to~C_{n+1}~=~K_{n+1}(W,M')~.}$$
The hessian $(C^n,-\psi_1 \in Q_{(-1)^{n+1}}(C^n))$ is the 
geometric self-intersection $(-1)^{n+1}$-quadratic form on the kernel 
$C^n=K_{n+1}(W)$ of the normal map $W^{2n+2}\to X\times [ 0,1]$, such that
$$\eqalign{&-(\psi_1+(-1)^{n+1}{\psi^*_1})~=~d\psi_0~=~{\rm inclusion}_*~:\ccr
&\hskip25pt
C^n~=~K_{n+1}(W) \to C_n~=~K_{n+1}(W,\partial W)~=~K_{n+1}(W)^*~.}$$
The chain equivalence $(1+T)\psi_0:C^{2n+1-*}\to C$ induces the 
Poincar\'e duality isomorphisms
$$[M]\cap -~:~H^{2n+1-*}(C)~=~K^{2n+1-*}(M)~\isa~H_*(C)~=~K_*(M)~.\eqno{\square}$$

\noindent {\sc Remark 6.4}
The $(2n+1)$-complex $(C,\psi)$ of 6.3 can also be obtained by working
inside $M$, assuming that $X$ has a single $(2n+1)$-cell
$$X~=~X_0 \cup D^{2n+1}$$
(as is possible by the Poincar\'e disc theorem of Wall [27])
so that there is defined a degree 1 map
$${\rm collapse}~:~X \to X/X_0 ~=~ S^{2n+1}~.$$
Let $U \subset M^{2n+1}$ be the disjoint union of the $k$ embeddings
$S^n \times D^{n+1} \hookrightarrow M$ with null-homotopies in $X$, so
that $(f,b)$ has a Heegaard splitting as a union of normal maps
$$\eqalign{
&(f,b)~=~(e,a) \cup (f_0,b_0)~:\ccr
&M~=~(U,\partial U) \cup (M_0,\partial M_0)
\to X~=~(D^{2n+1},S^{2n}) \cup (X_0,\partial X_0)}$$
with the inclusion (6.2) of the lagrangian
$$\pmatrix{\psi_0 \ccr d^* }~:~C^n \to C_{n+1}\oplus C^{n+1}$$
in the hyperbolic $(-1)^n$-quadratic form $H_{(-1)^n}(C_{n+1})$ given by
$${\rm inclusion}_* ~:~ K_{n+1}(M_0, \partial U) \to K_n( \partial U)~
=~ K_{n+1}(U, \partial U) \oplus K_n(U)~.$$
Wall obtained the surgery obstruction of $(f,b)$ using an extension (cf. 5.7)
of this inclusion to an automorphism
$$\alpha ~:~ H_{(-1)^n}(C_{n+1}) ~ \isa ~ H_{(-1)^n}(C_{n+1})~,$$
which will be discussed further in \S10 below. The presentation of $(f,b)$
used to obtain $(C,\psi)$ in 6.3 is the trace of the $k$ surgeries
on $U \subset M$
$$\eqalign{&(g,c)~=~(e_1,a_1) \cup (f_0,b_0) \times {\rm id}~:\ccr
&(W;M,M')~=~(V;U,U') \cup M_0 \times ([ 0,1];\{0\},\{1\}) \to
X \times ([ 0,1];\{0\},\{1\})}$$
with
$$(V;U,U') ~=~ \bigcup_k 
(D^{n+1} \times D^{n+1};S^n \times D^{n+1},D^{n+1} \times S^n)~.\eqno{\square}$$

\noindent {\sc Example 6.5}
There is also a relative version of 6.3. A presentation of an
$n$-connected normal map $(f,b):M^{2n}\to X$ from a
$(2n+1)$-dimensional manifold with boundary $(M,\partial M)$ to a
geometric Poincar\'e pair $(X,\partial X)$ with 
$\partial f =f\vert:\partial M \to \partial X$ a homotopy equivalence 
is a normal map of triads
$$\eqalign{&(W^{2n+2};M^{2n+1},M^{\prime 2n+1};\partial M\times [ 0,1])\ccr
&\hskip50pt
\to (X\times [ 0,1];X\times \{0\},X\times \{1\};\partial X\times [ 0,1])}$$
such that $W\to X \times [ 0,1]$ is $n$-connected.
Again, the presentation determines a $(2n+1)$-complex $(C,\psi)$ over 
$\bbz[\pi_1(X)]$ with
$$C_n~=~K_{n+1}(W,\partial W)~,~C_{n+1}~=~K_{n+1}(W,M')~,~
H_*(C)~=~K_*(M)~.\eqno{\square}$$

\noindent {\sc Remark 6.6}
(Realization of odd-dimensional surgery obstructions, Wall [28,\t 6.5])
The theorem of [28] realizing automorphisms of hyperbolic forms as
odd-dimensional surgery obstructions has the following interpretation 
in terms of complexes.
Let $(C,\psi)$ be a $(2n+1)$-complex over $\bbz[\pi]$, with  $\pi$ a
finitely presented group. 
Let $n \geq 2$, so that there exists a $2n$-dimensional manifold $X^{2n}$ 
with $\pi_1(X)=\pi$. For any such $n \geq 2$, $X$ there exists an 
$n$-connected $(2n+1)$-dimensional normal map 
$$(f,b)~:~(M^{2n+1};\partial_-M,\partial_+M)\to 
X^{2n} \times ([ 0,1];\{0\},\{1\})$$
with $\partial_-M=X\to X$ the identity and 
$\partial_+M\to X$ a homotopy equivalence, and with a presentation 
with respect to which $(f,b)$ has kernel $(2n+1)$-complex $(C,\psi)$.
Such a normal map is constructed from the identity $X\to X$ in two stages.
First, choose a basis $\{b_1,b_2,\dots,b_k\}$ for $C_{n+1}$, and perform 
surgeries on $k$ disjoint trivial embeddings 
$S^{n-1} \times D^{n+1} \hookrightarrow X^{2n}$ 
with trace
$$\eqalign{&(U;X,\partial_+U)~=~
(X\times [ 0,1]\cup \bigcup_k D^n \times D^{n+1};X\times \{0\},
X\# \#_k S^n \times S^n)\ccr
&\hskip150pt \to X \times ([ 0,1/2];\{0\},\{1/2\})~.}$$
The $n$-connected $2n$-dimensional normal map
$\partial_+U \to X \times \{1/2\}$ has kernel $(-1)^n$-quadratic form
$$(K_n(\partial_+ U), \lambda, \mu)~=~H_{(-1)^n}(\bbz[\pi]^k)~=~
H_{(-1)^n}(C_{n+1})~.$$
Second, choose a basis $\{c_1,c_2,\dots,c_k\}$ for $C^n$ and realize
the inclusion of the lagrangian in $H_{(-1)^n}(C_{n+1})$ by surgeries
on $k$ disjoint embeddings $S^n \times D^n \hookrightarrow \partial_+U$ 
with trace
$$(M_0;\partial_+ U, \partial_+ M) \to X\times ([ 0,1];\{0\},\{1\})$$
such that
$$\eqalign{\pmatrix{\psi_0 \ccr d^* }~
&=~\partial~:~C^n~=~K_{n+1}(M_0,\partial_+U)\ccr
&\to C_{n+1} \oplus C^{n+1}~=~
K_{n+1}(U,\partial_+U)\oplus K_n(U)~=~K_n (\partial_+U)~.}$$
The required $(2n+1)$-dimensional normal map realizing $(C,\psi)$ 
is the union
$$(M;\partial_-M,\partial_+M)~=~
(U;X,\partial_+U) \cup (M_0;\partial_+U,\partial_+M)
 \to X \times ([ 0,1];\{0\},\{1\})~.$$
The corresponding presentation is the trace of surgeries on $k$ disjoint 
embeddings $S^n \times D^{n+1} \hookrightarrow U\subset M^{2n+1}$.
This is the terminology (and result) of Wall [28,\t Chapter 6].
\hfill$\square$

The choice of presentation (6.3) for an $n$-connected 
$(2n+1)$-dimensional normal map $(f,b):M^{2n+1} \to X$
does not change the ``homotopy type'' of the associated $(2n+1)$-complex
$(C,\psi)$, in the following sense.

\noindent {\sc Definition 6.7}
(i) A {\it map} of $(2n+1)$-complexes over $\Lambda$
$$f~:~(C,\psi) \to (C',\psi')$$
is a chain map $f:C\to C'$ such that there exist $\Lambda$-module morphisms
$$\chi_0~:~C^{\prime n+1} \to C^{\prime}_{n+1}~~,~~
\chi_1~:~C^{\prime n} \to C^{\prime}_n$$
with
$$\eqalign{&f\psi_0f^* - \psi^{\prime}_0~=~(\chi_0 + (-1)^{n+1} \chi^*_0)
d^{\prime *}~:~C^{\prime n} \to C^{\prime}_{n+1}~, \ccr
&f\psi_1 f^* - \psi^{\prime}_1~=~ -d' \chi_0 d^{\prime *} +
 \chi_1 + (-1)^n \chi^*_1~:~C^{\prime n} \to C^{\prime}_n~.}$$
(ii) A {\it homotopy equivalence}
of $(2n+1)$-complexes is a map with $f:C\to C'$ a chain equivalence.
\hfil \break
(iii) An {\it isomorphism}
of $(2n+1)$-complexes is a map with $f:C\to C'$ an isomorphism of
 chain complexes.
\hfill$\square$

\noindent {\sc Proposition 6.8}
{\it Homotopy equivalence is an equivalence relation on $(2n+1)$-complexes.}
\hfil\break
{\sc Proof}\t :
For $m\ge 0$ let $E(m)$ be the contractible f.\ g.
 free $\Lambda$-module chain complex defined by
$$\eqalign{&d_{E(m)}~=~1~:~E(m)_{n+1}~= ~\Lambda^m  \to E(m)_n~=~\Lambda^m ~,\ccr
&E(m)_r~=~ 0 ~{\rm for}~r \ne n,n+1~.}$$
A map $f:(C,\psi)\to (C',\psi')$ is a homotopy equivalence if and only if
for some $m,m' \ge 0$ there exists an isomorphism
$$f'~:~(C,\psi) \oplus (E(m),0) ~ \isa ~ (C',\psi') \oplus (E(m'),0)$$
such that the underlying chain map $f'$ is chain homotopic to
$$f\oplus 0~:~C\oplus E(m) \to C' \oplus E(m')~.$$
Isomorphism is an equivalence relation on
$(2n+1)$-complexes, and hence so is homotopy equivalence.
\hfill$\square$

\noindent {\sc Example 6.9}
The $(2n+1)$-complexes $(C,\psi)$, $(C',\psi')$ associated by 6.3
to any two presentations 
$$(W;M,\widehat M)\to X \times ([ 0,1];\{0\},\{1\})~,~ 
(W';M,\widehat M') \to X\times ([ 0,1];\{0\},\{1\})$$
of an $n$-connected normal map $M^{2n+1}\to X$ are homotopy equivalent.
 Without loss of generality it may be assumed that $W$ and $W'$ are the traces
of surgeries on disjoint embeddings 
$$g^i~:~S^n \times D^{n+1} \hookrightarrow M~,~
g^{\prime j}~:~S^n \times D^{n+1}\hookrightarrow M~,$$
corresponding to two sets of $\bbz[\pi_1(X)]$-module generators of 
$K_n(M)$. Define a presentation of $M\to X$
$$(W'';M,M'')~=~(W;M,\widehat M) \cup (V;\widehat M,M'')~=~
(W';M,\widehat M') \cup (V';\widehat M',M'')$$
with $(V;\widehat M,M'')$ the presentation of $\widehat M\to X$
defined by the trace of the surgeries on the copies 
$\widehat g^{\prime j}:S^n\times D^{n+1} \hookrightarrow \widehat M$
of $g^{\prime j}:S^n\times D^{n+1} \hookrightarrow M$, and $(V';\widehat M',M'')$ the presentation of 
$\widehat M'\to X$ defined by the trace of the surgeries on the copies 
$\widehat g^i:S^n\times D^{n+1}\hookrightarrow \widehat M'$ of 
$g^i:S^n\times D^{n+1}\hookrightarrow M$.

$$\beginpicture
\setcoordinatesystem units <5pt,3pt> 
\setlinear
\plot -20 -10 -20 10 20 10 20 -10 -20 -10 /
\plot 0 10 0 -10 /
\plot -20 -20 -20 -40 20 -40 20 -20 -20 -20 /
\plot 0 -20 0 -40 /
\put{$M$} at -22 0
\put{$W$} at -10 0
\put{$\widehat{M}$} at 2 0
\put{$V$} at 10 0
\put{$M''$} at 22 0
\put{$M$} at -22 -30
\put{$W'$} at -10 -30
\put{$\widehat{M}'$} at 2 -30
\put{$V'$} at 10 -30
\put{$M''$} at 22 -30
\put{$W'' = W\cup_{\widehat M}V = W' \cup_{\widehat M'}V'$} at 0 -15
\endpicture$$

\vskip10mm

\noindent The projections $C'' \to C$, $C'' \to C'$
define homotopy equivalences of $(2n+1)$-complexes
$$(C'',\psi'') \to (C,\psi)~~,~~ (C'',\psi'') \to (C',\psi')~.\eqno{\square}$$

\noindent {\sc Definition 6.10}
A $(2n+1)$-complex $(C,\psi)$ over $\Lambda$ is {\it contractible}
if it is homotopy equivalent to the zero complex $(0,0)$,
or equivalently if $d:C_{n+1} \to C_n$ is a $\Lambda$-module isomorphism.
\hfill$\square$

\noindent {\sc Example 6.11}
A $(2n+1)$-complex $(C, \psi)$ associated to an $n$-connected
$(2n+1)$-dimensional normal map $(f,b):M^{2n+1} \to X$
is contractible if (and for $n \ge 2$ only if) $f$ is a
homotopy equivalence, by the theorem of J.H.C. Whitehead.
The $(2n+1)$-complexes $(C,\psi)$ associated to 
the various presentations of a homotopy equivalence
$(f,b):M^{2n+1} \to X$ are contractible, by 6.9. 
The zero complex $(0,0)$ is associated to the presentation
$$(f,b) \times {\rm id.} ~:~ M \times ([ 0,1];\{0\},\{1\})
 \to X \times ([ 0,1];\{0\},\{1\})~.\eqno{\square}$$

\noindent {\sectionfont \S7. Complex cobordism}

The cobordism of $(2n+1)$-complexes is the equivalence relation which
corresponds to the normal bordism of $n$-connected
$(2n+1)$-dimensional normal maps. The $(2n+1)$-dimensional surgery
obstruction group $L_{2n+1}(\Lambda)$ will be defined in \S8 below to
be the cobordism group of $(2n+1)$-complexes over $\Lambda$.

\noindent {\sc Definition 7.1}
A {\it cobordism}
of $(2n+1)$-complexes $(C,\psi)$, $(C',\psi')$ 
$$((j~ j'):C\oplus C'\to D,(\delta \psi,\psi\oplus - \psi'))$$
is a f.\ g. free $\Lambda$-module chain complex of the type
$$D ~:~ \dots \to 0 \to D_{n+1} \to 0 \to \dots $$
together with $\Lambda$-module morphisms
$$\displaylines{j~:~C_{n+1} \to D_{n+1}~~,~~j'~:~C'_{n+1} \to D_{n+1}~,\ccr
\delta \psi_0~:~D^{n+1}~=~(D_{n+1})^* \to D_{n+1}}$$
such that the {\it duality} $\Lambda$-module chain map
$$(1+T)(\delta \psi_0, \psi_0 \oplus -\psi'_0)~:~
{\scr C}(j')^{2n+2-*} \to {\scr C}(j)$$
defined by 
$$\eqalign{&
\hbox{$(1+T)(\delta \psi_0, \psi_0 \oplus -\psi'_0)~=~
\pmatrix{\delta \psi_0 + (-1)^{n+1} \delta \psi_0^* & j' \psi'_0 \ccr
\psi^*_0 j^* & 0}$}\ccr
& \indent :~{\scr C}(j')^{2n+1}~=D^{n+1}\oplus C'^n \to {\scr C}(j)_{n+1}~=~D_{n+1} 
\oplus C_n}$$
is a chain equivalence, with ${\scr C}(j)$, ${\scr C}(j')$ the algebraic
mapping cones of the chain maps $j : C \to D$, $j' : C' \to D$.
\hfill$\square$

The duality chain map ${\scr C}(j')^{2n+2-*} \to {\scr C}(j)$ is given by 
$$\xymatrix@C-10pt{
{\scr C}(j')^{2n+2-*} : 
\ar[d]_{(1+T)\delta\psi_0}\dots \ar[r] 
&0 \ar[r] \ar[d] &
0 \ar[r] \ar[d] & D^{n+1}\oplus C^{\prime n} \ar[r] 
\ar[d] & C^{\prime n+1}
\ar[r] \ar[d] & 0 \ar[r] \ar[d] & \dots \\
{\scr C}(j): \dots \ar[r] &0 \ar[r] &
C_{n+1} \ar[r] & D_{n+1}\oplus C_n \ar[r] & 0
\ar[r] & 0 \ar[r] & \dots}$$
The condition for it to be a chain equivalence is just that the 
$\Lambda$-module morphism
$$\eqalign{&\pmatrix{d & 0 & \psi^*_0 j^* \ccr
0 & d'^* & j'^* \ccr
(-1)^{n+1}j & j'\psi'_0 & \delta \psi_0 + (-1)^{n+1}\delta \psi^*_0 }\ccr
&\hskip50pt :~C_{n+1}\oplus C'^n \oplus D^{n+1} \to 
C_n \oplus C'^{n+1} \oplus D_{n+1} }$$
be an isomorphism.

\noindent {\sc Example 7.2}
Suppose given two $n$-connected $(2n+1)$-dimensional normal maps
$M^{2n+1}\to X$, $M'^{2n+1}\to X$ with presentations (6.3)
$$\eqalign{&(W^{2n+2};M^{2n+1},\widehat M^{2n+1}) \to X\times ([ 0,1];\{0\},\{1\})~,\ccr
&(W'^{2n+2};M'^{2n+1},\widehat M'^{2n+1}) \to X\times ([ 0,1];\{0\},\{1\})}$$
and corresponding $(2n+1)$-complexes $(C,\psi)$, $(C',\psi')$.
 An $n$-connected normal bordism
$$(V^{2n+2};M^{2n+1},M^{\prime 2n+1}) \to X\times 
([ 0,1];\{0\},\{1\})$$
determines a cobordism $((j~ j'):C\oplus C'\to D,(\delta \psi,\psi
\oplus - \psi'))$ (again, up to some choices) from 
$(C,\psi)$ to $(C',\psi')$. Define an $n$-connected normal bordism
$$\eqalign{&(V'; \widehat M, \widehat M')~=~
(W;\widehat M,M) \cup (V;M,M') \cup (W';M',\widehat M')\ccr
&\hskip100pt \to X \times ([ 0,1];\{0\},\{1\})~.}$$
The exact sequence of 
stably f.\ g. free $\bbz[\pi_1(X)]$-modules
$$\eqalign{ 0 \to K_{n+1}(V)~ &\to K_{n+1}(V',\partial V')\ccr
&\to K_{n+1}(W,\partial W) \oplus K_{n+1}(W',\partial W') \to 0}$$
splits. Choosing any splitting $K_{n+1}(V', \partial V') \to K_{n+1}(V)$
define $j$, $j'$ by
$$\eqalign{&(j~j')~:~C_{n+1} \oplus C'_{n+1}~=~K_{n+1}(W,\widehat M)
\oplus K_{n+1}(W',\widehat M')\ccr
&\hskip50pt 
\raise4pt\hbox{$ {\rm incl}_* \oplus {\rm incl}_* \atop \ra{5}$} ~
K_{n+1}(V',\partial V')~\to K_{n+1}(V)~=~D_{n+1}~.}$$
Geometric intersection numbers provide a $(-1)^{n+1}$-quadratic form
$(D^{n+1},\break\delta \psi_0)$ over $\bbz [\pi_1(X)]$ such that
the duality chain map ${\scr C}(j')^{2n+2-*} \to {\scr C}(j)$ 
is a chain equivalence inducing the Poincar\'e duality isomorphisms
$$[V]\cap -~:~H^{2n+2-*}(j')~=~K^{2n+2-*}(V,M')~ \isa ~H_*(j)~=~K_*(V,M)~.\eqno{\square}$$ 

\noindent {\sc Definition 7.3}
A {\it null-cobordism} of a $(2n+1)$-complex $(C,\psi)$ is a cobordism
$(j:C \to D, (\delta \psi, \psi))$ to $(0,0)$.
\hfill$\square$

\noindent {\sc Example 7.4}
Let $(W^{2n+2};M'^{2n+1},M^{2n+1}) \to X\times 
([ 0,1];\{0\},\{1\})$ be a presentation of an $n$-connected 
$(2n+1)$-dimensional normal map $M\to X$, with $(2n+1)$-complex
$(C,\psi)$.
For $n\ge 2$ there is a one-one correspondence between $n$-connected 
normal bordisms of $M \to X$
$$(V^{2n+2};M^{2n+1},N^{2n+1}) \to X \times ([ 0,1];\{0\},\{1\})$$
to homotopy equivalences $N\to X$ and null-cobordisms 
$(j:C\to D,\allowbreak(\delta \psi,\allowbreak \psi))$.
 (Every normal bordism of $n$-connected $(2n+1)$-dimensional normal 
 maps can be made $n$-connected by surgery below the middle dimension on 
 the interior).

$$
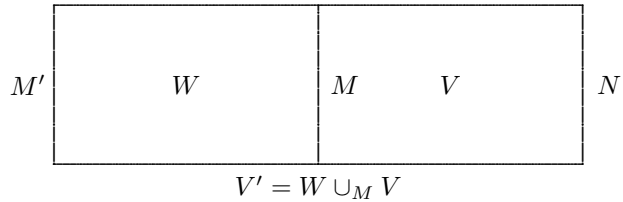
\beginpicture
\setcoordinatesystem units <5pt,3pt> 
\setlinear
\plot -20 -10 -20 10 20 10 20 -10 -20 -10 /
\plot 0 10 0 -10 /
\put{$M'$} at -22 0
\put{$W$} at -10 0
\put{$M$} at 2 0
\put{$V$} at 10 0
\put{$N$} at 22 0
\put{$V'=W\cup_MV$} at 0 -13
\endpicture$$

\noindent Any such $(V;M,N)\to X\times ([ 0,1];\{0\},\{1\})$ determines 
by 7.2 a null-cob\-ord\-ism $(j:C\to D, (\delta \psi,\psi))$ of $(C,\psi)$.
\hfill$\square$

Cobordisms of $(2n+1)$-complexes arise in the following way:
 
\noindent {\sc Construction 7.5}
An isomorphism of hyperbolic split $(-1)^n$-quadratic forms over $\Lambda$
$$\left(\pmatrix{\gamma & \widetilde \gamma \ccr \mu & \widetilde \mu}
,\pmatrix{\theta & 0 \ccr \widetilde \gamma^* \mu & \widetilde \theta}\right)~
:~ H_{(-1)^n}(G) ~ \isa ~ H_{(-1)^n}(F)$$
with $F,G$ f.\ g. free determines a cobordism of $(2n+1)$-complexes
$$((j ~j'):C \oplus C' \to D,(\delta \psi, \psi\oplus -\psi'))$$ 
by
$$\eqalign{
&d~=~\mu^*~:~C_{n+1}~=~F \to C_n~=~G^*~,\ccr
&\psi_0~=~\gamma~:~C^n~=~G \to C_{n+1}~=~F~,\ccr
&\psi_1~=~-\theta~:~C^n~=~G \to C_n ~=~G^*~,\ccr
&j~=~\widetilde \mu^* ~:~ C_{n+1}~=~F \to D_{n+1} ~=~G~,\ccr 
&d'~=~\gamma^*~:~C'_{n+1}~=~F^* \to C_n~=~G^*~,\ccr
&\psi'_0~=~\mu~:~C'^n~=~G \to C'_{n+1}~=~F^*~,\ccr
&\psi'_1~=~-\theta~:~C'^n~=~G \to C'_n ~=~G^*~,\ccr
&j'~=~\widetilde \gamma^* ~:~ C'_{n+1}~=~F^* \to D_{n+1} ~=~G~,\ccr 
&\delta \psi_0 ~=~0~:~D^{n+1}~=~G^* \to D_{n+1}~=~G~.}$$
\hfill$\square$

It can be shown that every cobordism of $(2n+1)$-complexes is
homotopy equivalent to one constructed as in 7.5. 

\noindent {\sc Example 7.6}
An $n$-connected $(2n+2)$-dimensional normal bordism
$$((g,c);(f,b),(f',b'))~:~(W^{2n+2};M^{2n+1},M'^{2n+1}) \to 
 X \times ([ 0,1];\{0\},\{1\})$$
with $g:W \to X \times [ 0,1]$ $n$-connected can be regarded both
as a presentation of $(f,b)$ and as a presentation of $(f',b')$. 
The cobordism of $(2n+1)$-complexes
$((j~j'):C \oplus C' \to D, (\delta \psi, \psi\oplus - \psi'))$
obtained in 7.2 with $W=V=W'$, $\widehat M = M'$, $\widehat M'=M'$
is the construction of 7.5 for an extension of the inclusion
of the lagrangian (6.2)
$$(\pmatrix{\gamma \ccr \mu},\theta)~=~
(\pmatrix{\psi_0 \ccr d^* },-\psi_1)~:~(C^n,0) \to 
H_{(-1)^n}(C_{n+1})$$
to an isomorphism of hyperbolic split $(-1)^n$-quadratic forms
$$\left(\pmatrix{\gamma & \widetilde \gamma \ccr \mu & \widetilde \mu}
,\pmatrix{\theta & 0 \ccr \widetilde \gamma^* \mu & \widetilde \theta}\right)~
:~ H_{(-1)^n}(C^n) ~ \isa ~ H_{(-1)^n}(C_{n+1})~,$$
with
$$\eqalign{&j~=~\widetilde \mu^* ~:~C_{n+1}~=~K_{n+1}(W,M')~
\to D_{n+1}=K_{n+1}(W)~, \ccr
&j' ~=~ \widetilde
\gamma^* ~:~ C'_{n+1}~=~K_{n+1}(W,M)\to D_{n+1}=K_{n+1}(W)~.}
\eqno{\lower10pt\hbox{$\square$}}$$

\noindent {\sc Remark 7.7}
Fix a $(2n+1)$-dimensional geometric Poincar\'e complex $X$
with reducible Spivak normal fibration, and choose
a stable vector bundle $\nu_X:X\to BO$ in the Spivak normal
class, e.g. a manifold with the stable normal bundle.
Consider the set of $n$-connected normal maps\break 
$(f:M^{2n+1}\to X,b:\nu_M\to \nu_X)$.
The relation defined on this set by
$$\eqalign{(M\to X) \sim (M'\to X)~&\hbox{\rm if there exists an 
$(n+1)$-connected normal}\cr
&{\rm bordism}~(W;M,M') \to X\times ([ 0,1];\{0\},\{1\})}$$
is an equivalence relation.
 Symmetry and transitivity are verified in the same way as for any geometric 
 cobordism relation. For reflexivity form the cartesian product 
 of an $n$-connected normal map $M^{2n+1}\to X$ with 
 $([ 0,1];\{0\},\{1\})$, as usual. 
 The product is an $n$-connected normal bordism
$$M\times ([ 0,1];\{0\},\{1\}) \to X\times ([ 0,1];\{0\},\{1\})$$
which can be made $(n+1)$-connected by surgery killing the $n$-dimensional
kernel $K_n(M \times [ 0,1])=K_n(M)$.
 The following verification that the cobordism of $(2n+1)$-complexes 
 is an equivalence relation uses algebraic surgery in exactly the same way.
\hfill$\square$

\noindent {\sc Proposition 7.8}
{\it Cobordism is an equivalence relation on $(2n+1)$-com\-plexes 
$(C,\psi)$ over $\Lambda$, such that $(C,\psi) \oplus (C,- \psi)$ is null-cobordant.
Homotopy equivalent complexes are cobordant.}
\hfil\break
{\sc Proof}\t :
Symmetry is easy: if 
$((j~j'):C\oplus C'\to D,(\delta \psi,\psi\oplus - \psi'))$ 
is a cobordism from $(C,\psi)$ to $(C',\psi')$ then
$$((j'~j):C' \oplus C\to D',(- \delta \psi,\psi' \oplus - \psi))$$ 
 is a cobordism from $(C',\psi')$ to $(C,\psi)$.
 For transitivity, suppose given adjoining cobordisms of $(2n+1)$-complexes
$$\eqalign{
&((j~j'):C\oplus C'\to D,(\delta \psi,\psi\oplus - \psi'))~, \ccr
&((\widetilde j'~j''):C' \oplus C''\to D',(\delta \psi',\psi'\oplus
- \psi''))~.}
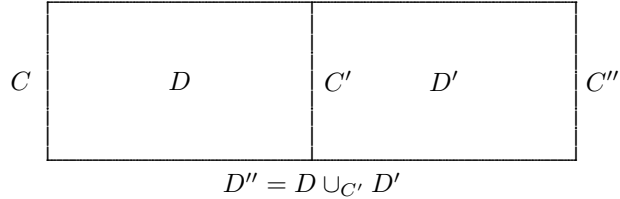
$$

$$\beginpicture
\setcoordinatesystem units <5pt,3pt> 
\setlinear
\plot -20 -10 -20 10 20 10 20 -10 -20 -10 /
\plot 0 10 0 -10 /
\put{$C$} at -22 0
\put{$D$} at -10 0
\put{$C'$} at 2 0
\put{$D'$} at 10 0
\put{$C''$} at 22 0
\put{$D''=D\cup_{C'}D'$} at 0 -13
\endpicture$$

\noindent Define the {\it union} cobordism between $(C,\psi)$ and
$(C'',\psi'')$
$$((\widetilde j ~ \widetilde j''):C \oplus C'' \to D'' ,
(\delta \psi'',\psi\oplus - \psi''))$$
by
$$\eqalign{
&D^{\prime \prime}_{n+1}~=~{\rm coker}(i= \pmatrix{j' \ccr d' \ccr \widetilde j'}
: C^{\prime}_{n+1} \to D_{n+1} \oplus C^{\prime}_n \oplus 
D^{\prime}_{n+1})~,\ccr
&\widetilde j~=~[j\oplus 0\oplus 0]~:~C_{n+1} \to 
D^{\prime \prime}_{n+1} ,\ccr
&\widetilde j''~=~[ 0\oplus 0 \oplus j'']~:~C^{\prime \prime}_{n+1} \to 
D^{\prime \prime}_{n+1}~,\ccr
&\delta \psi^{\prime \prime}_0~= ~
\left [\matrix{\delta \psi_0 & 0 & 0 \ccr
0 & 0 & 0 \ccr
0 & 0 & \delta \psi^{\prime}_0 }\right ]~:~D^{\prime \prime n+1}
 \to D^{\prime \prime}_{n+1}~.}$$
The $\Lambda$-module morphism $i:C^{\prime}_{n+1}\to D_{n+1} \oplus 
C^{\prime}_n \oplus D^{\prime}_{n+1}$ is a split injection since the dual
$\Lambda$-module morphism $i^*$ is a surjection, as follows from the 
Mayer-Vietoris exact sequence
$$H^{n+2}(D,C') \oplus H^{n+2}(D',C')~=~0 \oplus 0 \to 
H^{n+2}(D'') \to H^{n+2}(C')~=~0~.$$
\indent	Given any $(2n+1)$-complex $(C,\psi)$ let $(C',\psi')$ be the 
$(2n+1)$-complex defined by
$$\eqalign{&d'~=~(-1)^n \psi^*_0~:~C^{\prime}_{n+1}~=~C^{n+1} \to 
C^{\prime}_n~=~C_n~,\ccr
&\psi^{\prime}_0~=~d^*~:~C^{\prime n}~=~C^n \to 
C^{\prime}_{n+1}~=~C^{n+1}~,\ccr
&\psi'_1~=~ -\psi_1~:~C^{\prime n}~=~C^n \to C^{\prime}_n~=~
C_n~.}$$
Apply 5.7 to extend the inclusion of the lagrangian in $H_{(-1)^n}(C_{n+1})$
$$\pmatrix{\psi_0 \ccr d^*}~:C^n \to C_{n+1} \oplus C^{n+1}$$
to an isomorphism of $(-1)^n$-quadratic forms
$$\pmatrix{\psi_0 & \widetilde \psi_0 \ccr d^* & \widetilde d^*}~:~
 H_{(-1)^n}(C^n) ~ \isa ~ H_{(-1)^n}(C_{n+1})$$
with $\widetilde \psi_0 \in {\rm Hom}_\Lambda (C_n,C_{n+1})$, 
$\widetilde d \in {\rm Hom}_\Lambda (C_{n+1},C^n)$.
 Now apply 7.5 to construct from any such extension a cobordism 
$$((j~j'):C\oplus C'\to D,(\delta \psi,\psi\oplus - \psi'))$$ 
with
$$\eqalign{&j~=~\widetilde d~:~C_{n+1} \to D_{n+1}~=~C^n~,\ccr
&j'~=~\widetilde \psi^*_0~:~ C^{\prime}_{n+1}~=~C^{n+1} \to 
D_{n+1}~=~C^n~,\ccr
&d'~=~\psi^*_0~:~ C^{\prime}_{n+1}~=~C^{n+1} \to 
C^{\prime}_n~=~C_n~,\ccr
&\psi'_0 ~=~ d^* ~:~C'^n~=~C^n \to C'^{n+1}~=~C^{n+1}~,\ccr
&\delta \psi_0 ~=~0~:~ D^{n+1}~=~C_n \to D_{n+1}~=~C^n~.}$$
(This is the algebraic analogue of the construction of a presentation (6.3) 
$$(W^{2n+2};M^{2n+1},M'^{2n+1}) \to X \times ([ 0,1];\{0\},\{1\})$$
of an $n$-connected $(2n+1)$-dimensional normal map $M^{2n+1}\to X$
by surgery on a finite set of $\bbz[\pi_1(X)]$-module generators 
of $K_n(M)$). The union of the cobordisms
$$\eqalign{&((j~j'):C\oplus C'\to D,(\delta \psi,\psi\oplus - \psi'))~,\ccr
&((j'~j):C' \oplus C\to D,(- \delta \psi,\psi' \oplus - \psi))}$$
is a cobordism
$$((\widetilde j~ \widetilde j'):C\oplus C\to D',(\delta \psi',
\psi\oplus - \psi))~$$
with a $\Lambda$-module isomorphism
$$\eqalign{
&\left[\matrix{1 & 0 & -1 \ccr \psi_0 & (-1)^n \widetilde \psi_0 & 0}
\right]~:\ccr
&D'_{n+1} ~=~ {\rm coker}
(\pmatrix{ \widetilde \psi^*_0 \ccr \psi^*_0 \ccr \widetilde \psi^*_0}~
:~C^{n+1} \to C^n \oplus C_n \oplus C^n)~ \isa ~ C^n \oplus C_{n+1} ~.}$$
This verifies that cobordism is reflexive, and also that 
$(C,\psi) \oplus (C,- \psi)$ is null-cobordant.\hfil\break
\indent Suppose given a homotopy equivalence of $(2n+1)$-complexes
$$f~:~(C,\psi) \to (C',\psi')~,$$
with $\chi_0 : C'^{n+1} \to C'_{n+1}$ as in 6.6.
By reflexivity there exists a cobordism
$((j''~j'):C' \oplus C' \to D, (\delta \psi', \psi' \oplus -\psi'))$
from $(C',\psi')$ to itself. Define a cobordism
$((j~j'):C \oplus C' \to D, (\delta \psi, \psi\oplus -\psi'))$
from $(C,\psi)$ to $(C',\psi')$ by
$$\eqalign{
&j~=~j''f ~:~ C_{n+1}~ \raise4pt\hbox{$f \atop \ra{1.5} $}~ 
C'_{n+1}~ \raise4pt\hbox{$j'' \atop \ra{1.5} $}~D_{n+1}~,\ccr
&\delta\psi_0 ~=~ \delta \psi'_0 + j''\chi_0j''^* ~:~ 
D^{n+1} \to D_{n+1}~.}\eqno{\lower12pt\hbox{$\square$}}$$

\noindent {\sc Definition 7.9}
(i) A {\it weak map} of $(2n+1)$-complexes over $\Lambda$
$$f~:~(C,\psi) \to (C',\psi')$$
is a chain map $f:C\to C'$ such that there exist $\Lambda$-module morphisms
$$\chi_0~:~C^{\prime n+1} \to C^{\prime}_{n+1}~~,~~
\chi_1~:~C^{\prime n} \to C^{\prime n}$$
with
$$f\psi_0f^* - \psi^{\prime}_0~=~(\chi_0 + (-1)^{n+1} \chi^*_0)
d^{\prime *}~:~C^{\prime n} \to C^{\prime}_{n+1}~.$$
(ii) A {\it weak equivalence}
of $(2n+1)$-complexes is a weak map with $f:C\to C'$ a chain equivalence.
\hfil \break
(iii) A {\it weak isomorphism}
of $(2n+1)$-complexes is a weak map with $f:C\to C'$ an isomorphism of
 chain complexes.
\hfill$\square$

\noindent {\sc Proposition 7.10}
{\it Weakly equivalent $(2n+1)$-complexes are cobordant.}\hfil\break
{\sc Proof}\t :
The proof in 7.8 that homotopy equivalent $(2n+1)$-complexes
are cobordant works just as well for weakly equivalent ones.
\hfill$\square$

Given a $(2n+1)$-complex $(C,\psi)$ let
$$(\pmatrix{\psi_0 \ccr d^* },-\psi_1)~:~(C^n,0) \to 
 H_{(-1)^n}(C_{n+1})~$$
be the inclusion of a lagrangian in a hyperbolic split $(-1)^n$-quadratic
form given by 6.2. The result of 7.10 is that the cobordism class of
$(C,\psi)$ is independent of the hessian $(-1)^{n+1}$-quadratic form
$(C^n,-\psi_1)$.

\noindent {\sectionfont \S8. The odd-dimensional $L$-groups}

The odd-dimensional surgery obstruction groups $L_{2n+1}(\Lambda)$ of a
ring with involution $\Lambda$ will now be defined to be the cobordism
groups of $(2n+1)$-complexes over $\Lambda$.

\noindent {\sc Definition 8.1}
Let $L_{2n+1}(\Lambda)$ be the abelian group of cobordism classes of 
$(2n+1)$-complexes over $\Lambda$, with addition and inverses by
$$\eqalign{
&(C,\psi) + (C',\psi')~=~(C\oplus C',\psi\oplus \psi')~,\ccr
&-(C,\psi)~=~(C,- \psi) \in L_{2n+1}(\Lambda)~.}
\eqno{\lower10pt\hbox{$\square$}}$$

The groups $L_{2n+1}(\Lambda)$ only depend on the residue 
$n(\bmod\,2)$, so that only two $L$-groups have actually been defined, 
$L_1(\Lambda)$ and $L_3(\Lambda)$. Note that 8.1 uses 7.8 to justify 
$(C,\psi) \oplus (C,-\psi) =0 \in L_{2n+1}(\Lambda)$.

\noindent {\sc Example 8.2}
The odd-dimensional $L$-groups of $\Lambda=\bbz$ are trivial
$$L_{2n+1} (\bbz)~=~0~.\eqno{\square}$$

\indent	8.2 was implicit in the work of Kervaire and Milnor [7]
on the surgery classification of even-dimensional exotic spheres.

\noindent {\sc Example 8.3}
\rm	The surgery obstruction of an $n$-connected $(2n+1)$-dimen\-sional 
normal map $(f,b):M^{2n+1}\to X$ is the cobordism class
$$\sigma_*(f,b) ~=~ (C,\psi) \in L_{2n+1}(\bbz[\pi_1(X)])$$
of the $(2n+1)$-complex $(C,\psi)$ associated in 6.3 to any choice 
of presentation 
$$(W;M,M') \to X \times ([ 0,1];\{0\},\{1\})~.$$
The surgery obstruction vanishes $\sigma_*(f,b)=0$ if 
(and for $n \ge 2$ only if) $(f,b)$ is normal bordant to a homotopy
equivalence.
\hfill$\square$

\noindent {\sc Definition 8.4}
A {\it surgery}
$(j:C\to D,(\delta \psi,\psi))$ on a $(2n+1)$-complex 
$(C,\psi)$ is a $\Lambda$-module chain map $j:C\to D$ with 
$D_r=0$ for $r \ne n+1$ and $D_{n+1}$ a f.\ g. free $\Lambda$-module, 
together with a $\Lambda$-module morphism
$$\delta \psi_0~ :~ D^{n+1}~=~(D_{n+1})^* \to D_{n+1}~,$$
such that the $\Lambda$-module morphism 
$$\pmatrix{d & \psi^*_0j^*}~:~ C_{n+1} \oplus D^{n+1} \to C_n$$
is onto. The {\it effect}
of the surgery is the $(2n+1)$-complex $(C',\psi')$ defined by
$$\eqalign{&d'~=~ \pmatrix{d & \psi_0^* j^* \ccr (-1)^{n+1} j & 
\delta\psi_0 +(-1)^{n+1}\delta \psi^*_0 } \ccr
&\hskip50pt 
:~C'_{n+1}~=~C_{n+1} \oplus D^{n+1} \to C'_n~=~ C_n \oplus D_{n+1}~,\ccr
&\psi'_0~=~ \pmatrix{\psi_0 & 0 \ccr 0 & 1}~:~
C'^n~=~C^n \oplus D^{n+1} \to C'_{n+1}~=~ C_{n+1} \oplus D^{n+1}~,\ccr
&\psi'_1~=~\pmatrix{\psi_1 & - \psi^*_0 j^* \ccr 0 & - \delta \psi_0 }~:~
C'^n = C^n \oplus D^{n+1} \to C'_n~=~ C_n \oplus D_{n+1}~.}$$
The {\it trace} of the surgery is the cobordism of $(2n+1)$-complexes
$((j'~j''):C \oplus C' \to D',(0, \psi\oplus -\psi'))$, with
$$\eqalign{j''~=~(j' ~k) ~:~&C'_{n+1}~=~ C_{n+1} \oplus D^{n+1} \ccr
& \to D'_{n+1} ~=~{\rm ker}
(\pmatrix{d & \psi_0 j^*}:C_{n+1} \oplus D^{n+1} \to C_n)}$$
a splitting of the split injection 
$(d~\psi^*_0 j^*): C_{n+1} \oplus D^{n+1} \to C_n$.\hfill$\square$

\noindent {\sc Example 8.5}
Let
$$((e,a);(f,b),(f',b'))~:~
(V^{2n+2};M^{2n+1},M^{\prime 2n+1}) \to X \times ([ 0,1];\{0\},\{1\})$$
be the trace of a sequence of $k$ surgeries on an $n$-connected 
$(2n+1)$-dimensional normal map $(f,b):M \to X$ killing elements 
$x_1,x_2,\dots,x_k \in K_n(M)$, 
with $e$ $n$-connected and $f'$ $n$-connected. 
$V$ has a handle decomposition on $M$ of the type
$$V~=~M \times I \cup \bigcup_k (n+1)\hbox{\rm -handles}~D^{n+1} \times D^{n+1}~,$$
and also a handle decomposition on $M'$ of the same type
$$V~=~M' \times I \cup \bigcup_k (n+1)\hbox{\rm -handles}~D^{n+1} \times D^{n+1}~.$$
A presentation of $(f,b)$
$$((g,c);(\widehat f,\widehat b),(f,b))~:~
(W^{2n+2};\widehat M^{2n+1},M^{2n+1}) \to X \times ([ 0,1];\{0\},\{1\})$$
with $(2n+1)$-complex $(C,\psi)$ determines a presentation of $(f',b')$
$$\eqalign{
&((g',c');(\widehat f,\widehat b),(f',b'))~=~
((g,c);(\widehat f,\widehat b),(f,b)) \cup ((e,a);(f,b),(f',b'))~:\ccr
&\hskip10pt
(W';\widehat M,M')~=~(W;\widehat M,M) \cup (V;M,M') \to X \times ([ 0,1];\{0\},\{1\})}$$
\noindent such that the $(2n+1)$-complex $(C',\psi')$ is the effect of
a surgery $(j:C\to D,(\delta \psi,\psi))$ on $(C,\psi)$ with
$$\eqalign{&D_{n+1}\,=\,K_{n+1}(V,M')\,=\,\bbz[\pi_1(X)]^k \,, \ccr
&C'_{n+1}\,=\,K_{n+1}(W',\widehat M)\,=\,K_{n+1}(W,\widehat M) \oplus K_{n+1}(V,M)
\,=\,C_{n+1} \oplus D^{n+1}\,, \ccr
&C'_n\,=\,K_{n+1}(W',\partial W')\,=\,K_{n+1}(W,\partial W) \oplus K_{n+1}(V,M')
\,=\,C_n \oplus D_{n+1}\,.}$$
Also, the geometric trace determines the algebraic trace, with
$$D'_{n+1}~=~K_{n+1}(V)~.\eqno{\square}$$

It can be shown that $(2n+1)$-complexes $(C,\psi)$, $(C',\psi')$
are cobordant if and only if $(C',\psi')$ is homotopy equivalent to
the effect of a surgery on $(C,\psi)$. This result will only be
needed for $(C',\psi')=(0,0)$, so it will only be proved in this 
special case:

\noindent {\sc Proposition 8.6}
{\it A $(2n+1)$-complex $(C,\psi)$ represents 0 in 
$L_{2n+1}(\Lambda)$ if and only if there exists a surgery 
$(j:C \to D,(\delta \psi,\psi))$ with contractible effect.}\hfil\break
\noindent {\sc Proof}\t :
The effect of a surgery is contractible if and only if it is a 
null-cobordism.
\hfill$\square$

Given an $n$-connected $(2n+1)$-dimensional normal map 
$(f,b):M^{2n+1} \ra{1.4} \allowbreak X$ 
it is possible to kill every element $x \in K_n(M)$
by an embedding $S^n \times D^{n+1} \hookrightarrow M$
to obtain a bordant normal map 
$$(f',b') ~:~M'^{2n+1}~=~{\rm cl.}(M\backslash S^n \times D^{n+1}) \cup D^{n+1} \times S^n \to X~.$$
There are many ways of carrying out the surgery, which are quantified by
the surgeries on the kernel $(2n+1)$-complex $(C,\psi)$. In general, 
$K_n(M')$ need not be smaller than $K_n(M)$.

\noindent {\sc Example 8.7}
The kernel $(2n+1)$-complex $(C,\psi)$ over $\bbz$ of the identity normal map
$$(f,b)~=~{\rm id.}~:~M^{2n+1}~=~S^{2n+1} \to S^{2n+1}$$
is $(0,0)$. For any element 
$$\mu \in \pi_{n+1}(SO,SO(n+1))~=~Q_{(-1)^{n+1}}(\bbz)$$
let $\omega=\partial \mu\in\pi_n(SO(n+1))$, and
define a null-homotopic embedding of $S^n$ in $M$
$$e_{\omega}~:~S^n \times D^{n+1} \hookrightarrow M ~;~
(x,y) ~ \mapto (x,\omega (x)(y))/\,\Vert\,(x,\omega (x)(y))\,\Vert\,~ .$$
Use $\mu$ to kill $0 \in K_n(M)$ by surgery on $(f,b)$,
with effect a normal bordant $n$-connected $(2n+1)$-dimensional normal map
$$(f_{\mu},b_{\mu})~:~M_{\mu}^{2n+1}~=~
{\rm cl.}(M \backslash e_{\omega}(S^n \times D^{n+1})) \cup D^{n+1} \times S^n
 \to S^{2n+1}$$
exactly as in 2.18, with the kernel complex $(C',\psi')$ given by
$$d'~=~(1+T_{(-1)^{n+1}})(\mu)~:~C'_{n+1}~=~\bbz \to C'_n~=~\bbz~.$$
In particular, for $\mu =0,1$ this gives the $(2n+1)$-dimensional manifolds
$$\eqalign{&M'~=~M_0~=~ S^n \times S^{n+1} ~, \ccr 
&M''~=~M_1~=~S(\tau_{S^{n+1}}),~\hbox{\rm the tangent $S^n$-bundle of}~S^{n+1}\ccr
&\hphantom{M''~=~M_1~}=~O(n+2)/O(n)\ccr
&\hphantom{M''~=~M_1~}=~V_{n+2,2},\, 
\hbox{\rm the Stiefel manifold of orthonormal 2-frames in}\,\bbr^{n+2}\ccr 
&\hphantom{M''~=~M_1}(\hskip2.5pt =~SO(3)~=~\bbr \bbp^3~{\rm for}~n=1~)~,}$$
corresponding to the algebraic surgeries on $(0,0)$ 
$$(0:0 \to D,(\delta \psi',0))~,~(0:0 \to D,(\delta \psi'',0))$$
with 
$$D_{n+1}~=~\bbz~,~\delta \psi'_0 ~=~0~,~\delta \psi''_0 ~=~1~.\eqno{\square}$$ 

\noindent {\sectionfont \S9. Formations}

As before, let $\Lambda$ be a ring with involution, and let $\epsilon= \pm 1$.

\noindent {\sc Definition 9.1}
An {\it $\epsilon$-quadratic formation over $\Lambda$}
$(Q,\phi;F,G)$ is a non-singular $\epsilon$-quadratic form 
$(Q,\phi)$ together with an ordered pair of lagrangians $F$,$G$.
\hfill$\square$

Formations with $\epsilon= (-1)^n$ are essentially the
$(2n+1)$-complexes of \S6 expressed in the language of forms and
lagrangians of \S4. In the general theory it is possible to consider
formations $(Q,\phi;F,G)$ with $Q,F,G$ f.\ g. projective, but in view of
the more immediate topological applications only the f.\ g. 
free case is considered here. Strictly speaking, 9.1 defines a ``nonsingular
formation''. In the general theory a formation $(Q,\phi ;F,G)$ is a
nonsingular form $(Q,\phi)$ together with a lagrangian $F$ and a
sublagrangian $G$. The automorphisms of hyperbolic forms in the
original treatment due to Wall [28] of odd-dimensional surgery theory
were replaced by formations by Novikov [16] and Ranicki [18].

In dealing with formations assume that the ground ring $\Lambda$ is
such that the rank of f.\ g. free $\Lambda$-modules is well-defined
(e.g. $\Lambda=\bbz[\pi])$. The rank of a f.\ g. free $\Lambda$-module
$K$ is such that
$${\rm rank}_\Lambda(K)~=~k \in\bbz^+ $$
if and only if $K$ is isomorphic to $\Lambda^k $.
Also, since $\Lambda^k  \cong (\Lambda^k )^*$
$${\rm rank}_\Lambda (K)~=~{\rm rank}_\Lambda (K^*) \in\bbz^+ ~.$$

\noindent {\sc Definition 9.2}
An {\it isomorphism} of $\epsilon$-quadratic formations over $\Lambda$
$$f~ :~ (Q,\phi ;F,G) ~ \isa ~ (Q',\phi';F',G')$$
is an isomorphism of forms $f:(Q,\phi)\cong (Q',\phi')$ such that
$$f(F)~=~ F'~~ ,~~ f(G)~=~ G'~ .\eqno{\square}$$

\noindent {\sc Proposition 9.3}
{\it {\rm (i)} Every $\epsilon$-quadratic formation $(Q,\phi;F,G)$ is isomorphic 
to one of the type $(H_{\epsilon}(F);F,G)$. \hfil\break
{\rm (ii)} Every $\epsilon$-quadratic formation $(Q,\phi;F,G)$ is isomorphic to one 
of the type $(H_{\epsilon}(F);\-F,\alpha(F))$ for some automorphism 
$\alpha :H_{\epsilon}(F)\cong H_{\epsilon}(F)$.}\hfil\break
{\sc Proof}\t :
(i) By Theorem 5.7 the inclusion of the lagrangian $F\to Q$
extends to an isomorphism of forms $f:H_{\epsilon}(F) \cong (Q,\phi)$, 
defining an isomorphism of formations
$$f~:~(H_{\epsilon}(F);F,f^{-1}(G)) ~ \isa ~ (Q,\phi;F,G)~.$$
(ii) As in (i) extend the inclusions of the lagrangians to isomorphisms of 
forms
$$f~:~H_{\epsilon}(F)~\isa~(Q,\phi)~,~g~:~H_{\epsilon}(G) ~ 
\isa~(Q,\phi)~.$$
Then
$${\rm rank}_\Lambda(F)~=~{\rm rank}_\Lambda(Q)/2~=~{\rm rank}_\Lambda(G) \in\bbz^+ ~,$$
so that $F$ is isomorphic to $G$.
 Choosing a $\Lambda$-module isomorphism $\beta :G\cong F$
there is defined an automorphism of $H_{\epsilon}(F)$
$$\alpha~:~H_{\epsilon}(F)~\raise4pt\hbox{$\displaystyle{f} \atop \ra{1.5}$}~
(Q,\phi)~ \raise4pt\hbox{$\displaystyle{g^{-1}} \atop \ra3$}~H_{\epsilon}(G)~
\raise8pt\hbox{$\pmatrix{\beta & 0 \ccr0 & \beta^{*-1}} \atop \ra{4}$}~
H_{\epsilon}(F) $$
such that there is defined an isomorphism of formations
$$f~:~(H_{\epsilon }(F);F,\alpha (F)) ~ \isa ~ (Q,\phi ;F,G)~.\eqno{\square}$$

\noindent {\sc Proposition 9.4}
{\it The weak isomorphism classes of $(2n+1)$-complexes $(C,\psi)$ over $\Lambda$
are in natural one-one correspondence with the isomorphism classes of
$(-1)^n$-quadratic formations $(Q,\phi;F,G)$ over $\Lambda$, with
$$H_n(C)~=~Q/(F+G)~~,~~H_{n+1}(C)~=~F\cap G~.$$
Moreover, if the complex $(C,\psi)$ corresponds to the formation 
$(Q,\phi;F,G)$ then $(C, -\psi)$ corresponds to $(Q,-\phi;F,G)$.}
\hfil\break
{\sc Proof}\t : Given a $(2n+1)$-complex $(C,\psi)$ define a 
$(-1)^n$-quadratic formation
$$(Q,\phi;F,G)~=~(H_{(-1)^n}(C_{n+1});C_{n+1},{\rm im}
(\pmatrix{\psi_0 \ccr d^*}:C^n\to C_{n+1} \oplus C^{n+1}))~.$$
The formation associated in this way to the $(2n+1)$-complex $(C,-\psi)$ 
is isomorphic to $(Q,-\phi;F,G)$, by the isomorphism
$$\eqalign{\pmatrix{-1 & 0 \ccr 0 & 1}~:~&(Q,-\phi;F,G) ~ \isa \ccr
&(H_{(-1)^n}(C_{n+1});C_{n+1},{\rm im}(\pmatrix{-\psi_0 \ccr d^*}:
C^n\to C_{n+1} \oplus C^{n+1}))~.}$$
Conversely, suppose given an $(-1)^n$-quadratic formation $(Q,\phi;F,G)$.
 By 9.3 (i) this can be replaced by an isomorphic formation with 
 $(Q,\phi)=H_{(-1)^n}(F)$.
Let $\gamma \in {\rm Hom}_\Lambda(G,F)$, $\mu \in {\rm Hom}_\Lambda(G,F^*)$ be the 
components of the inclusion
$$i~=~\pmatrix{\gamma \ccr \mu}~:~G \to Q~=~F\oplus F^*~.$$
Choose any $\theta \in {\rm Hom}_\Lambda(G,G^*)$ such that
$$\gamma^* \mu~=~ \theta + (-1)^{n+1} \theta^* \in {\rm Hom}_\Lambda(G,G^*)~.$$
Define a $(2n+1)$-complex $(C,\psi)$ by
$$\eqalign{&d~=~\mu^*~:~C_{n+1}~=~F \to C_n~=~G^*~,\ccr
&\psi_0 ~=~\gamma~:~C^n~=~G \to C_{n+1}~=~F~,\ccr
&\psi_1~=~(-1)^n \theta~:~C^n~=~G \to C_n~=~G^*~.}$$
The exact sequence
$$0 \to G~ \raise4pt\hbox{$\displaystyle{i} \atop \ra{1.5} $}~ Q ~
\raise5pt \hbox{$\displaystyle{i^*(\phi + (-1)^n \phi^*)} \atop \ra{6}$}~ 
G^* \to 0$$
is the algebraic mapping cone
$$0 \to G~\raise8pt\hbox{$\pmatrix{\gamma \ccr \mu} \atop \ra{3}$}~
F\oplus F^*~ \raise4pt\hbox{$\pmatrix{\mu^* &(-1)^n\gamma^*} \atop \ra{6} $} ~
G^* \to 0$$
of the chain equivalence $(1+T)\psi_0:C^{2n+1-*}\to C$.
\hfill$\square$

\noindent {\sc Example 9.5}
An $n$-connected $(2n+1)$-dimensional normal map 
$M^{2n+1}\to \allowbreak X$ together with a choice of presentation 
$(W;M,M')\to X \times ([ 0,1];\{0\},\{1\})$
determine by 9.3 a $(2n+1)$-complex $(C,\psi)$, and hence by 9.4 a 
$(-1)^n$-quadratic formation $(Q,\phi;F,G)$ over $\bbz[\pi_1(X)]$
such that
$$\eqalign{
&Q/(F+G)~=~H_n(C)~=~K_n(M)~,\ccr
&F\cap G~=~H_{n+1}(C)~=~K_{n+1}(M)~.}\eqno{\lower10pt\hbox{$\square$}}$$

The following equivalence relation on formations 
corresponds to the weak equivalence (7.9) of $(2n+1)$-complexes.
 
\noindent {\sc Definition 9.6}
(i) An $\epsilon$-quadratic formation $(Q,\phi;F,G)$ is {\it trivial}
if it is isomorphic to $(H_{\epsilon}(L);L,L^*)$ for some f.\ g.
 free $\Lambda$-module $L$. \hfil\break
(ii) A {\it stable isomorphism} of $\epsilon$-quadratic formations
$$[f]~:~(Q,\phi;F,G) ~ \isa ~ (Q',\phi';F',G')$$
is an isomorphism of $\epsilon$-quadratic formations of the type
$$f~:~(Q,\phi;F,G)\oplus ({\rm trivial}) ~ \isa ~ (Q',\phi';F',G') 
\oplus ({\rm trivial}')~.\eqno{\square}$$

\noindent {\sc Example 9.7}
The $(-1)^n$-quadratic formations associated in 9.5 to all 
the presentations of an $n$-connected $(2n+1)$-dimensional normal map 
$M^{2n+1}\to X$ define a stable isomorphism class.
\hfill$\square$

\noindent {\sc Proposition 9.8}
{\it The weak equivalence classes of $(2n+1)$-complexes over $\Lambda$ are 
in natural one-one correspondence with the stable isomorphism classes of 
$(-1)^n$-quadratic formations over $\Lambda$.} \hfil\break
{\sc Proof}\t : The $(2n+1)$-complex $(C,\psi)$ associated (up to weak 
equivalence) to a $(-1)^n$-quadratic formation $(Q,\phi;F,G)$ in 9.4 
is contractible if and only if the formation is trivial.
\hfill$\square$

The following formations correspond to the null-cobordant complexes.

\noindent {\sc Definition 9.9}
The {\it boundary} of a $(- \epsilon)$-quadratic form $(K,\lambda,\mu)$ is the 
$\epsilon$-quadratic formation
$$\partial (K,\lambda,\mu)~=~(H_{\epsilon }(K);K,\Gamma_{(K,\lambda)})$$
with $\Gamma_{(K,\lambda)}$ the {\it graph} lagrangian
$$\Gamma_{(K,\lambda)}~=~\{(x,\lambda(x))
\in K\oplus K^* \,\vert\, x\in K\}~.\eqno{\square}$$
\indent
Note that the form $(K,\lambda,\mu)$ may be singular, that is the
$\Lambda$-module morphism $\lambda:K\to K^*$ need not be an
isomorphism. The graphs $\Gamma_{(K,\lambda)}$ of $(-\epsilon)$-quadratic 
forms $(K,\lambda,\mu)$ are precisely the lagrangians of
$H_{\epsilon}(K)$ which are direct complements of $K^*$.

\noindent {\sc Proposition 9.10}
{\it A $(-1)^n$-quadratic formation $(Q,\phi;F,G)$ is stably isomorphic 
to a boundary $\partial (K,\lambda,\mu)$ if and only if the corresponding
$(2n+1)$-complex $(C,\psi)$ is null-cobordant.}\hfil\break
{\sc Proof}\t :
Given a $(-1)^{n+1}$-quadratic form $(K,\lambda,\mu)$ choose a split
form $\theta:K \to K^*$ (4.2) and let $(C,\psi)$ be the
$(2n+1)$-complex associated by 9.4 to the boundary formation
$\partial(K,\lambda,\mu)$, so that
$$\eqalign{
&d~=~\lambda~=~\theta + (-1)^{n+1} \theta^*~:~C_{n+1}~=~K \to C_n~=~ K^*~,\ccr
&\psi_0~=~ 1~:~C^n~=~K \to C_{n+1}~=~ K~,\ccr
&\psi_1 ~=~- \theta ~:~ C^n~=~ K \to C_n~=~ K^*~.}$$
Then $(C,\psi)$ is null-cobordant, with a null-cobordism
$(j:C\to D,(\delta \psi,\psi))$ defined by
$$\eqalign{&j~=~1~:~C_{n+1}~=~K \to D_{n+1}~=~K~,\ccr
&\delta \psi_0~=~ 0~ :~ D^{n+1}~=~K^* \to D_{n+1}~=~K~.}$$
\indent	Conversely, suppose given a $(2n+1)$-complex $(C,\psi)$ with a 
null-cobordism $(j:C\to D,(\delta \psi,\psi))$ as in 8.1.
The $(2n+1)$-complex $(E,\theta)$ defined by
$$\eqalign{&d\,=\,\pmatrix{\psi_1 + (-1)^{n+1}\psi^*_1 & d & \psi^*_0j^* \ccr
(-1)^{n+1}d^* & 0 & -j^* \ccr
(-1)^{n+1}j\psi^*_0 & (-1)^nj & \delta \psi_0 + (-1)^{n+1}
\delta \psi^*_0}\,:\ccr
&\hskip25pt E_{n+1}\,=\,C^n \oplus C_{n+1} \oplus D^{n+1} \, \to \, 
E_n\,=\,C_n \oplus C^{n+1} \oplus D_{n+1}\,,\ccr
&\theta_0\,=\,1\,:\,E^n\,=\,C^n \oplus C_{n+1} \oplus D^{n+1}
\, \to \, E_{n+1}\,=\,C^n \oplus C_{n+1} \oplus D^{n+1}\,,\ccr
&\theta_1\,=\,\pmatrix{-\psi_1 & -d & -\psi^*_0 j^* \ccr
 0 & 0 & j^* \ccr 0 & 0 & -\delta \psi_0}\,:\ccr
&\hskip25pt E^n\,=\,C^n \oplus C_{n+1} \oplus D^{n+1} \, \to \, 
E_n\,=\,C_n \oplus C^{n+1} \oplus D_{n+1}}$$
corresponds to the boundary $(-1)^n$-quadratic formation 
$\partial (E^n,\lambda_1,\mu_1)$ of the $(-1)^{n+1}$-quadratic form 
$(E^n,\lambda_1,\mu_1)$ determined by the split form $\theta_1$, 
and there is defined a homotopy equivalence $f:(E,\theta) \to (C,\psi)$
with
$$\eqalign{&f_n~=~(1~~\psi^*_0~~0)~:~E_n~=~C_n \oplus C^{n+1} \oplus 
D_{n+1} \to C_n~,\ccr
&f_{n+1}~=~(0~~ 1~~ 0)~:~E_{n+1}~=~C^n \oplus C_{n+1}
 \oplus D^{n+1} \to C_{n+1}~.}\eqno{\lower10pt\hbox{$\square$}}$$

\noindent {\sc Proposition 9.11}
{\it The cobordism group $L_{2n+1}(\Lambda)$ of $(2n+1)$-complexes is naturally 
isomorphic to the abelian group of equivalence classes of $(-1)^n$-quadratic
formations over $\Lambda$, subject to the equivalence relation
$$\eqalign{&(Q,\phi;F,G) \sim (Q',\phi';F',G')~ 
\hbox{\it if there exists a stable isomorphism}\ccr
&[f]~:~(Q,\phi;F,G) \oplus (Q',- \phi' ;F',G') ~ \isa ~ \partial 
(K,\lambda,\mu)\ccr
&\hbox{\it for some $(-1)^{n+1}$-quadratic form $(K,\lambda,\mu)$ over $\Lambda$}~,}$$
with addition and inverses by}
$$\eqalign{
(Q,\phi ;F,G) + (Q',\phi';F',G')~&=~(Q\oplus Q',\phi \oplus \phi';
F\oplus F',G\oplus G')~,\ccr
-(Q,\phi;F,G) ~&=~(Q,- \phi ;F,G) \in L_{2n+1}(\Lambda)~.}$$
{\sc Proof}\t :
This is just the translation of the definition (8.1) of $L_{2n+1}(\Lambda)$
into the language of $(-1)^n$-quadratic formations, using 9.4, 9.8 and 9.10.
\hfill$\square$

Use 9.11 as an identification of $L_{2n+1}(\Lambda)$
with the group of equivalence classes of $(-1)^n$-quadratic formations 
over $\Lambda$.

\noindent {\sc Corollary 9.12}
{\it A $(-1)^n$-quadratic formation $(Q,\phi;F,G)$ over $\Lambda$ is such that
$(Q,\phi;F,G)=0 \in L_{2n+1}(\Lambda)$ if and only if it is stably isomorphic
to the boundary $\partial (K,\lambda,\mu)$ of a $(-1)^{n+1}$-quadratic form
$(K,\lambda,\mu)$ on a f.\ g. free $\Lambda$-module $K$.} \hfil\break
\noindent {\sc Proof}\t : Immediate from 9.10. \hfill$\square$

Next, it is necessary to establish the relation
$$(Q,\phi ;F,G)\oplus (Q,\phi;G,H) ~=~(Q,\phi;F,H) \in L_{2n+1}(\Lambda)~.$$
This is the key step in the identification in \S10 below of $L_{2n+1}(\Lambda)$
 with a stable unitary group.

\noindent {\sc Lemma 9.13}
{\it {\rm (i)} An $\epsilon$-quadratic formation $(Q,\phi;F,G)$ is trivial if 
and only if the lagrangians $F$ and $G$ are direct complements in $Q$.
\hfil\break
{\rm (ii)} An $\epsilon$-quadratic formation $(Q,\phi;F,G)$ is isomorphic to a 
boundary if and only if $(Q,\phi)$ has a lagrangian $H$ which is a direct 
complement of both the lagrangians $F$,\t $G$.} \hfil\break
{\sc Proof}\t :
(i) If $F$ and $G$ are direct complements in $Q$ express any representative 
$\phi \in {\rm Hom}_\Lambda (Q,Q^*)$ of $\phi \in Q_{\epsilon}(Q)$ as
$$\phi ~=~\pmatrix{\lambda -\epsilon \lambda^* & \gamma \ccr
\delta & \mu -\epsilon \mu^*}~:~Q~=~F \oplus G \to Q^*~=~F^*
 \oplus G^*~.$$
Then $\gamma + \epsilon \delta^* \in {\rm Hom}_\Lambda (G,F^*)$ is an 
$\Lambda$-module isomorphism, and there is defined an isomorphism of 
$\epsilon$-quadratic formations
$$\pmatrix{1 & 0 \ccr 0 & (\gamma + \epsilon \delta^*)^{-1}}~
:~(H_{\epsilon }(F);F,F^*) ~ \isa ~ (Q,\phi;F,G)$$
so that $(Q,\phi;F,G)$ is trivial.
 The converse is obvious. \hfil\break
(ii) For the boundary $\partial (K,\lambda,\mu)$ of a 
$(- \epsilon)$-quadratic form $(K,\lambda,\mu)$ the lagrangian $K^*$ 
of $H_{\epsilon }(K)$ is a direct complement of both the lagrangians 
$K$,\t $\Gamma_{(K,\lambda)}$. Conversely, suppose that $(Q,\phi;F,G)$ is 
such that there exists a lagrangian $H$ in $(Q,\phi)$ which is a direct 
complement to both $F$ and $G$.
 By the proof of (i) there exists an isomorphism of formations
$$f~:~(H_{\epsilon}(F);F,F^*) ~ \isa ~ (Q,\phi;F,H)$$
which is the identity on $F$.
 Now $f^{-1}(G)$ is a lagrangian of $H_{\epsilon }(F)$ which is a direct 
complement of $F^*$, so that it is the graph $\Gamma_{(F,\lambda)}$ of 
a $(-\epsilon)$-quadratic form $(F,\lambda,\mu)$, and $f$ defines an isomorphism 
of $\epsilon$-quadratic formations
$$f~:~\partial (F,\lambda,\mu)~=~
(H_{\epsilon}(F);F,\Gamma_{(F,\lambda)})~\isa~(Q,\phi;F,G)~.\eqno{\square}$$

\noindent {\sc Proposition 9.14}
{\it For any lagrangians $F,G,H$ in a nonsingular $(-1)^n$-quadratic form 
 $(Q,\phi)$ over $\Lambda$}
$$(Q,\phi;F,G)\oplus (Q,\phi;G,H) ~=~(Q,\phi;F,H) \in L_{2n+1}(\Lambda)~.$$
{\sc Proof}\t :
Choose lagrangians $F^*,G^*,H^*$ in $(Q,\phi)$ complementary to $F,G,H$
respectively.  The $(-1)^n$-quadratic formations $(Q_i,\phi_i;F_i,G_i)$ 
$(1\le i \le 4)$ defined by
$$\eqalign{
&(Q_1,\phi_1;F_1,G_1)~=~(Q,-\phi;G^*,G^*)~,\ccr
&(Q_2,\phi_2;F_2,G_2)~=~(Q\oplus Q,\phi \oplus - \phi;F\oplus F^*,H\oplus G^*)\ccr
&\hphantom{(Q_2,\phi_2;F_2,G_2)~=~(Q\oplus Q,\phi}
\oplus (Q\oplus Q,- \phi \oplus \phi ;\Delta_{Q},H^* \oplus G)~,\ccr
&(Q_3,\phi_3;F_3,G_3)~=~(Q\oplus Q,\phi \oplus - \phi,F\oplus F^*,
G\oplus G^*)~,\ccr
&(Q_4,\phi_4;F_4,G_4)~=~(Q\oplus Q,\phi \oplus - \phi;
G\oplus G^*,H\oplus G^*)\ccr
&\hphantom{(Q_2,\phi_2;F_2,G_2)~=~(Q\oplus Q,\phi} 
\oplus (Q\oplus Q,- \phi \oplus \phi ; \Delta_{Q},H^* \oplus G)}$$
are such that
$$\eqalign{
(Q,\phi ;F,G)~&\oplus (Q,\phi ;G,H)\oplus (Q_1,\phi_1;F_1,G_1) 
\oplus (Q_2,\phi_2;F_2,G_2)\ccr
&=~(Q,\phi ;F,H)\oplus (Q_3,\psi_3;F_3,G_3)
 \oplus (Q_4,\phi_4;F_4,G_4)~.}$$
 Each of $(Q_i,\phi_i;F_i,G_i)~ (1\le i\le 4)$
is isomorphic to a boundary, since there exists a lagrangian $H_i$
in $(Q_i,\phi_i)$ complementary to both $F_i$ and $G_i$, 
so that 9.13 (ii) applies and $(Q_i,\phi_i;F_i,G_i)$ represents $0$
in $L_{2n+1}(\Lambda)$. Explicitly, take
$$\eqalign{&H_1~=~G \subset Q_1~=~ Q~ ,\ccr
&H_2~=~\Delta_{Q\oplus Q} \subset Q_2~=~(Q\oplus Q)\oplus (Q\oplus Q)~,\ccr
&H_3~=~\Delta_{Q} \subset Q_3~=~ Q\oplus Q~,\ccr
&H_4~=~\Delta_{Q\oplus Q} \subset Q_4~=~(Q\oplus Q)\oplus (Q\oplus Q)~.}
\eqno{\lower28pt\hbox{$\square$}}$$

\noindent {\sc Remark 9.15}
It is also possible to express $L_{2n+1}(\Lambda)$ as the abelian group of 
equivalence classes of $(-1)^n$-quadratic formations over $\Lambda$ subject to
the equivalence relation generated by
{\parindent=20pt
\parskip=4pt
\item{\rm (i)} $(Q,\phi ;F,G) \sim (Q',\phi';F',G')$
if $(Q,\phi;F,G)$ is stably isomorphic to $(Q',\phi';F',G')$,
\item{\rm (ii)} $(Q,\phi ;F,G)\oplus (Q,\phi ;G,H) \sim (Q,\phi ;F,H)$,
with addition and inverses by
$$\eqalign{
(Q,\phi ;F,G) + (Q',\phi';F',G')~&=~
 (Q\oplus Q',\phi \oplus \phi';F\oplus F',G\oplus G') ~,\ccr
-(Q,\phi;F,G)~&=~ (Q,\phi;G,F) \in L_{2n+1}(\Lambda)~.}$$\par}
\vskip-4pt

\noindent This is immediate from 9.13 and the observation that for any 
$(-1)^{n+1}$-quadratic form $(K,\lambda,\mu)$ on a f.\ g.
 free $\Lambda$-module $K$ the lagrangian $K^*$ in $H_{(-1)^n}(K)$ is a complement 
 to both $K$ and the graph $\Gamma_{(K,\lambda)}$, so that
$$\eqalign{\partial (K,\lambda,\mu)~& \sim~ (H_{(-1)^n}(K);K,\Gamma_{(K,\lambda)})
 \oplus (H_{(-1)^n}(K);\Gamma_{(K,\lambda)},K^*)\ccr
 &\sim~ (H_{(-1)^n}(K);K,K^*)~\sim~ 0~ .}\eqno{\lower8pt\hbox{$\square$}}$$

\noindent {\sectionfont \S10. Automorphisms}

The $(2n+1)$-dimensional $L$-group $L_{2n+1}(\Lambda)$ of a ring with
involution $\Lambda$ is identified with a quotient of the stable
automorphism group of hyperbolic $(-1)^n$-quadratic forms over
$\Lambda$, as in the original definition of Wall [28].

Given a $\Lambda$-module $K$ let Aut$_\Lambda (K)$ be the group of
automorphisms $K\to K$, with the composition as group law.

\noindent {\sc Example 10.1}
The automorphism group of the f.\ g. free $\Lambda$-module $\Lambda^k $
is the general linear group $GL_k (\Lambda)$ of invertible $k \times k$
matrices in $\Lambda$
$${\rm Aut}_\Lambda(\Lambda^k )~=~GL_k(\Lambda)$$
with the multiplication of matrices as group law (cf. Remark 1.12).
The general linear group is not abelian for $k \geq 2$, since
$$\pmatrix{1&1\ccr0&1} \pmatrix{1&0\ccr 1&1}~\ne~
\pmatrix{1&0\ccr1&1} \pmatrix{1&1\ccr 0&1}~.\eqno{\square}$$

\noindent {\sc Definition 10.2}
For any $\epsilon$-quadratic form $(K,\lambda,\mu)$ let 
${\rm Aut}_\Lambda(K,\lambda,\mu)$ be the
subgroup of ${\rm Aut}_\Lambda(K)$ consisting of the automorphisms
$f:(K,\lambda,\mu)\to (K,\lambda,\mu)$.\hfill$\square$

\noindent {\sc Definition 10.3}
The $(\epsilon,k)$-{\it unitary group} of $\Lambda$
is defined for $\epsilon=\pm 1$, $k\ge 0$ to be the automorphism group
$$U_{\epsilon ,k}(\Lambda)~=~{\rm Aut}_\Lambda (H_{\epsilon}(\Lambda^k ))$$
of the $\epsilon$-quadratic hyperbolic form $H_{\epsilon}(\Lambda^k )$.
\hfill$\square$

\noindent {\sc Proposition 10.4}
$U_{\epsilon ,k}(\Lambda)$ is the group of invertible $2k\times 2k$
matrices 
$\pmatrix{\alpha & \beta \ccr \gamma & \delta } \in GL_{2k}(\Lambda)$
such that
$$\alpha^* \delta + \epsilon \gamma^* \beta ~=~ 1 \in M_{k,k}(\Lambda)~,~
\alpha^* \gamma ~=~ \beta^* \delta ~ =~ 0 \in Q_{\epsilon}(\Lambda^k )~.$$
\noindent {\sc Proof}\t :
This is just the decoding of the condition
$$\pmatrix{\alpha^* & \gamma^* \ccr\beta^* & \delta^*}
\pmatrix{0&1\ccr0&0}
\pmatrix{\alpha & \beta \ccr\gamma & \delta}~=~
\pmatrix{0&1\ccr0&0} \in Q_{\epsilon}(\Lambda^k \oplus (\Lambda^k )^*)$$
for $\pmatrix{\alpha &\beta \ccr\gamma & \delta }$ to define an 
automorphism of the hyperbolic (split) $\epsilon$-quadratic form 
$$H_{\epsilon}(\Lambda^k )~=~(\Lambda^k  \oplus (\Lambda^k )^*,
\pmatrix{0&1\ccr0&0})~.\eqno{\square}$$

\indent	Use 10.4 to express the automorphisms of $H_{\epsilon}(\Lambda^k )$ 
as matrices.

\noindent {\sc Example 10.5}
$U_{\epsilon ,1}(\Lambda)$ is the subgroup of $GL_2(\Lambda)$ consisting of the
$2\times 2$ matrices $\pmatrix{a & b \ccr c & d }$ such that
$$d \bar a + \epsilon b \bar c~=~1 \in \Lambda~,~
c \bar a~=~d \bar b ~=~0 \in Q_{\epsilon}(\Lambda)~.\eqno{\square}$$

\noindent {\sc Definition 10.6}
The {\it elementary $(\epsilon,k)$-quadratic unitary group} of $\Lambda$
is the normal subgroup 
$$EU_{\epsilon,k}(\Lambda) \subseteq U_{\epsilon,k}(\Lambda)$$
of the full $(\epsilon,k)$-quadratic unitary group
generated by the elements of the following two types:
{\parindent=20pt
\parskip=0pt
\item{\rm (i)} $\pmatrix{\alpha & 0\ccr0 & \alpha^{*-1}}$ 
for any automorphism $\alpha \in GL_k (\Lambda)$~,
\item{\rm (ii)} $\pmatrix{1& 0 \ccr \theta -\epsilon \theta^* & 1}$ for any 
split $(-\epsilon)$-quadratic form $(\Lambda^k,\theta)$.
\hfill$\square$\par}

\noindent {\sc Lemma 10.7}
{\it For any $(-\epsilon)$-quadratic form 
$(\Lambda^k,\theta \in Q_{- \epsilon}(\Lambda^k))$}
$$\pmatrix{1 & \theta - \epsilon \theta^* \ccr 0 & 1}
\in EU_{\epsilon,k}(\Lambda)~.$$
{\sc Proof}\t :
This is immediate from the identity
$$\pmatrix{1 & \theta - \epsilon \theta^* \ccr 0 & 1}~
=~\pmatrix{0 & 1 \ccr 1 & 0}^{-1} 
\pmatrix{1 & 0 \ccr \theta - \epsilon \theta^* & 1}
\pmatrix{0 & 1 \ccr 1 & 0}.\eqno{\lower6pt\hbox{$\square$}}$$

\indent	Use the identifications
$$\Lambda^{k+1}~=~\Lambda^k  \oplus \Lambda~,~H_{\epsilon}(\Lambda^{k+1})~=~H_{\epsilon}(\Lambda^k )
 \oplus H_{\epsilon}(\Lambda)$$
to define injections of groups
$$U_{\epsilon ,k}(\Lambda) \to U_{\epsilon ,k+1}(\Lambda)~;~f ~ \mapto f\oplus 1 ,$$
such that $EU_{\epsilon,k}(\Lambda)$ is sent into $EU_{\epsilon,k+1}(\Lambda)$.

\noindent {\sc Definition 10.8}
(i) The {\it stable $\epsilon$-quadratic unitary group}
of $\Lambda$ is the union
$$U_{\epsilon}(\Lambda)~=~\bigcup_{k=1}^{\infty}U_{\epsilon ,k}(\Lambda)~.$$
(ii) The {\it elementary stable $\epsilon$-quadratic unitary group}
of $\Lambda$ is the union
$$EU_{\epsilon}(\Lambda)~=~\bigcup^{\infty}_{k=1}EU_{\epsilon,k}(\Lambda)~,$$
a normal subgroup of $U_{\epsilon}(\Lambda)$. \hfil\break
(iii) The {\it ${\epsilon}$-quadratic $M$-group}
of $\Lambda$ is the quotient
$$M_{\epsilon}(\Lambda)~=~ U_{\epsilon}(\Lambda)/
\{EU_{\epsilon}(\Lambda),\sigma_{\epsilon}\}$$
with $\sigma_{\epsilon}=\pmatrix{ 0 & 1 \ccr \epsilon & 0 } 
\in U_{\epsilon,1}(\Lambda) \subseteq U_{\epsilon}(\Lambda)$.
\hfill$\square$

The automorphism group
$M_{\epsilon}(\Lambda)$ is the original definition due to Wall [28,\t Chap.\t 6]
of the odd-dimensional $L$-group $L_{2n+1}(\Lambda)$, with $\epsilon = (-1)^n$.
The original verification that $M_{\epsilon}(\Lambda)$ is abelian 
used a somewhat complicated matrix identity ([28,\t p.66]),
corresponding to the formation identity 9.14.
Formations will now be used to identify $M_{(-1)^n}(\Lambda)$ with 
the a priori abelian $L$-group $L_{2n+1}(\Lambda)$ defined in \S8. 

Given an automorphism of a hyperbolic $(-1)^n$-quadratic form 
$$\alpha~=~\pmatrix{\gamma & \widetilde \gamma \ccr \mu & \widetilde \mu}~:~
H_{(-1)^n}(\Lambda^k ) ~ \isa ~ H_{(-1)^n}(\Lambda^k)~$$
define a $(2n+1)$-complex $(C,\psi)$ by
$$\eqalign{&d~=~\mu^*~:~C_{n+1}~=~\Lambda^k  \to C_n~=~\Lambda^k ~,\ccr
&\psi_0~=~\gamma~:~C^n~=~\Lambda^k  \to C_{n+1}~=~\Lambda^k ~,}$$
corresponding to the $(-1)^n$-quadratic formation
$$\Phi_k (\alpha)~=~(H_{(-1)^n}(\Lambda^k );\Lambda^k ,
{\rm im}(\pmatrix{\gamma \ccr \mu}:\Lambda^k \to \Lambda^k \oplus (\Lambda^k)^*))~.$$

\noindent {\sc Lemma 10.9}
{\it The formations $\Phi_k(\alpha_1)$, $\Phi_k(\alpha_2)$
associated to two automorphisms
$$\alpha_i~=~\pmatrix{\gamma_i & \widetilde \gamma_i \ccr
\mu_i & \widetilde \mu_i}~:~H_{(-1)^n}(\Lambda^k) ~ \isa ~ H_{(-1)^n}(\Lambda^k) ~~ (i=1,2)$$
are isomorphic if and only if there exist $\beta_i \in GL_k (\Lambda)$,
$\theta_i \in S(\Lambda^k)$ such that}
$$\eqalign{\pmatrix{\beta_1 & 0 \ccr 0 & \beta^{*-1}_1}&
\pmatrix{1 & \theta_1+(-1)^{n+1}\theta_1^* \ccr 0 & 1 }
\pmatrix{\gamma_1 & \widetilde \gamma_1 \ccr \mu_1 & \widetilde \mu_1}\ccr
&=~\pmatrix{\gamma_2 & \widetilde \gamma_2 \ccr \mu_2 & \widetilde \mu_2}
\pmatrix{\beta_2 & 0 \ccr 0 & \beta^{*-1}_2}
\pmatrix{1 & \theta_2 + (-1)^{n+1}\theta_2^* \ccr 0 & 1}\ccr
&\hskip75pt :~H_{(-1)^n}(\Lambda^k ) ~ \isa ~ H_{(-1)^n}(\Lambda^k)~.}$$
{\sc Proof}\t : 
An automorphism $\alpha$ of the hyperbolic $(-1)^n$-quadratic form
$H_{(-1)^n}(\Lambda^k)$ preserves the lagrangian $\Lambda^k \subset
\Lambda^k \oplus (\Lambda^k)^*$ if and only if there exist 
$\beta \in GL_k(\Lambda)$, $\theta \in S(\Lambda^k)$ such that
$$\alpha~=~\pmatrix{\beta & 0 \ccr 0 & \beta^{*-1}}
\pmatrix{1 & \theta + (-1)^{n+1}\theta^* \ccr 0 & 1}~ :~
H_{(-1)^n}(\Lambda^k ) ~ \isa ~ H_{(-1)^n}(\Lambda^k)~.$$ 
\line{\hfill$\square$}

\noindent {\sc Proposition 10.10}
{\it The function
$$\Phi~:~M_{(-1)^n}(\Lambda) \to L_{2n+1}(\Lambda)~;~\alpha ~ \mapto \Phi_k(\alpha)~~
(\alpha \in U_{(-1)^n,k}(\Lambda))~$$
is an isomorphism of groups.}\hfil\break
\noindent {\sc Proof}\t : The function
$$\Phi_k ~:~U_{(-1)^n,k}(\Lambda) \to L_{2n+1}(\Lambda)~;~
\alpha \mapto \Phi_k  (\alpha)$$
is a group morphism by 9.14. Each of the generators (10.6) of the
elementary subgroup $EU_{(-1)^n,k}(\Lambda)$ is sent to $0$ with

{\parindent=23pt
\item{\rm (i)} $\Phi_k\pmatrix{ \beta & 0 \ccr 0 & \beta^{*-1} }=
(H_{(-1)^n}(\Lambda^k);\Lambda^k,\Lambda^k)
=\partial (\Lambda^k,0,0) =0 \in L_{2n+1}(\Lambda),$
\item{\rm (ii)} $\Phi_k\pmatrix{ 1 & 0 \ccr \theta +(-1)^{n+1}\theta^* & 1}
=\partial (\Lambda^k ,\theta+(-1)^{n+1}\theta^*,\theta) =0 \in L_{2n+1}(\Lambda).$}

\noindent Also, abbreviating $\sigma_{(-1)^n}$ to $\sigma$
$$\eqalign{
&\Phi_1 (\sigma)~=~ (H_{(-1)^n}(\Lambda);\Lambda,\Lambda^*)~=~ 0~,\ccr
&\Phi_{k+1} (\alpha \oplus 1) ~=~\Phi_k  (\alpha) \oplus 
(H_{(-1)^n}(\Lambda);\Lambda,\Lambda) ~=~\Phi_k  (\alpha) \in L_{2n+1}(\Lambda)~.}$$
Thus the morphisms $\Phi_k$ $(k\ge 0)$ fit together to define a group morphism
$$\Phi ~:~M_{(-1)^n}(\Lambda) \to L_{2n+1}(\Lambda) ~;~
 \alpha \mapto \Phi_k  (\alpha)~{\rm if}~ \alpha \in U_{(-1)^n,k}(\Lambda)$$
such that
$$\Phi (\alpha_1 \alpha_2) ~=~\Phi (\alpha_1 \oplus \alpha_2) ~=~
\Phi (\alpha_1) \oplus \Phi (\alpha_2) \in L_{2n+1}(\Lambda)~ .$$
$\Phi$ is onto by 9.3 (ii). It remains to prove that $\Phi$ is one-one.\hfil\break
\indent
For any $\alpha_i \in U_{(-1)^n,k_i}(\Lambda)~(i=1,2)$ 
$$\alpha_1 \oplus \alpha_2 ~=~\alpha_2 \oplus \alpha_1 \in M_{(-1)^n}(\Lambda)~,$$
since
$$\eqalign{\pmatrix{\alpha_1 & 0 \ccr 0 & \alpha_2}~&=~
\pmatrix{0 & 1 \ccr 1 & 0}^{-1} \pmatrix{\alpha_2 & 0 \ccr 0 & \alpha_1}
\pmatrix{0 & 1 \ccr 1 & 0}\ccr
&:~H_{(-1)^n}(\Lambda^{k_1 + k_2}) \to H_{(-1)^n}(\Lambda^{k_1+k_2})~.}$$
Now $\sigma = 1 \in M_{(-1)^n}(\Lambda)$ (by construction), 
so that for any $\alpha \in U_{(-1)^n,k}(\Lambda)$
$$\alpha \oplus \sigma~=~\sigma \oplus \alpha~=~(\sigma \oplus 1)
(1\oplus \alpha)~=~\alpha \in M_{(-1)^n}(\Lambda)~.$$
It follows that for every $m\ge 1$
$$\sigma \oplus \sigma \oplus \dots \oplus \sigma~=~
1 \in M_{(-1)^n}(\Lambda)~ (m\hbox{\rm -fold sum})~.$$
\indent If $\alpha \in U_{(-1)^n,k}(\Lambda)$ is such that
$\Phi (\alpha) = 0 \in L_{2n+1}(\Lambda)$ then by 9.12 the $(-1)^n$-quadratic 
formation $\Phi_k  (\alpha)$ is stably isomorphic to the boundary 
$\partial (\Lambda^{k'},\lambda,\mu)$ of a $(-1)^{n+1}$-quadratic form 
$(\Lambda^{k'},\lambda,\mu)$. Choosing a split form $\theta \in S(\Lambda^{k'})$ for $(\lambda,\mu)$
this can be expressed as
$$\partial (\Lambda^{k'},\lambda,\mu) ~=~\Phi_{k'}
\pmatrix{ 1 & 0 \ccr \theta + (-1)^{n+1}\theta^* & 1}~.$$
Thus for a sufficiently large $k'' \ge 0$ there exist by 10.9
$\beta_i \in GL_{k''}(\Lambda)$, $\theta_i \in S(\Lambda^{k''})$ 
$(i=1,2)$ such that
$$\eqalign{
&\pmatrix{\beta_1 & 0 \ccr 0 & \beta^{*-1}_1}
\pmatrix{1 & \theta_1+(-1)^{n+1}\theta_1^* \ccr 0 & 1 }
(\alpha \oplus \sigma \oplus \dots \oplus \sigma)\ccr
&=~\bigg(\pmatrix{1 & 0 \ccr \theta + (-1)^{n+1} \theta^* & 1}
\oplus\sigma \oplus \dots \oplus \sigma\bigg)
\pmatrix{\beta_2 & 0 \ccr 0 & \beta^{*-1}_2}\cr
&\hskip50pt
\pmatrix{1 & \theta_2 + (-1)^{n+1}\theta_2^* \ccr 0 & 1}~
:~H_{(-1)^n}(\Lambda^{k''}) \to H_{(-1)^n}(\Lambda^{k''})}$$
so that by another application of 10.7
$$\alpha~=~\pmatrix{1 & 0 \ccr \theta + (-1)^{n+1}\theta^* & 1}~
=~1 \in M_{(-1)^n}(\Lambda)~,$$
verifying that $\Phi$ is one-one.\hfill$\square$

\noindent{\sectionfont References}
\bigskip

{\parskip=0pt
\parindent=23pt
\item{[1]} W.\t Browder, {\it Surgery on simply-connected manifolds},
Springer (1972)
\item{[2]} \bysame,
{\it Differential topology of higher dimensional manifolds}, in
{\it Surveys on surgery theory, Volume 1}, Ann. of Maths. Studies 145, 41--71,
Princeton (2000)
\item{[3]} J.\t Bryant, S.\t Ferry, W.\t Mio and S.\t Weinberger,
{\it Topology of homology manifolds}, Ann. of Maths. (2) 143, 435--467 (1996)
\item{[4]} I.\t Hambleton and L.\ Taylor,
{\it A guide to the calculation of surgery obstruction groups for finite
groups}, in
{\it Surveys on surgery theory, Volume 1}, Ann. of Maths. Studies 145, 225--274,
Princeton (2000)
\item{[5]} M.\t Kervaire,
{\it A manifold which does not admit a differentiable structure},
Comm. Math. Helv. 34, 257--270 (1960)
\item{[6]} \bysame,
{\it Les noeuds de dimensions sup\'erieures},
Bull. Soc. Math. France 93, 225--271 (1965)
\item{[7]} \bysame\ and J.\t Milnor, 
{\it Groups of homotopy spheres I.}, Ann. of Maths. 77, 504--537 (1963)
\item{[8]} R.\t Kirby and L.\t Siebenmann, 
{\it Foundational essays on topological manifolds, smoothings, 
and triangulations},
Ann. of Maths. Studies 88, Princeton (1977)
\item{[9]} J.\t Levine, {\it Lectures on groups of homotopy spheres},
Algebraic and Geometric Topology, Rutgers 1983, Lecture Notes in Mathematics
1126, 62--95, Springer (1983)
\item{[10]} J.\t Milnor, {\it On manifolds homeomorphic to the 7-sphere},
Ann. of Maths. 64, 399--405 (1956)
\item{[11]} \bysame, {\it Differentiable manifolds which are homotopy spheres},
notes (1959)
\item{[12]} \bysame, {\it Differentiable structures on spheres},
Amer. J. of Math. 81,  962--972 (1959)
\item{[13]} \bysame and D.\t Husemoller,
{\it Symmetric bilinear forms},
Ergebnisse der Mathematik und ihrer Grenzgebiete 73, Springer (1973)
\item{[14]} \bysame and J.\t Stasheff,
{\it Characteristic classes}, Ann. of Maths. Studies 76, Princeton (1974)
\item{[15]} A.\t S.\t Mishchenko,
{\it Homotopy invariants of non--simply connected manifolds,
III. Higher signatures},
Izv. Akad. Nauk SSSR, ser. mat. 35, 1316--1355 (1971) 
\item{[16]} S.\t P.\t Novikov, 
{\it The algebraic construction and properties of
hermitian analogues of $K$-theory for rings with involution, from the point 
of view of the hamiltonian formalism. Some applications to differential 
topology and the theory of characteristic classes},
Izv. Akad. Nauk SSSR, ser. mat. 34, 253--288, 478--500 (1970)
\item{[17]} F.\t Quinn, 
{\it A geometric formulation of surgery}, in
{\it Topology of manifolds}, Proceedings 1969 Georgia Topology
Conference, Markham Press, 500--511 (1970)
\item{[18]} A.\t A.\t Ranicki, {\it Algebraic $L$--theory}, 
Proc. L.M.S. (3) 27, I. 101--125, II. 126--158 (1973)
\item{[19]} \bysame, {\it The algebraic theory of surgery}, 
Proc. L.M.S. (3) 40, I. 87--192, II. 193--287 (1980)
\item{[20]} \bysame, {\it The total surgery obstruction},
Proc. 1978 Arhus Topology Conference, 
Lecture Notes in Mathematics 763, 275--316, Springer (1979)
\item{[21]} \bysame, {\it Exact sequences in the algebraic theory of
surgery}, Mathematical Notes 26, Princeton (1981)
\item{[22]} \bysame, {\it Algebraic $L$-theory and topological manifolds},
Cambridge Tracts in Mathematics 102, CUP (1992)
\item{[23]} \bysame (ed.), {\it The Hauptvermutung Book},
Papers in the topology of manifolds by Casson, Sullivan, Armstrong, Rourke,
Cooke and Ranicki, K-Monographs in Mathematics 1, Kluwer (1996)
\item{[24]} \bysame, {\it High-dimensional knot theory}, 
Springer Mathematical Monograph, Springer (1998)
\item{[25]} C.\t W.\t Stark, {\it Surgery theory and infinite fundamental groups},
in {\it Surveys on surgery theory, Volume 1}, Ann. of Maths. Studies 145,
275--305, Princeton (2000)
\item{[26]} C.\t T.\t C.\t Wall,
{\it Classification of $(n-1)$-connected $2n$-manifolds},
Ann. of Maths. 75, 163--189 (1962)
\item{[27]} \bysame,
{\it Poincar\'e complexes}, Ann. of Maths. 86, 213--245 (1967)
\item{[28]} \bysame, 
{\it Surgery on compact manifolds}, Academic Press (1970);
2nd Edition (ed. A.\t A.\t Ranicki), 
Mathematical Surveys and Monographs 69, A.M.S. (1999)
\par}

\noindent
Dept. of Mathematics and Statistics\hfil\break
University of Edinburgh \hfil\break
Edinburgh EH9 3JZ\hfil\break
Scotland, UK\hfil\break
~\hfil\break
e-mail\t : \tt aar@maths.ed.ac.uk

\bye